\documentclass[english,11pt]{amsart}

\usepackage{amssymb}
\usepackage{hyperref} 
\usepackage{mathrsfs}
\usepackage{amsmath}
\usepackage{amsthm}
\usepackage{enumerate}
\usepackage{dsfont}
\usepackage{color}
\usepackage[a4paper]{geometry}
\usepackage[utf8]{inputenc}
\usepackage{stmaryrd}
\usepackage{titlesec}
\usepackage{wasysym}
\usepackage{appendix}
\usepackage{bigints}

\let\oldtocsection=\tocsection
\let\oldtocsubsection=\tocsubsection
\let\oldtocsubsubsection=\tocsubsubsection

\renewcommand{\tocsection}[2]{\hspace{0em}\oldtocsection{#1}{#2}}
\renewcommand{\tocsubsection}[2]{\hspace{1em}\oldtocsubsection{#1}{#2}}
\renewcommand{\tocsubsubsection}[2]{\hspace{2em}\oldtocsubsubsection{#1}{#2}}

\titleformat{\section}[block]
  {\normalfont\sffamily\bfseries\large}{{\hspace{-0.4cm}}\thesection}{.5em}{\titlerule\\[.8ex]} 

\titleformat{\subsection}[wrap]{\normalfont\sffamily\bfseries}{\hspace{-0.7cm}\thesubsection}{.5em}{}
                  \titlespacing{\subsection}{20pc}{1.5ex plus .1ex minus .2ex}{1pc}

\titleformat{\subsubsection}[wrap]{\normalfont\sffamily\bfseries}{\hspace{-1cm}\thesubsubsection}{.5em}{}
                  \titlespacing{\subsubsection}{15pc}{1ex plus .1ex minus .2ex}{1pc}

\newtheoremstyle{mystyle}
{3pt}                           
{3pt}							
{\sffamily}					
{-1.3em}						
{\bfseries\sffamily}		
{.}									
{.5em}							
{}									

\theoremstyle{mystyle}

\newtheorem{theorem}{$\RHD\,\,$ Theorem}[section]
\newtheorem{proposition}[theorem]{$\RHD\,\,$ Proposition}
\newtheorem{coro}[theorem]{$\RHD\,\,$ Corollary}
\newtheorem{lemma}[theorem]{$\RHD\,\,$ Lemma}
\newtheorem{rem}[theorem]{{\hspace{0.37cm}} Remark}

\newtheorem{definition}[theorem]{$\RHD\,\,$ Definition}
            \newtheorem*{definition*}{$\RHD\,\,$ Definition}

\numberwithin{theorem}{section}
\numberwithin{equation}{section}                  
\numberwithin{subsection}{section}
\numberwithin{subsubsection}{subsection}

\newenvironment{Dem}{%
    \begin{list}{\hspace{0.5cm}{\sf \textbf{Proof --}}}{%
        \setlength{\topsep}{0pt}%
        \setlength{\leftmargin}{0pt}%
        \setlength{\rightmargin}{0pt}%
        \setlength{\listparindent}{0pt}%
        \setlength{\itemindent}{0pt}%
        \setlength{\parsep}{0pt}%
        \addtolength{\leftmargin}{20pt}%
        \addtolength{\rightmargin}{0pt}%
    } \item }{\hfill{\begin{flushright} $\rhd$ \end{flushright}}\end{list}\smallskip}

\renewcommand{\epsilon}{\varepsilon}

\newcommand\N{\mathbb{N}}
\newcommand\R{\mathbb{R}}
\newcommand\E{\mathbb{E}}

\newcommand\G{\mathcal{G}}

\newcommand\LL{\mathcal L}

\newcommand\SSS{{\mathcal S}_o'}
\newcommand{\fr}[2]{{\scriptstyle \frac{#1}{#2}}}

\newcommand{\RR}{\mathbb{R}}  
\newcommand{\EE}{\mathbb{E}}
\newcommand{\mcE}{\mathcal{E}}
\newcommand{\NN}{\mathbb{N}}
\newcommand{\mcS}{\mathcal{S}}   
\newcommand{\mcA}{\mathcal{A}}
\newcommand{\mcR}{\mathcal{R}}
\newcommand{\mcP}{\mathcal{P}}
\newcommand{\mcQ}{\mathcal{Q}}

\newcommand{\bbX}{\mathbb{X}}
\newcommand{\bfX}{{\bf X}}

\newcommand{\calM}{\ensuremath{\mathcal{M}}}
\newcommand{\calC}{\ensuremath{\mathcal{C}}}

\newcommand{\abs}[1]{\left|#1\right|}

\DeclareMathOperator{\osc}{Osc}

\def\Xint#1{\mathchoice
   {\XXint\displaystyle\textstyle{#1}}%
   {\XXint\textstyle\scriptstyle{#1}}%
   {\XXint\scriptstyle\scriptscriptstyle{#1}}%
   {\XXint\scriptscriptstyle\scriptscriptstyle{#1}}%
   \!\int}

\def\XXint#1#2#3{{\setbox0=\hbox{$#1{#2#3}{\int}$}
     \vcenter{\hbox{$#2#3$}}\kern-.5\wd0}}

\def\aver#1{\Xint-_{#1}}

\newcommand{\ssk}{\smallskip}

\author[I. Bailleul]{Ismael Bailleul}
\address{Institut de Recherche Mathematiques de Rennes, 263 Avenue du General Leclerc, 35042 Rennes, France}
\email{ismael.bailleul@univ-rennes1.fr}

\author[F. Bernicot]{Fr\'ed\'eric Bernicot}
\address{CNRS - Universit\'e de Nantes \\ Laboratoire de Math\'ematiques Jean Leray \\ 2, Rue de la Houssini\`ere 44322 Nantes Cedex 03, France}
\email{frederic.bernicot@univ-nantes.fr}
 
\author{Dorothee Frey}
\address{CNRS - Universit\'e Paris-Sud, Laboratoire de Math\'ematiques, UMR 8628, 91405  Orsay, France }
\email{dorothee.frey@univ-nantes.fr}

\date{\today}

\begin{document}

\date{\today}

\vspace*{3ex minus 1ex}
\begin{center}
\huge\sffamily{Heat semigroup and singular PDEs \vspace{0.5cm}}
\end{center}
\vskip 5ex minus 1ex

\begin{center}
{\sf I. BAILLEUL\footnote{{\sf I.B. was partly supported by the ANR project ''Retour Post-doctorant'', no. 11-PDOC-0025; I.B. also thanks the U.B.O. for their hospitality, part of this work was written there.}} and F. BERNICOT\footnote{{\sf F. Bernicot's research is partly supported by ANR projects AFoMEN no. 2011-JS01-001-01 and HAB no. ANR-12-BS01-0013.}}}
\end{center}
\vskip 3ex minus 1ex

\begin{center}
\large \sffamily{(with an Appendix by F. Bernicot \& D. Frey)}
\end{center}

\vspace{1cm}

\begin{center}
\begin{minipage}{0.8\textwidth}
\renewcommand\baselinestretch{0.7} \sffamily {\scriptsize {\sf \noindent \textbf{Abstract.}} We provide in this work a semigroup approach to the study of singular PDEs, in the line of the paracontrolled approach developed recently by Gubinelli, Imkeller and Perkowski. Starting from a heat semigroup, we develop a functional calculus and introduce a paraproduct based on the semigroup, for which commutator estimates and Schauder estimates are proved, together with their paracontrolled extensions. This machinery allows us to investigate singular PDEs in potentially unbounded Riemannian manifolds under mild geometric conditions. As an illustration, we study the generalized parabolic Anderson model equation and prove, under mild geometric conditions, its well-posed character in H\"olders spaces, in small time on a potentially unbounded 2-dimensional Riemannian manifold, for an equation driven by a weighted noise, and for all times for the linear parabolic Anderson model equation  in 2-dimensional unbounded manifolds. This machinery can be extended to an even more singular setting and deal with Sobolev scales of spaces rather than H\"older spaces.   }
\end{minipage}
\end{center}

\bigskip

\setcounter{tocdepth}{2}
\begin{quote}
\footnotesize  \textsf{\tableofcontents}
\end{quote}

\section[\hspace{0.6cm} {\sf Introduction}]{Introduction}
\label{SectionIntro}

\subsection[\hspace{-1cm} {\sf State of the art}]{State of the art}
\label{SubsectionStateOfArt}

Following the recent breakthrough of Hairer \cite{H} and Gubinelli, Imkeller, Perkowski \cite{GIP}, there has been recently a tremendous activity in the study of parabolic singular partial differential equations (PDEs), such as the KPZ equation
$$
\big(\partial_t - \partial_x^2\big) u = \big(\partial_x u\big)^2 + \xi,
$$
the stochastic quantization equation
$$
\big(\partial_t - \Delta\big) u = -u^3 + \xi,
$$ 
or the generalized Parabolic Anderson Model equation
$$
\big(\partial_t - \Delta\big) u = \textrm{F}(u)\xi
$$
in all of which $\xi$ stands for a space or space-time white noise. Each of these equations involves, under the form of a product, a term which does not make sense a priori, given the expected regularity of the solution in terms of the regularity of the noise $\xi$. Hairer's theory of regularity structures is built on the insights of earlier works \cite{HCPAM,HW,HKPZ} on $(1+1)$-dimensional space-time problems where he used the framework of rough paths theory, under the form of Gubinelli's controlled paths, to make sense of previously ill-posed singular PDEs and give a meaningful solution theory. Rough paths theory was used in this approach as a framework for studying the properties in the $1$-dimensional \textit{space variable} of potential solutions. However, the very notion of a rough path is intimately linked with the $1$-dimensional time axis that parametrizes paths.

\medskip

To by-pass this barrier, both the theory of regularity structures and the paracontrolled approach developed in \cite{GIP} take as a departure point the fact that, like in rough paths theory, to make sense of the equation, one needs to enrich the noise $\xi$ into a finite collection of objects/distributions, and that one should try and describe the potential solution of a singular PDE in terms of that enriched noise. The latter depends on the equation under study and plays in the theory of regularity structures the role plaid by polynomials in the usual $C^k$ world to give local descriptions of functions under the form of Taylor expansions at every space-time point. The description of a solution in the paracontrolled approach is of a different nature and rests on a global comparison with the solution to a linear equation, $\big(\partial_t - \Delta\big) u = \xi$, in the above examples, via the use of Bony's paraproduct. In both approaches, the use of an ansatz for the solution space allows for fixed point arguments to give a robust solution theory where the solution becomes a continuous function of all the parameters of the problem.

\medskip

So far, both theories have only been formulated and tested on some singular PDEs on the torus, to the exception of the works \cite{HL,HLR3} on the parabolic Anderson model equation in $\RR^2$ and $\RR^3$, and our follow up work \cite{BBF16}. We introduce in the present work a functional analytic setting in which we are able to extend the paracontrolled approach of \cite{GIP} to investigate singular PDEs of the form 
$$
\big(\partial_t + L\big) u = \textrm{F}(u,\xi)
$$
for a second order differential operator $L$, and a nonlinear term $\textrm{F}(u,\xi)$, on potentially unbounded, Riemannian or even sub-Riemannian, manifolds or graphs. (The change of sign $-$ to $+$ in the operator is irrelevant.) This is a priori far from obvious as the main analytic tools used in the paracontrolled approach in the torus involve technics from Fourier analysis that do not make sense on manifolds or graphs. We develop to that end a functional calculus adapted to the heat semigroup associated with the operator $\big(\partial_t + L\big)$, which we use to define a paraproduct enjoying the same regularity properties as its Euclidean analogue. Such paraproducts adapted to a semigroup, as well as a paralinearization theory, have already been studied in recent works \cite{B-T1,BS}. However, the irregular character of the noises $\xi$ involved in the above motivating equations requires us to improve the definition of such paraproducts so as to build a framework where to consider regularity with a negative exponent; such an extension will be provided here. Building on these tools, one can set up, as in \cite{GIP}, a framework where to investigate the well-posed character of a whole class of parabolic singular PDEs. It is especially nice that all the objects in our framework are defined uniquely in terms of semigroups, unlike the notions of H\"older spaces used in the theory of regularity structures that involve a metric structure unrelated to the equation under study. As a by-product, we are able to handle some general classes of operators $L$ whose treatment seem to be beyond the present-day scope of the theory of regularity structures, as illustrated in some examples given in section \ref{SubsectionHeatSemiGroup}.

\medskip

It is unclear presently how one can adapt the different notions and tools of the theory of regularity structures to extend them to a (Lipschitz) manifold or graph setting, or to other second order operators (other than the Laplace operator), or how to work with Sobolev spaces instead of H\"olders spaces. Apart from the very definition of a regularity structure on a manifold, the existence of the reconstruction operator in this setting seems in particular challenging, as its proof in $\RR^d$ involves some deep results on wavelets that were not proved so far to hold true on generic manifolds, not even on all open sets of $\RR^d$. Their extension to a non-smooth setting also seems higly non-trivial. So it comes as a good news that one can use some reasonably elaborate theory of semigroups to implement the alternative machinery of paracontrolled calculus in that setting; as described below, it also allows us to have much flexibility on the operator $L$ and also on the geometry of the ambiant space. Roughly speaking, we could say that the point of view of the theory of regularity structures relies on the metric and differential properties of the underlying space, while the present extension of the paracontrolled calculus corresponds to a functional point of view adapted to the operator $L$ involved in the parabolic singular PDEs. We link here these two sides of the medal by requiring from the heat semigroup $(e^{-tL})_{t>0}$ to have a kernel  together with its gradient, that satisfies pointwise Gaussian bounds; this describes in some sense the link between the functional calculus and the ambiant space, with its metric and its differential geometry.

\ssk

We explain in Appendix B how this approach can be used in the context of Sobolev spaces rather than H\"older spaces. The former setting is slightly more difficult to handle, from a technical point of view, since Sobolev spaces involve simultaneously all the frequencies, whereas for H\"older spaces we can work at a fixed frequency scale. We do not know how such extension could be implemented within the setting of the regularity structure.

\medskip

The first part of this work is devoted to a precise study of the so-called paracontrolled calculus in a very abstract setting, given by a doubling ambiant space, equipped with a self-adjoint operator $-L$ generating a semigroup with Gaussian bounds for its kernel and its gradient. A suitable definition of paraproducts is given, and the main rules of calculus for paracontrolled distributions are described. This general theory is all we need to study a number of parabolic singular PDEs on manifolds.

\bigskip

\subsection[\hspace{-1cm} {\sf Paracontrolled calculus}]{Paracontrolled calculus}
\label{SubsectionParacontrolledCalculus}

The mechanics of paracontrolled calculus \cite{GIP} is elementary and easy to use; we describe it here as the present work only extends its scope while keeping its structure untouched -- unlike our work \cite{BBF16}. We use a somewhat informal style in this section and take the stochastic PDE given by the parabolic Anderson model equation (PAM)
\begin{equation}
\label{EqPAMParaC}
(\partial_t+\Delta) u = u\,\xi
\end{equation}
as an illustration. The symbol  $\xi$ stands here for a $2$-dimensional spatial white noise, of H\"older regularity $-1^-$. At the beginning of this story is the fact that one expects a solution $u$ to that equation to be $(-1^-+2)$-parabolic H\"older regular, as a consequence of the regularizing properties of the heat semigroup, while this regularity is not sufficient to make sense of the product $u\,\xi$, as the sum of their parabolic H\"older regularity is not positive. The point, however, is that $u$ is not expected to be any kind of $1^-$-H\"older function, rather it is expected to behave, at small space-time scale, like the solution $X$ of the elementary well-posed equation
$$
(\partial_t+\Delta) X =\xi,
$$
with null initial condition. The paracontrolled approach to solving the $2$-dimensional (PAM) equation consists in the following three step process.

\medskip

{\bf   \begin{enumerate}
   \item[{\sf \textbf{(1)}}] \textsf{Set yourself a $\xi$-dependent ansatz for the solution space, made up of functions that behave like $X$ at small space-time scale, and equipped with a Banach space structure.}   \vspace{0.1cm}
   
   \item[{\sf \textbf{(2)}}] \textsf{Show that the product $u\,\xi$ is well-defined for any element $u$ of the ansatz space.  }   \vspace{0.1cm}

   \item[{\sf \textbf{(3)}}] \textsf{Solve the equation via a fixed point argument.  }
\end{enumerate}  }

\medskip

The subtle point here is that the data of the distribution $\xi$ itself is not sufficient to give sense to the product $u\,\xi$, and that we really need that $\xi$ be random to build on the probability space where it is defined another distribution together with which one can make sense of all the products $u\,\xi$, for any $u$ in the ansatz space. Once this enrichment of $\xi$ has been done by purely probabilistic means, the above three step process is run in a deterministic setting.

\ssk

From a technical point of view, a $1^-$-H\"older function $u$ will be said to behave like $X$ at small scale if it is of the form $u\simeq \Pi_v(X)$, with $v$ bounded, up to some term more regular than $X$; write $u = \Pi_v(X) + u^\sharp$, with a remainder $u^\sharp$ of H\"older regularity strictly greater than $1$. The bilinear operator $\Pi_\cdot(\cdot)$, which will be a generalized paraproduct, that appears here has good continuity properties on large classes of distributions and satisfies the identity
$$
ab = \Pi_a(b) + \Pi_b(a) +\Pi(a,b)
$$
for any bounded functions $a,b$, for a continuous operator $\Pi(\cdot,\cdot)$ on $L^\infty\times L^\infty$ that happens to extend continuously to pairs of H\"older regular distributions for which the sum of their regularity is positive. In the torus, the Littlewood-Paley decomposition of $a$ and $b$ as an infinite sum of smooth functions whose Fourier transforms have supports in dyadic annuli can be used to define $\Pi_{\cdot}(\cdot)$ and $\Pi(\cdot,\cdot)$, by writing
$$
ab = \sum_{i,j} a_ib_j = \sum_{i<j-1}a_ib_j + \sum_{j<i-1}b_ja_i + \sum_{|i-j|\leq 1} a_ib_j.
$$
This definition justifies that we call $\Pi(\cdot,\cdot)$ the diagonal operator. The following formal analogy with the rules of stochastic calculus will enlighten the core technical tool of paracontrolled calculus described in a second. Recall that if $M$ and $N$ are two continuous martingales one has
$$
d(MN) = MdN + NdM + d\langle M,N\rangle.
$$
The above space of functions $u = \Pi_v(X)+u^\sharp,$ can be turned into a Banach space. Once this ansatz for the solution space has been chosen, remark that the product $u\,\xi$ can formally be written as
\begin{equation*}
\begin{split}
u\,\xi &= \Pi_u(\xi) + \Pi_\xi(u) + \Pi(u,\xi)  \\
		  &= \Pi_u(\xi) + \Pi_\xi(u) + \Pi\big(\Pi_v(X),\xi\big) + \Pi(u^\sharp,\xi).
\end{split}
\end{equation*}
Since $u^\sharp$ has H\"older regularity strictly bigger than $1$ and $\xi$ is $-1^-$-regular, the sum of their regularity indices is positive, and the term $\Pi(u^\sharp,\xi)$ is perfectly well-defined. This lives us with $\Pi\big(\Pi_v(X),\xi\big)$ as the only undefined term. The following fact is the workhorse of paracontrolled calculus. \textit{The trilinear map 
$$
C(a,b,c) := \Pi\big(\Pi_a(b),c\big) - a\,\Pi(b,c)
$$ 
happens to depend continuously on $a,b$ and $c$ provided they are H\"older distributions, with the sum of their H\"older exponents positive.} Note the paralell between the continuity of this 'commutator' and the rule for stochastic differentials, for which, given another continuous martingale $P$, we have
$$
d\Big\langle \int_0^\cdot M dN \,, P \Big\rangle = M\,d\langle N,P\rangle.
$$
The formal product $u\,\xi$ can thus be written as a sum of well-defined terms plus the formal product $v\,\Pi(X,\xi)$, with a diagonal term $\Pi(X,\xi)$ still undefined on a purely analytic basis. This is where probability comes into play. If one regularizes $\xi$ into $\xi^\epsilon$, with $X^\epsilon$ defined accordingly, one can prove that there exists a function/constant $C^\epsilon$ such that the renormalized quantity $\xi^{(2),\epsilon} := \Pi\big(X^\epsilon,\xi^\epsilon\big) - C^\epsilon$ converges in probability to some limit distribution $\xi^{(2)}$ of H\"older regularity $0^- = 1^- + (-1^-)$; this is enough to make sense of the product $v\,\xi^{(2)}$ on an analytical basis; but replacing $\Pi\big(X^\epsilon,\xi^\epsilon\big)$ by $\Pi\big(X^\epsilon,\xi^\epsilon\big)-C^\epsilon$ in the decomposition of $u\,\xi^\epsilon$ amounts to looking at the product $u\,\big(\xi^\epsilon-C^\epsilon\big)$. The enhancement $\widehat{\xi} := \big(\xi,\xi^{(2)}\big)$ of $\xi$ is called a \textbf{rough, or enhanced, distribution}, and one can use it to define the product $u\,\xi$ from the above formulas. At that point, it does not come as a surprise that one can then set (PAM) equation as a fixed point problem in the ansatz space, and that the unique solution to the problem (as it happens to be) is the limit of the solutions to the elementary problems
$$
(\partial_t + \Delta)u^\epsilon = u^\epsilon\big(\xi^\epsilon-C^\epsilon\big).
$$
The (PAM) equation \eqref{EqPAMParaC} is said to have been \textit{renormalized}. 

\medskip

More complicated problems are treated along these lines of reasoning in \cite{GIP, CC15, ZZ1, ZZ2}, to cite but a few. All use the above paraproduct machinery in the setting of the torus, where it is defined via Paley-Littlewood decomposition, such as described in this section. We introduce in this work a far more flexible paraproduct, defined intrinsically in terms of the semigroup associated with the operator $L$ that plays the role of $\Delta$ in the equation, in a very general geometrical setting. This offers the possibility to investigate stochastic PDEs in a manifold setting, which is our primary motivation. We gain much flexibility along the way, in terms of operators that can be used in place of $\Delta$, and even in Euclidean domains, the scope of the present work seems to be beyond the present day knowledge provided by the theory of regularity structures.

\bigskip

\subsection[\hspace{-1cm} {\sf A generalized parabolic Anderson model}]{A generalized parabolic Anderson model}
\label{SubsectionResultgPAM}

As an illustration of our machinery, we study the stochastic PDE given by the generalized parabolic Anderson model equation (gPAM)
$$
\partial_t u + Lu = \textrm{F}(u)\,\xi,\qquad u(0)=u_0,
$$
on some possibly unbounded $2$-dimensional Riemannian manifold $M$ satisfying some mild geometric conditions. One can take as operator $L$ the Laplace-Beltrami operator or some sub-elliptic diffusion operator; see section \ref{SubsectionHeatSemiGroup} for examples.  The nonlinearity F is $\calC^3_b$, and $\xi$ stands here for a weighted Gaussian noise with weight in $L^2\cap L^\infty$ -- see the definition in section \ref{SubsectionRenormalization}. The deterministic counterpart of the (gPAM) equation can be set once we are given a rough distribution $\widehat{\zeta} = \big(\zeta,\zeta^{(2)}\big)$; we write informally $u\,\widehat{\zeta}$ for the product operation in the ansatz space -- a slightly different and more precise notation will be adopted later on. The following results involve some parabolic H\"older spaces $\calC^\gamma$, with negative exponents $\gamma$, that are defined in section \ref{SubsectionHolder} in terms only of the semigroup $P_t$ generated by $-L$. We refer to section \ref{sec:pam} for a full statement and proof of these results. The geometrical assumptions stated here are introduced and explained in the next sections. The next two statements hold true provided the heat kernel of the semigroup generated by $-L$, together with its gradient, satisfy some Gaussian bounds given in the begining of section \ref{SubsectionHeatSemiGroup}; conditions on $L$ ensuring that these bounds hold are also given there. The letter $\calC^\alpha$ stands for a spatial (that is time-independent) H\"older space.

\medskip

\begin{theorem}    {\sf
\label{thm:TheoreticPam}
Assume that the measured metric manifold $(M,d,\mu)$ is equipped with a volume doubling measure, and that the heat semigroup generated by $-L$ has a kernel that satisfies some Gaussian bounds (UE), together with its gradient (Lip). Let $\alpha\in\big(\frac{2}{3},1\big)$, an initial data $u_0\in \calC^{2\alpha}$, a nonlinearity $\textrm{F}\in C^3_b$, and a positive time horizon $T$. Let $\widehat\zeta = \big(\zeta,\zeta^{(2)}\big)$ be a rough distribution, with $\zeta\in\calC^{\alpha-2}$ and $\zeta^{(2)} \in C_T\calC^{2\alpha-2}$.  \vspace{0.1cm}
\begin{enumerate}
   \item {\bf Local well-posedness for (gPAM).} For a small enough time horizon $T$, the generalized PAM equation
\begin{equation}
\label{EqgPAM}
\partial_t u + Lu = \textrm{F}(u)\cdot\widehat\zeta,\qquad u(0)=u_0 
\end{equation}
has a unique solution.  \vspace{0.1cm}

   \item {\bf Global well-posedness for (PAM).}  Under the assumption that the rough distribution takes values in some space of weighted distributions, the PAM equation
$$ 
\partial_t u + Lu = u\cdot\widehat\zeta,\qquad u(0)=u_0 
$$
has a unique global in time solution in some function space.  
\end{enumerate}   }
\end{theorem}

\medskip

The implementation of this result in the case where $\zeta = \xi$ is a random Gaussian spatial noise takes the following form, for a precise version of which we refer to theorem \ref{thm:Renormalization}; it holds in the same geometrical setting as the above result. The property of the measure $\mu$ put forward in the statement is called \textbf{Ahlfors regularity}; write $B_r(x)$ for a metric ball of center $x$ and radius $r$.

\medskip

\begin{theorem}   \label{thm:pam}   {\sf 
Suppose, in addition to the Gaussian bounds (UE) and (Lip) satisfied by the heat kernel and its derivative, that the reference measure $\mu$ on the manifold $M$ is doubling and satisfies the uniform lower bound $\mu\big(B_r(x)\big)\geq c_1 r^\nu$, for all $x\in M$, for some positive constant $c _1$ and the homogeneous dimension $\nu$. Let $\xi$ stand for a time-independent weighted noise in space, and set $\xi^\epsilon := P_\epsilon\xi$, and $X^\epsilon(t) = \int_0^t P_{t-s}\big(\xi^\epsilon\big)\,ds$.  \vspace{0.1cm}

\begin{enumerate}
   \item There exists a time-independent function $C^\epsilon (\cdot):= \E\Big[\Pi\big(L^{-1}\xi^\epsilon,\xi^\epsilon\big)(\cdot)\Big]$ on $M$ such that the pair $\big(\xi^\epsilon , X^\epsilon- C^\epsilon\big)$ converges in probability to a random rough distribution $\widehat\xi$.   \vspace{0.1cm}

   \item If $u^\epsilon$ stands for the solution of the renormalized equation 
   \begin{equation*} 
   \partial_t u^\epsilon + L u^\epsilon = \textrm{F}\big(u^\epsilon\big)\,\xi^\epsilon - C^\epsilon\,\textrm{F}'\big(u^\epsilon\big)\,\textrm{F}(u^\epsilon) ,\qquad u^\epsilon(0)=u_0  
   \end{equation*}
it converges in probability to the solution $u$ of equation \eqref{EqgPAM} driven by $\widehat\xi$.
\end{enumerate}   }
\end{theorem}

\bigskip

Note that one cannot expect the renormalizing function $C^\epsilon$ to be constant unless the manifold $M$ is homogeneous and the operator $L$ commutes with the group action -- which holds in the torus when working with the Laplacian. Note also that we do \textit{not} assume $M$ to be bounded. Working with a weighted noise rather than with white noise allows us to by-pass the somewhat heavy use of weighted H\"older spaces, such as done in \cite{HL,HLR3} and \cite{BBF16}; the latter work deals, among other things, with paracontrolled calculus in weighted H\"older spaces.

\bigskip

We have organized our work as follows. Section \ref{sec:preli} presents the functional setting in which our theory is set. The main geometrical assumptions on the geometric background are given in section \ref{SubsectionHeatSemiGroup}, where examples are given; these assumptions involve the properties of the heat kernel of the semigroup $\big(e^{-tL}\big)_{t\geq 0}$ generated by $L$. A family of operators is introduced in section \ref{SubsectionGamma}, which will play in the sequel the role played by Fourier projectors in the classical Littlewood-Paley theory. We introduce in section \ref{SubsectionHolder} a scale of H\"older spaces, defined uniquely in terms of the semigroup $\big(e^{-tL}\big)_{t\geq 0}$. A paraproduct is introduced in section \ref{SubsectionParaproducts} and is shown in section \ref{SubsectionParapdctEstimates} to enjoy the same continuity properties as its Euclidean analogue. A crucial commutator estimate between paraproduct and resonant terms is proved in section \ref{SubsectionCommutatorEstimates}, together with some paralinearization and composition estimates in section \ref{SubsectionParalinearization}. Following \cite{GIP}, we then introduce in section \ref{SubsectionParaControlledDistributions} what plays the role in our setting of paracontrolled distributions, and prove some fundamental Schauder estimates in section \ref{SubsectionSchauder}. Sections 2 to 4 give us all the material needed to investigate a number of singular PDEs on manifolds from the point of view of paracontrolled distributions. Section \ref{sec:pam} is dedicated to the proof of theorems \ref{thm:TheoreticPam} and \ref{thm:pam}.

\medskip

We end this work by Appendix B, jointly written with Dorothee Frey, in which we explain how we can weaken our assumption of Lipschitz regularity of the heat kernel (Lip), which we make in the main body of this work, in terms of more geometrical properties. We also show that one can prove results in Sobolev spaces similar to those proved in H\"older spaces in the main body of that work.

\bigskip
\bigskip

We collect here a number of \textbf{notations} that will be used throughout that work. \vspace{0.2cm}

\begin{itemize}
   \item For a ball $B$ of radius $r$ and a real $\lambda>0$, denote by $\lambda B$   the ball concentric  with $B$ and with radius $\lambda r$. We shall use $u\lesssim v$ to say that there exists a constant $C$ (independent of the important parameters) such that $u\leq Cv$ and $u\simeq v$ to say that $u\lesssim v$ and $v\lesssim u$. We also adopt the non-conventional notation $\gamma_a$ for the classical gamma function, defined for $a>0$ by the formula
$$ 
\gamma_a := \int_0^\infty x^a e^{-x} \, \frac{dx}{x};
$$
the capital letter $\Gamma$ will be used to denote the carr\'e du champ operator of some other operator.  \vspace{0.1cm}

 \item For $p\in[1,\infty]$ and every $f\in L^p$, the $L^p$-norm, with respect to the measure $\mu$, is denoted by $\|f\|_p$. For $p,q\in[1,\infty]$, we write $\|T\|_{p\to q}$ for the operator norm of an operator $T$ from $L^p$ to $L^q$.   \vspace{0.1cm}

   \item For an integer $k\geq 0$, we write $C^k_b$ for the set of functions continuously differentiable $k$-times $f:\R \rightarrow \R$, equipped with the norm
$$ 
\|f\|_{C^k_b} := \|f\|_{\infty} + \sup_{1\leq i\leq k} \Big\| f^{(i)} \Big\|_\infty.
$$
\end{itemize}

\bigskip

\section[\hspace{0.6cm} {\sf Functional calculus adapted to the heat semigroup}]{Functional calculus adapted to the heat semigroup} 
\label{sec:preli}

As announced in the introduction, this section is dedicated to describing the functional framework where we shall set our study. Section \ref{SubsectionHeatSemiGroup} sets the geometrical framework needed for what we want to do, in terms of a semigroup. We introduce in section \ref{SubsectionGamma} some operators that will play the role of 'localizers' in frequency space. These operators are used in section \ref{SubsectionHolder} to define a scale of H\"older spaces which will be instrumental in the sequel.

\ssk

\subsection[\hspace{-1cm} {\sf Heat semigroup on a doubling space}]{Heat semigroup on a doubling space}
\label{SubsectionHeatSemiGroup}

Let denote by $(M,d)$ be a locally compact separable metrisable space, equipped with a Radon measure $\mu$, strictly positive on any non-empty open set. Given a ball $B(x,r)$ of center $x$ and radius $r$, the notation $V(x,r)$ will stand in the sequel for $\mu\big(B(x,r)\big)$. To make things concrete, the space $(M,d)$ will mainly be for us  smooth Riemannian manifold or a (possibly infinite) metric graph. We shall assume  that the metric measure space  $(M,d,\mu)$ satisfies the following \textbf{volume doubling property}
\begin{equation} \label{d} \tag{VD}
V(x,2r)\lesssim  V(x,r),
\end{equation}
for all $x\in M$ and positive $r$, which can be stated equivalently under the form 
\begin{equation} \label{dnu}
V(x,r)\lesssim \left(\frac{r}{s}\right)^{\nu} V(x,s),
\end{equation}
for some positive scaling factor $\nu$, for all $x\in M$, and all $0<s\leq r$; it implies the inequality
\begin{equation*}
V(x,r)\lesssim \left(\frac{d(x,y)+r}{s}\right)^{\nu} V(y,s),
\end{equation*}
for any two points $x,y$ in $M$ and $0<s\leq r$. (Another easy consequence of the volume doubling property is that  balls with a non-empty intersection and comparable radii have comparable measures.)

\medskip

Let also be given a non-negative self-adjoint operator $L$ on $L^2(M,\mu)$ with dense domain ${\mathcal D}_2(L) \subset L^2(M,\mu)$. Denote by $\mathcal{E}$ its associated quadratic form, defined by the formula
$$
{\mathcal E}(f,g) := \int_Mf Lg\,d\mu,
$$  
on a domain $\mathcal{F}$ which contains ${\mathcal D}_2(L)$. We shall assume that the \textit{Dirichlet form} $\mathcal{E}$ is  \textit{strongly local} and \textit{regular}; we refer the reader to the books \cite{FOT, GSC} of Fukushima \& co. and Gyrya--Saloff-Coste for precise definitions and background on Dirichlet forms. (The reader unfamiliar with this setting may think of the Laplace operator in a compact Riemannian manifold.) These two properties will be obviously satisfied in the examples we shall work with. It follows from these conditions that the operator $L$ generates a strongly continuous semigroup  $\big(e^{-tL}\big)_{t>0}$ of contractions on $L^2(M,\mu)$ which is conservative, in the sense that $e^{-tL}{\bf 1} = {\bf 1}$, for all $t\geq 0$; see e.g. Subsection 2.2.7 in the book \cite{GSC}. We shall also assume that the semigroup $\big(e^{-tL}\big)_{t>0}$ has a  kernel, given for all positive times $t$ by a non-negative measurable real-valued function $p_t$ on $M\times M$, such that 
$$
\big(e^{-tL}f\big)(x) = \int_M p_t(x,y)f(y)\,d\mu(y), 
$$ 
for $\mu$-almost all $x$ in $M$, and every $f\in {\mathcal D}_2(L)$. The kernel $p_t$ is called the \textbf{heat kernel} associated with $L$. We assume that it satisfies for all $0<t\leq 1$ and $\mu$-almost all $x,y$, the following typical upper estimates
\begin{equation*}
 p_{t}(x,y)\lesssim \frac{1}{\sqrt{V\big(x,\sqrt{t}\big)V\big(y,\sqrt{t}\big)}}.
\end{equation*}
Under the volume doubling condition \eqref{d}, the previous estimate self-improves into a Gaussian upper estimate (UE) for the heat kernel and its time derivatives
\begin{equation}\tag{UE}
\label{UE} 
\Big| \partial_t^a p_{t}(x,y)\Big| \lesssim \frac{t^{-a}}{\sqrt{V\big(x,\sqrt{t}\big)V\big(y,\sqrt{t}\big)}}\,\exp\left(-c\;\frac{d(x,y)^2}{t}\right).
\end{equation}
that holds for a fixed positive constant $c$, for all integers $a$, all times $0<t\leq 1$, and $\mu$-almost every $x,y\in M$; see for instance the article \cite[Theorem 1.1]{Gr1} for the Riemannian case, and the work \cite[Section 4.2]{CS} for a metric measure space setting. We also assume that the heat kernel satisfies the following Lipschitz condition
\begin{equation}  \label{Lipschitz}  \tag{$\textrm{Lip}$}
\Big| p_{t}(x,y)-p_t(z,y) \Big|\lesssim \left(\frac{d(x,z)}{\sqrt{t}}\right) \frac{1}{\sqrt{V\big(x,\sqrt{t}\big)V\big(y,\sqrt{t}\big)}}\,\exp
\left(-c \;\frac{d(x,y)^2}{t}\right). 
\end{equation}
Let insist here that \textit{inequalities \eqref{UE} and \eqref{Lipschitz} are assumed to hold only for $0<t\leq 1$}, rather than for all positive times. It follows classically from the Gaussian estimates \eqref{UE} and the volume doubling property that the heat semigroup $\big(e^{-tL}\big)_{t>0}$ is uniformly  bounded  on $L^p(M,\mu)$ for every $p\in[1,\infty]$, and strongly continuous for $p\in[1,\infty)$. Last, note that $\big(e^{-tL}\big)_{0<t\leq 1}$ is, under these conditions, bounded analytic on $L^p(M,\mu)$, for every $1<p<+\infty$, which means in particular that the time-derivatives $\big((tL)^n e^{-tL}\big)_{0<t\leq 1}$ are bounded on $L^p(M,\mu)$ uniformly in $0<t\leq 1$, for every integer $n\geq 0$; see \cite{topics}.

\medskip

A comment is in order here, about our two assumptions \eqref{UE} and \eqref{Lipschitz}. In the theory of regularity structures or Euclidean theory of  paracontrolled calculus, regularity at any order may be considered because of the implicit use of the very nice differential geometry of Euclidean space, or the torus. In our current and far more general framework, since we only have a pointwise assumption on the heat kernel and its gradient, it is natural to expect that one cannot quantify the regularity of some objects to an order greater than $1$. That is why in the different statements proved in the next sections some extra mild conditions on the regularity exponents will appear, as compared with their Euclidean analogue. Since we aim to work within the present optimal / minimal setting, these new limitations cannot be removed without additional assumptions, and we shall be restricted to study regularity properties at order at most $1$, including negative orders; this is not restrictive as far as applications are concerned in the present work.

\medskip

Here are four representative classes of \textbf{examples} of doubling metric measure spaces and Dirichlet forms satisfying the above conditions. This list of examples emphasizes that we have much flexibility in terms of the operator $L$ as well as in terms of the underlying space $(M,d,\mu)$.

\medskip

\begin{enumerate}
   \item {\sf \textbf{Markov chains.}} Let $X$ be a countable set equipped with a Markov chain, specified by a symmetric Markov kernel $ k : X \times X \rightarrow \RR_+$, and let $m$ be a non-negative function on $X$, used to define a measure ${\sf m}$ on $X$, with density $m$ with respect to the counting measure $\mu$. Denote by $\langle \cdot,\cdot\rangle_{\sf m}$ the scalar product on $\ell^2({\sf m})$. Consider also for integers $n \geq 1$  the iterated kernel $k^n$ defined recursively by $k^n(x,y) := \int k^{n-1}(x,z) k(z,y)\, \mu(dz)$. Denoting by $K$ the symmetric Markov operator with kernel $k$ -- with respect to $\mu$, the formula
\begin{equation*}
\begin{split}
\mcE(f,g) &= \frac{1}{2}\sum_{x,y\in X}k_{xy}\big(f_x-f_y\big)\big(g_x-g_y\big)  \\
                &= \sum_{x\in X}f_x\,\frac{1}{m_x}\Big(g_x - \sum_{y\in X} k_{xy}g_y\Big)m_x  \\
                &= \big\langle f,Lg \big\rangle_{\sf m}
\end{split}
\end{equation*}   
associated with the non-negative self-adjoint operator
$$
\big(Lg\big)(x) = \frac{1}{m_x}\Big(g_x - \sum_{y\in X} k_{xy}g_y\Big) = \frac{1}{m_x} \big(g_x - \big(Kg\big)_x\big),
$$
defines a (strongly local) regular Dirichlet form and allows us to generate the continuous heat semigroup $\big(e^{-tL}\big)_{t\geq 0}$. (The above sum in $x$ is implicitly restricted to those $x$ for which $m_x>0$, so there is no loss of generality in assuming that $m>0$. ) The map $k$ induces a distance $d$ on $X$ by letting be equal to 
$$
\min\big\{n\geq 1\,;\,\exists\,z_0,\dots,z_n, \textrm{with }z_0=x,\,z_n=y \textrm{ and } k\big(z_i,z_{i+1}\big)>0, \textrm{ for } i=0..n-1\big\},
$$
for $y\neq x$. Following Grigor'yan's result \cite{Gr2}, one can give growth conditions on the $\sf m$-volume of $d$-balls that ensure the conservative character of the semigroup generated by $L$ in $\ell^2(\sf m)$. Then it is classical that getting Gaussian upper estimates for the semigroup $\big(e^{-tL}\big)_{t\geq 0}$ is very closely related to getting discrete-time versions of Gaussian estimates for the iterated Markov chains $\big(K^n\big)_{n\geq 1}$, and similarly for the Lipschitz regularity of their kernels. Usually, given such a discrete framework, one prefers to work with the discrete-time Markov chains rather than the continuous heat semigroup. To obtain upper Gaussian estimates and a Lipschitz regularity for the iterated Markov chains on a graphs is the topic of a huge literature to which we refer the reader; see for instance works by Hebisch and Saloff-Coste \cite{HSc1} for discrete groups and by Ischiwata \cite{Is} for an extension to nilpotent covering graphs and more recently \cite{Is2} for a perturbation of these previous results. For example, the regular graphs ${\mathbb Z}^d$ and $({\mathbb Z}/{N {\mathbb Z}})^d$ have heat semigroups satisfying the Gaussian estimates \eqref{UE} and the Lipschitz property \eqref{Lipschitz}. Needless to say, for a (large) finite graph $(X,E)$, with edge set $E$, and $b_{xy} = 1$ if $(x,y)\in E$, and $m_x = \sum_{y\in X} b_{xy}$, the previous results hold with the graph distance in the role of $d$. \vspace{0.2cm}

   \item {\sf \textbf{Second order differential operators on Riemannian manifolds.}} Let $(M,d,\mu)$ be a doubling, possibly non-compact," complete Riemannian manifold with Ricci curvature bounded from below. Then the heat semigroup $\big(e^{-t\Delta}\big)_{t\geq 0}$ generated by the Riemannian Laplace operator satisfies both the upper Gaussian estimates \eqref{UE} and the Lipschitz regularity \eqref{Lipschitz} for small time $0<t\leq 1$, and for every time $t>0$ if the Ricci curvature is nonnegative; see \cite{Var1} and \cite{LY} for references. Particular examples are smooth compact Riemannian manifolds, or unbounded Riemannian manifolds with pinched negative Ricci curvature, such as hyperbolic spaces. 

\ssk   
   
Even on the Euclidean space $\RR^d$, we may consider a second order divergence form operator $L = -\textrm{div}(A\nabla)$ given by a map $A$ taking values in real symmetric matrices and satisfying the usual ellipticity/accretivity condition. Then if $A$ is H\"older continuous, it is known that $-L$ generates a self-adjoint semigroup satisfying the properties \eqref{UE} and \eqref{Lipschitz}; see \cite{AMT}. Similarly, consider an open bounded subset $\Omega\subset \RR^d$ (with Lipschitz boundary for example to ensure the Ahlfors regularity) and consider the self-adjoint Laplace operator $L$ associated with Dirichlet or Neumann boundary conditions. There is an extensive literature to describe assumptions on $\Omega$ such that \eqref{UE} and \eqref{Lipschitz} are satisfied. The present setting may well be beyond the present scope of regularity structures, for which the Green function of the operator needs to satisfy some regularity assumptions that were not proved to hold true under a sole H\"older continuity assumption for $A$, and whose formulation on a manifold is a real problem outside the realm of Lie groups or homogeneous spaces. On the other hand, the theory developed here works well in that relatively minimal setting.  \vspace{0.1cm}

The estimates \eqref{UE} and \eqref{Lipschitz} also hold when working on a convex or $C^2$-regular bounded subset of the Euclidean space, with $L$ given by Laplace operator with Neumann boundary conditions;   see \cite{WY}.  \vspace{0.2cm}
  
   \item {\sf \textbf{Sub-elliptic left invariant diffusions on groups.}} Let $G$ be a unimodular connected Lie group, endowed with its left-right Haar measure $\mu$. Consider a family ${\bf X}:=\{X_1,...,X_\ell\}$ of left-invariant vector fields on $G$ satisfying H\"ormander condition. They define a class of admissible paths $\gamma_\bullet$, characterized by the existence, for each of them, of measurable functions $a_1,...,a_k$ such that one has 
$$
\gamma'(t) = \sum_{i=1}^k a_i(t) X_i(\ell(t)).
$$
The length of such a curve is defined as
$$ 
\big|\gamma\big| := \frac{1}{2}\int_0^1 \left( \sum_{i=1}^\ell |a_i(t)|^2 \right)^\frac{1}{2}\, dt,
$$
and the (Carnot-Caratheodory) distance $d(x,y)$ between any two points $x,y$ of $G$ is defined as the infimum of the lengths of all admissible curves joining $x$ to $y$. We then consider the sublaplacian $\Delta$ defined by 
$$ 
\Delta := - \sum_{i=1}^k X_i^2.
$$
Then the operator $\Delta$ generates a heat semigroup satisfying both the upper Gaussian estimates \eqref{UE} and the Lipschitz regularity \eqref{Lipschitz} for small time $t\in(0,1]$; see for instance Chapter 8 in the book \cite{VSC}. If the group is nilpotent then it is also globally doubling \cite{G} and so the heat semigroup satisfies the Gaussian upper bound \eqref{UE} and enjoys the Lipschitz property \eqref{Lipschitz} for every $t>0$; see \cite{Var0,Sa}. Particular examples of such groups, are stratified Lie groups, and so Heisenberg groups. For such Heisenberg-type Lie groups, a kind of Fourier transform may be defined involving their irreducible unitary representations, which can be used to define an analog of the Euclidean paraproducts / paradifferential calculus, such as done is \cite{GS}. We shall see, as a by-product of the present work, that the structure of heat semigroup is sufficient to construct similar tools with greater scope. \vspace{0.2cm}
 
   \item The general case given by a subelliptic operator is more difficult. Let $(M,d,\mu)$  be a complete and smooth connected manifold endowed with a self-adjoint smooth locally subelliptic diffusion operator $L$ satisfying $L1 = 0$. Then Baudoin and Garofalo introduced in \cite{BG} a property, called ``a generalized curvature-dimension inequality", which has to be thought of as a lower bound on a sub-Riemannian generalization of the Ricci tensor. Under such a condition, the heat kernel generated by $L$ satisfies \eqref{UE} as well as \eqref{Lipschitz}; see \cite{Qian}. We refer the reader to \cite{BG} for some examples of such sub-elliptic settings and the fact that the heat kernel also satisfies in that case some Gaussian lower bound.
\end{enumerate}

\bigskip

Throughout that work, a point $o\in M$ will be fixed, which we shall use to define a class of test functions, together with its 'dual' class of distributions.

\medskip

\begin{definition*}  {\sf 
We define a Fr\'echet space of test functions setting
$$ 
{\mathcal S}_o := \Bigg\{ f\in \bigcap_{n\geq 0} {\mathcal D}_2\big(L^n\big)\,;\, \forall \, a_1,a_2\in \NN, \quad  \Big\|\big(1+d(o,\cdot)\big)^{a_1} L^{a_2}f\Big\|_{2} <\infty\Bigg\},
$$
with
$$
\| f \| := \sup_{a_1,a_2\in \NN}\; 1\wedge \Big\|\big(1+d(o,\cdot)\big)^{a_1} L^{a_2}f\Big\|_{2}.
$$
A \emph{\textbf{distribution}} is a continuous linear functional on ${\mathcal S}_o$; we write ${\mathcal S}_o'$ for the set of all distributions.  }
\end{definition*}
 
\medskip 
 
(We point out that the arbitrary choice of point $o\in M$ is only relevant in the case of a unbounded ambiant space $M$; even in that case, the set ${\mathcal S}_o$ does not depend on $o$, for $o$ ranging inside a bounded subset of $M$.) Every bounded function defines for instance an element of $\mathcal{S}'_o$. Examples of test functions are provided by the $p_t(x,\cdot)$, for every fixed $x\in M$ and $0<t\leq 1$. Indeed for integers $a_1, a_2$, the upper bound \eqref{UE} with the analyticity of the semigroup yield that $(tL)^{a_2}p_t(x,\cdot)$ satisfies the same upper Gaussian estimates than the heat kernel itself and so we deduce that
\begin{align*}
 \Big| \big(1+d(o,y)\big)^{a_1} \big(L^{a_2} p_t(x,\cdot)\big)(y)\Big| & \lesssim \frac{t^{-a_2}}{V(x, \sqrt{t})} (1+d(o,y))^{a_1} e^{-c\frac{d(x,y)^2}{t}} \\
  & \lesssim \frac{t^{-a_2}}{V(x, \sqrt{t})} (1+d(o,x))^{a_1} e^{-c'\frac{d(x,y)^2}{t}} 
\end{align*}
for some positive constants $c$ and $c'$. Note that the heat semigroup acts not only on functions, but also on distributions, by setting 
$$
\big\langle\big(e^{tL}\phi\big) , f\big\rangle := \big\langle\phi , \big(e^{tL}f\big)\big\rangle
$$
for $\phi\in \mcS_o'$ and $f\in\mcS_o$. We refer the reader to \cite{BDY} and \cite{LYY} for more details on the extension of the semigroup to distributions.
 
\ssk

For a linear operator $T$ acting from ${\mathcal S}_o$ to ${\mathcal S}_o'$, it will be useful below, to denote by $K_T$ its Schwarz kernel, characterized by the   identity 
$$ 
\langle T(f),g\rangle = \int K_T(x,y) f(y) g(x) \, \mu(dy) \mu(dx),
$$
giving an integral representation for every $f,g\in\mcS_o$.

\medskip

\subsection[\hspace{-1cm} {\sf Time derivatives and Carr\'e du champ of the semigroup}]{Time derivatives and Carr\'e du champ of the semigroup}
\label{SubsectionGamma}

Let us introduce here a family of operators that will play the role in our setting of the Fourier multipliers used in the classical Littlewood-Paley theory, that localize a function in frequency space. These will be the building blocks used to define a convenient paraproduct for our needs, such as done below in section \ref{SubsectionParaproducts}.

\medskip

\begin{definition*} \label{def:Qt-Pt}  {\sf
Given a fixed positive integer $a$, set
\begin{equation} \label{def:Qt}
Q_t^{(a)} := (tL)^a e^{-tL}
\end{equation}
and 
\begin{equation}  \label{def:Pt}
P_t^{(a)} := \phi_a(tL), \qquad \text{ where } \qquad \phi_a(x) := \frac{1}{\gamma_a} \int_x^\infty s^a e^{-s} \,\frac{ds}{s}, \;\;\; x\geq 0,
\end{equation}
for every $t>0$.   }
\end{definition*}

\medskip

So we have for instance $P_t^{(1)} = e^{-tL}$, and $Q_t^{(1)} = tLe^{-tL}$. The two families of operators $\big(P_t^{(a)}\big)_{t>0}$ and $\big(Q_t^{(a)}\big)_{t>0}$ are defined so as to have the relation
\begin{equation}
\label{EqRelationPtQt}
t\partial_t P_t^{(a)} =  tL\phi_a'(tL)= - \gamma_a^{-1}\,Q_t^{(a)},
\end{equation}
so $Q_t^{(a)}= (-1)^a t^a \partial_t^a e^{-tL}$, and $P_t^{(a)} = p_a(tL)e^{-tL}$, for some polynomial $p_a$ of degree $a-1$, with $p_a(0)=1$. The analyticity of the semigroup provides a direct control of the operators $P^{(a)}_t$ and $Q^{(a)}_t$.

\medskip

\begin{proposition}   \label{PropGaussianPtQt}    {\sf
For any integer $a\geq 0$, the operators $P_t^{(a)}$ and $Q_t^{(a)}$ have kernels satisfying the Gaussian estimate \eqref{UE}, and the Lipschitz regularity property \eqref{Lipschitz}; as a consequence, they are bounded in every $L^p$ spaces for $p\in[1,\infty]$, uniformly with respect to $t\in (0,1]$.  }
\end{proposition}

\medskip

Following the above interpretation of the operators $Q^{(a)}$ and $P^{(a)}$, the following Calder\'on reproducing formula provides a decomposition of a function $f$ in $L^p(M,\mu)$ into a low frequency part and a high frequency part very similar to the Littlewood-Paley decomposition of a distribution in terms of frequencies; see e.g. \cite{BCD}.

\medskip

\begin{proposition}[Calder\'on reproducing formula]  \label{prop:reproducing}   {\sf
Given $p \in (1,{+\infty})$ and $f\in L^p(M,\mu)$, we have 
\begin{align*}
\lim_{t\to {0^+}} P_t^{(a)}f & =  f \quad \textrm{ in $L^p(M,\mu)$ }
\end{align*}
for every positive integer $a$, and so
\begin{equation}\label{calde}
 f = \gamma_a^{-1} \int_0^{1} Q_t^{(a)}f \, \frac{dt}{t} +  P^{(a)}_{1}(f).  
\end{equation}   }
\end{proposition}

\medskip

\begin{Dem}
One knows from theorem 3.1 in \cite{DR}, that the operator $L$ has a bounded $H^\infty$ functional calculus in $L^p(M,\mu)$ under the volume doubling condition on $(M,d,\mu)$, and the assumption that the heat kernel satisfies the upper estimate \eqref{UE}. Since this implies in particular sectoriality of  $L$ in $L^p(M,\mu)$, Theorem 3.8 in \cite{CDMY} yields the decomposition of $L^p(M,\mu)$ into nullspace and range of $L$. Using this decomposition, the Convergence Lemma implies for every $f\in L^p(M,\mu)$
\begin{align*}
f  & =  \lim_{t\to 0} P_t^{(a)}f =  -  \int_0^1 \partial_t P_t^{(a)}f \, dt+  P^{(a)}_1(f) \\
 & = \gamma_a^{-1} \int_0^{1} Q_t^{(a)}f \, \frac{dt}{t} +  P^{(a)}_1(f),
\end{align*}
where the limit is taken in $L^p(M,\mu)$ and where we have used identity \eqref{EqRelationPtQt}; see e.g. \cite[Theorem D]{ADM} or \cite[Lemma 9.13]{KW}.
\end{Dem}

\bigskip

We shall also make an extensive use in the sequel of the square-root of $L$, given by its \textbf{carr\'e du champ} operator $\Gamma$, defined for all $(f,g)\in {\mathcal D}_2(L) \times {\mathcal D}_2(L)$ as a bilinear operator satisfying the identity
$$ 
\mathcal{E}(f,g) := \int_M f L(g) \, d\mu = \int_M g L(f) \, d\mu = \int_M  \Gamma(f,g) \, d\mu.
$$
It is also given by the explicit formula 
$$
\Gamma(f,g) = -\frac{1}{2} \Big(L(fg)-fL(g) - gL(f)\Big);
$$
we shall write ${\mathcal D}_2(\Gamma)\subset L^2$ for its domain, which contains ${\mathcal D}_2(L)$. 

As a shorthand, we write $\Gamma(f)$ for $\Gamma(f,f)^\frac{1}{2}$ in the sequel, which can be thought as the length of the intrinsic gradient of $f$. It follows from the conservative property of $L$ and its non-negative character, that the bilinear map $\Gamma$ is positive and satisfies the identity
$$ 
\big\| \Gamma(f) \big\|_{L^2}^2 = \int_M \Gamma(f,f) \, d\mu = \int_M f L(f) \, d\mu = \mathcal{E}(f,f).
$$
From its positive property, a generalized Cauchy-Schwarz inequality yields that for every $f,g\in {\mathcal D}_2(L)$ then
\begin{equation}
\left|\Gamma(f,g)\right|^2 \leq \Gamma(f,f) \Gamma(g,g)=\Gamma(f) \Gamma(g).
\label{eq:CS}
\end{equation}

According to the Beurling-Deny-Le Jan formula, the carr\'e du champ satisfies a \textit{Leibniz rule}
\begin{equation}  \label{eq:leibniz}
 \Gamma(fg,h) = f\,\Gamma(g,h) + g\,\Gamma(f,h),   
\end{equation}
for all $f,g,h\in\mathcal{D}_2(\Gamma)$, and a \textit{chain rule} 
\begin{equation}  \label{eq:chain-rule}
L\big(F(f)\big) = F'(f)\,L(f) + F''(f)\,\Gamma(f,f).  
\end{equation}
for every function $F\in C^2_b(\R)$ and every $f\in {\mathcal D}_2(L)$; the function $F(f)$ is automatically in ${\mathcal D}_2(L)$ -- see e.g. \cite[Section 3.2]{FOT} and \cite[Appendix]{ST1} for these points.

\medskip

The following pointwise and $L^p$-estimate for the intrinsic gradient of the semigroup will be used several times in a crucial way; its proof is given in Appendix \ref{Appendix}. It says that the carr\'e du champ of the semigroup satisfies also some Gaussian pointwise estimates, as given by the following claim.

\medskip

\begin{proposition}  \label{prop:gradientsemigroup}   {\sf
The following inequality holds
\begin{equation}  \label{eq:gradientpt-0}
\Big|\big(\sqrt{t}\,\Gamma)\big(e^{-tL} f\big)(x_0)\Big| \lesssim \int_M \frac{1}{\sqrt{V\big(x_0,\sqrt{t}\big)V\big(y,\sqrt{t}\big)}}\,\exp\left(-c\; \frac{d(x_0,y)^2}{t}\right) \big|f(y)\big| \, d\mu(y),
\end{equation}
for every $t>0$, every function $f\in L^2$, and almost every $x_0\in M$. Consequently, we have 
$$ 
\sup_{t>0} \, \Big\| \big(\sqrt{t} \Gamma\big)\big(e^{-tL}\cdot\big) \Big\|_{p \to p} < \infty,
$$
for every $p\in[1,\infty]$. We may replace the semigroup $e^{-tL}$ in the above equations by any of the operators $P_t^{(a)}$, for any $a\geq 0$.   }
\end{proposition}

\bigskip

\subsection[\hspace{-1cm} {\sf H\"older and Besov spaces through the heat semigroup}]{H\"older and Besov spaces through the heat semigroup}
\label{SubsectionHolder}

Let us recall as a start that given a parameter $\sigma\in (0,1]$, a bounded function $f\in L^\infty$ is said to belong to the H\"older space $\Lambda^\sigma$ if
$$
\|f\|_{\Lambda^ \sigma} := \|f\|_\infty + \sup_{0<d(x,y)\leq 1} \frac{|f(x)-f(y)|}{d(x,y)^\sigma} <\infty.
$$
Recall on the other hand the definition of the inhomogeneous Besov spaces associated to a semigroup; they were precisely studied in several works, such as \cite{BDY} or \cite{GL}, to name but a few. We shall make an extensive use of these spaces.

\medskip

\begin{definition} \label{def:besov}  {\sf
Fix a positive integer $a$, an exponent $p,q\in (1,\infty)$, and $\sigma\in{\mathbb R}$. A distribution $f\in \SSS$, is said to belong to the \emph{\textbf{Besov space}} $B^\sigma_{p,q}$ if
$$
\|f\|_{B^\sigma_{p,q}} := \big\|e^{-L}f\big\|_p + \left(\int_0^1 t^{-q\,\frac{\sigma}{2}} \, \Big\|Q_t^{(a)}f\Big\|^q_p \, \frac{dt}{t}\right)^\frac{1}{q} <\infty.
$$
This definition of the space does not depend on the integer $a\geq 1$, provided $a$ is big enough.   }
\end{definition}

\medskip

We refer the reader to \cite{BDY} for details about such spaces and a proof of the fact that they do not depend on the parameter $a$ used to define them, provided $a$ is sufficiently large with respect to $\sigma$. The limiting case $p=q=\infty$ leads to the following definition.

\medskip

\begin{definition*}    {\sf
Let a positive integer $a$ be given. For $\sigma\in (-\infty,2)$, a distribution $f\in \SSS$ is said to belong to the space $\calC^\sigma$ if
$$
\|f\|_{\calC^ \sigma} := \Big\|e^{-L} f\Big\|_\infty + \sup_{0<t\leq 1}\,\Big\|Q_t^{(a)} f\Big\|_\infty t^{- \frac{\sigma}{2}} <\infty.
$$
This definition of the space does not depend on the integer $a\geq 1$.   }
\end{definition*}

\medskip

We give in Appendix \ref{Appendix} a simple and self-contained proof that the space $\calC^\sigma$ does not depend on $a$, and that any two norms $\|\cdot\|_{\calC^ \sigma}$, defined with two different values of $a$, are equivalent. The following proposition justifies that we call the spaces $\calC^\sigma$ H\"older space, for all $\sigma<2$, possibly non-positive.

\medskip

\begin{proposition}  \label{prop:caraholder}   {\sf
For $\sigma\in(0,1)$, the spaces $\Lambda^\sigma$ and $\calC^\sigma$ are the same and the two corresponding norms are equivalent.   }
\end{proposition}

\medskip

\noindent We give here a complete proof of this proposition as it provides an elementary illustration of how the properties of the operators $Q_t^{(a)}$ are used to make actual computations. This kind of reasoning and computations will be used repeatedly in the sequel, when working with our paraproduct. Recall that the operators $Q_t^{(a)}$ have kernels $K_{Q_t^{(a)}}$ satisfying Gaussian pointwise estimates, by proposition \ref{PropGaussianPtQt}.

\bigskip

\begin{Dem} 
We divide the proof in two steps, by showing successively that $\Lambda^\sigma$ is continuously injected in $\calC^\sigma$, and that, conversely, $\calC^\sigma$ is continuously injected in $\Lambda^\sigma$.

\medskip 

{\bf Step 1. $\Lambda^\sigma \hookrightarrow\calC^\sigma$.} Note first that since the H\"older space $\Lambda^\sigma$ is made up of bounded functions, it is included in $\mathcal{S}'_o$. Fix an integer $a\geq 1$; then for every $t\in(0,1)$, we have 
\begin{align*} 
\Big(Q_t^{(a)}f\Big)(x)  = \Big(Q_t^{(a)}\big(f(\cdot)-f(x)\big)\Big)(x) = \int K_{Q_t^{(a)}}(x,z) \big(f(z)-f(x)\big)\,\mu(dz).
\end{align*}
For the points $z\in M$, with $d(x,z)\leq \sqrt{t}<1$, we have
$$
\big|f(z)-f(x)\big| \leq d(x,z)^\sigma \|f\|_{\Lambda^ \sigma} \leq  t^{\fr{\sigma}{2}} \|f\|_{\Lambda^ \sigma}
$$
so that
\begin{align*}
\left|\int_{d(x,z)\leq \sqrt{t}} K_{Q_t^{(a)}}(x,z) \big(f(z)-f(x)\big)\,\mu(dz) \right| &\leq t^{\frac{\sigma}{2}} \|f\|_{\Lambda^ \sigma} \int \Big| K_{Q_t^{(a)}}(x,z)\Big|\,\mu(dz) \\
& \lesssim  t^{\frac{\sigma}{2}} \|f\|_{\Lambda^ \sigma},
\end{align*}
since $Q_t^{(a)}$ has a kernel satisfying Gaussian pointwise bounds. The same bounds show that
\begin{align*}
\left|\int_{\sqrt{t} \leq d(x,z)\leq 1} K_{Q_t^{(a)}}(x,z)\right. &\left.\big(f(z)-f(x)\big) \mu(dz) \right| \\  
&\leq \|f\|_{\Lambda^ \sigma}\left(\int_{\sqrt{t} \leq d(x,z)\leq 1} \Big| K_{Q_t^{(a)}}(x,z)\Big|  d(x,z)^\sigma \, \mu(dz) \right)  \\
&\lesssim t^{\sigma \over 2} \|f\|_{\Lambda^ \sigma} \,\int_{\sqrt{t} \leq d(x,z)\leq 1}  \frac{1}{V(x,\sqrt{t})} e^{-c\,\frac{d(x,z)^2}{t}} \left(\frac{d(x,z)}{\sqrt{t}}\right)^\sigma \, \mu(dz)  \\
&\lesssim t^{\sigma \over 2} \|f\|_{\Lambda^ \sigma}.
\end{align*}
Similarly, we have
\begin{align*}
\left|\int_{1 \leq d(x,z)} K_{Q_t^{(a)}}(x,z) \big(f(z)-f(x)\big) \mu(dz) \right| & \leq \|f\|_{\infty} \left(\int_{1 \leq d(x,z)} \Big| K_{Q_t^{(a)}}(x,z) \Big| \, \mu(dz)\right)  \\
& \lesssim e^{-\frac{c}{t}} \, \|f\|_{\Lambda^ \sigma} \\
& \lesssim t^{\sigma \over 2} \, \|f\|_{\Lambda^ \sigma},
\end{align*}
so it comes that the inequality
$$ 
\left|\Big(Q_t^{(a)}f\Big)(x)\right| \lesssim t^{\sigma \over 2} \|f\|_{\Lambda^ \sigma} 
$$
holds uniformly in $t\in(0,1)$, and for every $x\in M$, which proves that $\|f\|_{\calC^\sigma} \lesssim \|f\|_{\Lambda^ \sigma}$.

\medskip

\noindent {\bf Step 2. $\calC^\sigma \hookrightarrow \Lambda^\sigma$.} Let $f\in \calC^ \sigma$ be given. Using the decomposition of the identity provided by Calder\'on reproducing formula
$$
f = e^{-L}f + \int_0^1 Q_t^{(1)}f \,\frac{dt}{t},
$$
we first deduce that $f$ is bounded, with
$$ 
\|f\|_\infty \lesssim \|f\|_{\calC^ \sigma} \left(1 + \int_0^1 t^{\sigma \over 2} \, \frac{dt}{t} \right) \lesssim  \|f\|_{\calC^ \sigma}.
$$
Moreover, for any two points $x,y$, with $0<d(x,y)\leq 1$, we have 
\begin{align*}
 f(x)-f(y) & = \Big\{\big(e^{-L}f\big)(x) - \big(e^{-L}f\big)(y)\Big\} + \int_0^1 \left\{\Big(Q^{(1)}_tf\Big)(x) - \Big(Q^{(1)}_tf\Big)(y) \right\}\,\frac{dt}{t}  \\
 & = \Big\{e^{-L}f(x)-e^{-L}f(y)\Big\} + \Big\{\Big(Q^{(1)}_1f\Big)(x) - \Big(Q^{(1)}_1f\Big)(y)\Big\}  \\
 & \qquad\qquad\qquad\qquad\qquad\hspace{0.3cm}\; + \int_0^1 \Big\{\Big(Q^{(2)}_t f\Big)(x) - \Big(Q^{(2)}_t f\Big)(y)\Big\}\,\frac{dt}{t}. 
\end{align*}
One can use the Lipschitz regularity \eqref{Lipschitz} of the heat kernel to bound the first term in the above sum, giving
\begin{align*}
\Big|e^{-L}f(x)-e^{-L}f(y)\Big| & \lesssim \int \big|p_1(x,z) - p_1(y,z)\big|\,\big|f(z)\big|\,\mu(dz) \\
& \lesssim d(x,y) \,\|f\|_\infty.
\end{align*}
As similar bounds hold for $Le^{-L}$, by analyticity of the heat kernel, the second term admits a similar upper bound. Let now focus on the third term, using  a similar reasoning and noting that $Q^{(2)}_t = 16\,Q_{\frac{t}{2}}^{(2)} Q_{\frac{t}{2}}^{(2)}$. So, for $d(x,y)\leq \sqrt{t}$, we can write
\begin{align*}
\Big| \Big(Q^{(2)}_tf\Big)(x) - \Big(Q^{(2)}_tf\Big)(y)\Big| & \lesssim \int \Big| K_{Q_{\frac{t}{2}}^{(2)}}(x,z)-K_{Q_{\frac{t}{2}}^{(2)}}(y,z)\Big|\,\Big| K_{Q_{\frac{t}{2}}^{(2)}} f(z)\Big| \, \mu(dz) \\
& \lesssim \frac{d(x,y)}{\sqrt{t}}\; \Big\|K_{Q_{\frac{t}{2}}^{(2)}} f\Big\|_\infty \\
& \lesssim \frac{d(x,y)}{\sqrt{t}}\;  t^{\frac{\sigma}{2}} \, \|f\|_{\calC^\sigma}.
\end{align*}
If $\sqrt{t} \leq d(x,y)$, then we directly have
\begin{align*}
\Big| \Big(Q^{(2)}_tf\Big)(x) - \Big(Q^{(2)}_tf\Big)(y)\Big| & \lesssim \Big\| Q_t^{(2)} f \Big\|_\infty  \lesssim  t^{\frac{\sigma}{2}} \|f\|_{\calC^\sigma}.
\end{align*}
Hence, 
\begin{align*}
\left| \int_0^1\Big\{\Big(Q^{(2)}_tf\Big)(x) - \Big(Q^{(2)}_tf\Big)(y)\Big\}\frac{dt}{t} \right| & \lesssim \left(\int_0^{d(x,y)^2} t^{\frac{\sigma}{2}} \, \frac{dt}{t} + \int_{d(x,y)^2}^1 \left(\frac{d(x,y)}{\sqrt{t}}\right) t^{\frac{\sigma}{2}}  \frac{dt}{t} \right) \|f\|_{\calC^\sigma} \\
& \lesssim  d(x,y)^\sigma \, \|f\|_{\calC^\sigma},
\end{align*}
since $\sigma\in(0,1)$. Consequently, we have obtained 
$$
\big| f(x)-f(y)\big| \lesssim d(x,y)^\sigma \|f\|_{\calC^\sigma} 
$$
uniformly for every $x\neq y$ with $d(x,y)\leq 1$, so indeed $\|f\|_{\Lambda^ \sigma} \lesssim \|f\|_{\calC^\sigma}$.
\end{Dem}

\medskip

Our main example of a $\calC^\sigma$ distribution with negative H\"older exponent $\sigma$ will be given by typical realizations of a (possibly weighted) noise over $(M,\mu)$ -- see Proposition \ref{prop:whitenoise}. To prove that regularity property, it will be convenient to assume that the metric measure space $(M,d,\mu)$ has the following property, called \textbf{Ahlfors regularity}. There exists a positive constant $c_1$ such that
$$
V(x,1) \geq c_1,
$$
for all $x\in M$, which, by the doubling property, implies that we have
\begin{equation}
\label{EqDoublingProperty}
V(x,r)\geq c_1 r^\nu,
\end{equation}
for some positive exponent $\nu$, all $x\in M$ and all $0<r\leq 1$. The constant $\nu$ is $d$ on a $d$-dimensional manifold equiped with a smooth measure. This is a relatively weak assumption that is essentially satisfied in a Riemannian setting for closed manifolds without boundary and injectivity radius bounded below by a positive constant or in a sub-domain of the Euclidean space provided that the boundary is Lipschitz. \textit{Under that additional non-degeneracy assumption} on the volume measure, we have the following Besov embedding, proved in Appendix \ref{Appendix}.

\medskip

\begin{lemma}[Besov embedding]   \label{LemBesovEmbedding} {\sf   
Given $-\infty<\sigma<2$, and $1<p<\infty$, we have the following continuous embeddings.
$$ 
B^{\sigma}_{p,p} \hookrightarrow B^{\sigma}_{p,\infty} \hookrightarrow B^{\sigma - \frac{\nu}{p}}_{\infty,\infty}=\calC^{\sigma - \frac{\nu}{p}}
$$   }
\end{lemma}

\medskip

Besov embedding can be used in a very efficient way to investigate the regularity properties of random Gaussian fields, as will be illustrated in section \ref{SubsectionRenormalization}.

\medskip

\noindent {\sf \textbf{Remark.}} {\it Let us point out here that our H\"older spaces $C^\sigma$, with $\sigma < 0$, coincide in the Euclidean setting with those used by Hairer \cite{H}, and, roughly speaking, defined in terms of scaling properties of the pairing of a distribution with a one parameter family of rescaled functions. Indeed, on the Euclidean space it is known that to define Besov spaces or H\"older spaces through Littlewood-Paley functionals, we may chose any good Fourier multipliers satisfying suitable conditions; the latter are satisfied by the derivatives $\big(Q_t^{(a)}\big)_{t}$ of the heat semigroup. So our spaces correspond to the standard inhomogeneous spaces defined by any Littlewood-Paley functionals. From wavelet or frame characterization \cite{Meyer}, we can conclude that our H\"older space coincides with those used in \cite{H} or \cite{HL}.}

\medskip

Before turning to the definition of our paraproduct, we close this section with two continuity properties involving the H\"older spaces $\calC^\sigma$, which we shall use in the sequel.

\medskip

\begin{proposition}  \label{prop:low}   {\sf 
For any $\sigma\in(-\infty,2)$, and every integer $a\geq 0$, we have 
$$
\Big\|P_1^{(a)} f\Big\|_{\infty} \lesssim \big\|f\big\|_{\calC^\sigma}.
$$  }
\end{proposition}

\medskip

\begin{Dem} 
We have by construction $P_1^{(a)} = \big(1+ \alpha_1L+ \cdots + \alpha_{a-1} L^{a-1}\big) e^{-L}$, for some coefficients $\alpha_1, \dots, \alpha_{a-1}$. As we have by definition $\big\| e^{-L} f\big\|_{\infty} \lesssim \|f\|_{\calC^\sigma}$, and $L^\ell e^{-L} = Q^{(\ell)}_1$, for $\ell=1\dots (a-1)$, we see that $\big\| L^\ell e^{-L} f \big\|_\infty \lesssim  \|f\|_{\calC^\sigma}$, since we have seen above that we can choose the parameter $a$ in the definition of the H\"older space.
\end{Dem}

\medskip

\begin{proposition} \label{prop:gradient}  {\sf  
For $\sigma\in(-\infty,1)$, we have
$$ 
\sup_{t\in(0,1]} t^{-\frac{\sigma}{2}} \Big\|\big(\sqrt{t}\,\Gamma\big)\Big(e^{-tL}f\Big)\Big\|_{\infty} \lesssim \|f\|_{\calC^\sigma}.
$$
The same conclusion holds with any of the operators $P_t^{(a)}$ in the role of $e^{-tL}$.  }
\end{proposition}

\medskip

\begin{Dem} 
Given $t\in(0,1]$, use Calder\'on reproducing formula to write
$$ 
\Big|\big(\sqrt{t}\,\Gamma\big)\big(e^{-tL}f\big)\Big|  \lesssim \int_0^1\Big| \big(\sqrt{t}\,\Gamma\big)\big(e^{-tL} Q_s^{(1)}f\big)\Big| \,\frac{ds}{s}  + \Big|\big(\sqrt{t}\,\Gamma\big)\big(e^{-(1+t)L} f\big)\Big|.
$$
We divide the integration interval in the above-right hand side into $(0,t)$ and $[t,1]$ to bound that term. For $s<t$, we have $e^{-tL} Q_s^{(1)} = \frac{s}{\frac{t}{2}+s}\,e^{-\frac{t}{2} L}\,Q_{s+\frac{t}{2}}^{(1)}$, so we can use Proposition  \ref{prop:gradientsemigroup} to get
\begin{align*} 
\Big\|\big(\sqrt{t}\,\Gamma\big)\big(e^{-tL} Q_s^{(1)} f\big)\Big\|_{\infty}  &  \lesssim \frac{s}{t} \, \Big\|\big(\sqrt{t}\,\Gamma\big)\big(e^{-\frac{t}{2}L}\big)\Big\|_{\infty \to \infty} \, \Big\| Q_{s+\frac{t}{2}}^{(1)} f\Big\|_\infty \\
& \lesssim \frac{s}{t} \,t^{\frac{\sigma}{2}} \, \|f\|_{\calC^\sigma}.
\end{align*}
Similarly for $t\leq s$, then $e^{-tL} Q_s^{(1)} = e^{-\frac{s}{2}\,L}\,\frac{s}{t+\frac{s}{2}}\,Q_{\frac{s}{2}+t}^{(1)}$, and we have
\begin{align*} 
\Big\|\big(\sqrt{t}\,\Gamma\big)\big(e^{-tL} Q_s^{(1)} f\big)\Big\|_{\infty} &  \lesssim  \left(\frac{t}{s}\right)^{\frac{1}{2}} \, \Big\|\big(\sqrt{s}\,\Gamma\big)\big(e^{-\frac{s}{2}\,L}\big)\Big\|_{\infty \to \infty} \, \Big\| Q_{\frac{s}{2}+t}^{(1)}f\Big\|_\infty \\
& \lesssim \left(\frac{t}{s}\right)^{\frac{1}{2}} s^{\frac{\sigma}{2}} \, \|f\|_{\calC^\sigma}.
\end{align*}
Similar computations give the estimate
$$
\Big\|\big(\sqrt{t}\,\Gamma\big)\big(e^{-(1+t)L} f\big) \Big\|_{\infty} \lesssim  \sqrt{t}\,\|f\|_{\calC^\sigma}.
$$
We conclude by integrating with respect to $s\in(0,1)$, using here the fact that $\sigma<1$.
\end{Dem}

\bigskip

\section[\hspace{0.6cm} {\sf Paraproduct and commutator estimates in H\"older spaces}]{Paraproduct and commutator estimates in H\"older spaces} 
\label{sec:esti}

\subsection[\hspace{-1cm} {\sf Paraproducts based on the semigroup}]{Paraproducts based on the semigroup}
\label{SubsectionParaproducts}

Bony's paraproduct machinery has its roots in the Littlewood-Paley decomposition of any distribution $f$ as a sum of smooth functions $\Delta_i f$ localized in the frequency space, so a product $fg$ of any two distributions can formally be decomposed as
\begin{equation}
\label{EqProductLPDecomposition}
fg = \sum_{i,j} \Delta_i f \,\Delta_j g = \sum_{|i-j|\geq 2} \Delta_i f \,\Delta_j g + \sum_{|i-j|\leq 1} \Delta_i f \,\Delta_j g =: (1) + \Pi(f,g)
\end{equation}
into a sum of products of two functions oscillating on different scales, and an a priori resonant term $\Pi(f,g)$. This decomposition draws its usefulness from some relatively elementary a priori estimates that show that the term $(1)$ above makes sense and is well-controlled under extremely general conditions, while the resonant term $\Pi(f,g)$ can be shown to define a continuous map from $\calC^\alpha\times\calC^\beta$ to $\calC^{\alpha+\beta}$, provided $\alpha+\beta > 0$. These estimates rely crucially on some properties inherited from the very definition of the Littlewood-Paley blocks as Fourier projectors. These properties cannot be grasped so easily in our semigroup setting; however, we shall use the operators $P^{(a)}_t, Q^{(a)}_t$ and $\sqrt{t}\Gamma$ or $(tL)P^{(a)}_t$ as frequency projectors, with $P^{(a)}_t$ projecting on frequencies lower than or equal to $t^{-\frac{1}{2}}$, and $Q^{(a)}_t, \sqrt{t}\Gamma$ or $(tL)P^{(a)}_t$ as localizing at frequencies of order $t^{-\frac{1}{2}}$. Indeed, in the torus, and working with the Euclidean Laplacian, the operator $Q^{(a)}_t$ has for instance a Fourier transform equal to 
$$
\widehat{Q^{(a)}_t}(\lambda) = \big(t|\lambda|^2\big)^ae^{-t|\lambda|^2};
$$
it is essentially localized in an annulus $|\lambda|\simeq t^{-\frac{1}{2}}$. Similar explicit Fourier pictures for the other operators can be given in the setting of the torus. This 'frequency' interpretation of these operators will be our main guide in the definition of our paraproduct given below. This paraproduct will depend on a choice of a positive integer-valued parameter $b$ that can be tuned on demand in any given problem; the bigger $b$ is, the more we can do some 'integration by parts' -- it will be fixed at some point, and the reader is invited not to bother about its value. To clarify notations, we shall repeatedly use below the notation $f\cdot g$ for the usual product of two functions.

\medskip

Rather than starting with Bony's decomposition \eqref{EqProductLPDecomposition}, we take as a starting point Calderon's reproducing formula
\begin{equation}
\label{EqDecompCalderonFG}
\begin{split}
fg &= \lim_{t\to 0} P_t^{(b)} \left( P_t^{(b)} f \cdot P_t^{(b)} g \right) = -\int_0^1 t \partial_t \,\left\{ P_t^{(b)}\left( P_t^{(b)} f \cdot P_t^{(b)} g \right)\right\} \, \frac{dt}{t} + \Delta_{-1}(f,g)  \\
& = \frac{1}{\gamma_b} \bigintssss_0^1 \Bigg\{P_t^{(b)}\left(Q_t^{(b)} f\cdot P_t^{(b)} g\right) + P_t^{(b)} \left(P_t^{(b)} f\cdot Q_t^{(b)} g\right) + Q_t^{(b)} \left(P_t^{(b)} f\cdot P_t^{(b)} g\right) \Bigg\}\, \frac{dt}{t} \\
&\hspace{0.48cm}+ \Delta_{-1}(f,g),
\end{split}
\end{equation}
where
$$
\Delta_{-1}(f,g):= P_1^{(b)}\left(P_1^{(b)} f \cdot P_1^{(b)} g \right) 
$$
stands for the low-frequency part of the product of $f$ and $g$, and where we implicitly make the necessary assumptions on $f$ and $g$ for the above formula to make sense.

\medskip

Guided by the above heuristic argument about the role of the operators $P^{(a)}_t, Q^{(a)}_t$, etc. as frequency projectors, we decompose the term involving the product of $P^{(a)}_tf$ and $P^{(a)}_tg$, by using the definition of the carr\'e du champ operator $\Gamma$
$$
L\big(\phi_1 \cdot \phi_2\big) = L(\phi_1)\,\phi_2 + L(\phi_2)\,\phi_1 - 2 \Gamma(\phi_1,\phi_2)
$$
and write
{\small   \begin{equation*}
\begin{split}
&Q_t^{(b)} \left(P_t^{(b)} f\cdot P_t^{(b)} g\right) \\ 
&= Q_t^{(b-1)} \left((tL)P_t^{(b)} f\cdot P_t^{(b)} g\right) + Q_t^{(b-1)} \left(P_t^{(b)} f\cdot (tL)P_t^{(b)} g\right) - 2 Q_t^{(b-1)} \Gamma\left(\sqrt{t} P_t^{(b)} f, \sqrt{t} P_t^{(b)} g\right) \\
&=: B_g(f) + B_f(g) + R(f,g).
\end{split}
\end{equation*}   }
If one rewrites identity \eqref{EqDecompCalderonFG} under the form 
$$
fg =: \int_0^1\big\{(1) + (2) + (3) \big\}\,\frac{dt}{t} + \Delta_{-1}(f,g) 
$$
with obvious notations, this suggest to decompose it as
$$
fg = \int_0^1 \Big(\big\{(1) + B_g(f)\big\} + \big\{(2) + B_f(g)\big\} + R(f,g) \Big)\,\frac{dt}{t} + \Delta_{-1}(f,g)
$$
and to identify the integral of the terms into brackets in the above formula as paraproducts, and by defining the resonant term as the integral of $R(f,g)$. This is what was done in \cite{BS} where this notion of paraproduct, introduced in \cite{B-T1}, was shown to have nice continuity properties in H\"older spaces $\calC^\alpha$, provided one deals with positive exponents $\alpha$. Given our needs to deal with negative exponents, a refinement of this decomposition seems to be needed to get some continuity properties for negative exponent as well. We thus use the carr\'e du champ formula in each term $(1)$ and $(2)$, and decompose $(1)$ under the form
\begin{equation*}
\begin{split}
&(tL)P_t^{(b)} \left(Q_t^{(b-1)} f\cdot P_t^{(b)} g\right) + \Big\{2P_t^{(b)} \big(t\Gamma\big)\left(Q_t^{(b-1)} f, P_t^{(b)} g\right) - P_t^{(b)} \left(Q_t^{(b-1)} f\cdot (tL) P_t^{(b)} g\right) \Big\}  \\
      &=: A_g(f) + S(f,g),
\end{split}
\end{equation*}
with $S(f,g)$ the sum of the two terms into bracket, and 
$$
(2) = A_f(g) + S(g,f).
$$
Note that the functions $A_f(g), S(f,g),\dots$ all depend implicitly on time. This decomposition leads to the following definition.

\medskip

\begin{definition*} {\sf
Given an integer $b\geq 2$ and $f\in \bigcup_{s\in (0,1)} \calC^s$ and $g\in L^\infty$, we define their paraproduct by the formula
\begin{align*}
\Pi^{(b)}_g(f) &= \frac{1}{\gamma_b} \int_0^1 \Big\{A_g(f) + B_g(f)\Big\} \,\frac{dt}{t} \\
                        & = \frac{1}{\gamma_b} \int_0^1 \left\{ (tL)P_t^{(b)} \left(Q_t^{(b-1)} f\cdot P_t^{(b)} g\right) + Q_t^{(b-1)} \left((tL)P_t^{(b)} f\cdot P_t^{(b)} g\right) \right\}\, \frac{dt}{t}.
\end{align*} }
\end{definition*}

\medskip

The well-defined character of this integral is proved in proposition \ref{prop:para} below. With this notation, Calderon's formula becomes
$$ 
fg= \Pi^{(b)}_g(f)+\Pi^{(b)}_f(g)+\Pi^{(b)}(f,g)+\Delta_{-1}(f,g)
$$
with the 'low-frequency part'
$$ 
\Delta_{-1}(f,g):= P_1^{(b)}\left(P_1^{(b)} f \cdot P_1^{(b)} g \right) 
$$
and the 'resonant term'
\begin{align*} 
\Pi^{(b)}(f,g) &= \frac{1}{\gamma_b} \int_0^1 \Big\{S(f,g) + S(g,f) + R(f,g) \Big\} \,\frac{dt}{t}  \\
 & = \frac{1}{\gamma_b} \int_0^1 \left\{-P_t^{(b)} \left(Q_t^{(b-1)} f \cdot (tL) P_t^{(b)} g\right) + 2P_t^{(b)} \Gamma \Big(\sqrt{t}\,Q_t^{(b-1)} f, \sqrt{t}\,P_t^{(b)} g\Big)\right\}\,\frac{dt}{t} \\
  &\quad + \frac{1}{\gamma_b} \int_0^1 \left\{-P_t^{(b)} \left((tL)P_t^{(b)} f \cdot Q_t^{(b-1)} g\right) + 2P_t^{(b)}\Gamma \left(\sqrt{t}\,P_t^{(b)} f, \sqrt{t}\, Q_t^{(b-1)} g\right) \right\}\, \frac{dt}{t}\\
  &\quad - \frac{1}{\gamma_b} \int_0^1 2 Q_t^{(b-1)} \Gamma\left(\sqrt{t} P_t^{(b)} f, \sqrt{t} P_t^{(b)} g\right) \, \frac{dt}{t}.
\end{align*}
Note that we have $\Pi^{(b)}_{{\bf 1}}(\cdot) = \textrm{Id}-\Delta_{-1}(\cdot,{\bf 1})$, as a consequence of our choice or renormalizing constant.

\medskip

\subsection[\hspace{-1cm} {\sf Paraproduct estimates}]{Paraproduct estimates}
\label{SubsectionParapdctEstimates}

We prove in this paragraph the basic continuity estimates satisfied by the maps defined by the low frequency part, the paraproduct and the resonant term. 

\medskip

\begin{proposition} \label{prop:low2}  {\sf 
Fix an integer $b\geq 2$. For any real numbers $\alpha,\beta$ and every positive $\gamma$, we have
\begin{equation} 
\big\| \Delta_{-1}(f,g)\big\|_{\calC^\gamma} \lesssim  \|f\|_{\calC^\alpha} \|g\|_{\calC^\beta}.  
\end{equation}
for every $f \in \calC^\alpha$ and $g\in \calC^\beta$.   }
\end{proposition}   

\medskip

\begin{Dem}
Consider the collection $\Big(Q^{(a)}_s\Big)_{0<s\leq 1}$ for a large enough integer $a\geq \gamma$. Then
$$ 
Q_s^{(a)} \Delta_{-1}(f,g) = Q_s^{(a)} P_1^{(b)}\left(P_1^{(b)} f \cdot P_1^{(b)} g \right) .
$$
Since $s\leq 1$, we have $Q_s^{(a)}P_t^{(b)} = s^a e^{-sL} L^a P_1^{(b)}$, with the operator $L^a P_1^{(b)}$ bounded on $L^\infty$. We obtain the conclusion from Proposition \ref{prop:low} as we have
\begin{align*} 
\big\|Q_s^{(a)} \Delta_{-1}(f,g)\big\|_{\infty} & \lesssim s^a \big\|P_1^{(b)} f\big\|_\infty \big\|P_1^{(b)} g \big\|_\infty \\
  & \lesssim s^\gamma \|f\|_{\calC^\alpha} \|g\|_{\calC^\beta}.
\end{align*}
\end{Dem}

\ssk

The continuity properties of the paraproduct are given by the following statement; they are the exact analogue of their classical counterpart, based on Littlewood-Paley decomposition, as can be found for instance in the textbook \cite{BCD} of Bahouri, Chemin and Danchin.

\medskip

\begin{proposition}  \label{prop:para}   {\sf 
Fix an integer $b\geq 2$. For any $\alpha\in (-2,1)$ and $f\in \calC^\alpha$, we have
\begin{itemize}
 \item for every $g\in L^\infty$
\begin{equation} 
\Big\| \Pi^{(b)}_g(f)\Big\|_{\calC^\alpha} \lesssim \|g\|_\infty \|f\|_{\calC^\alpha} 
\label{eq:para1} \end{equation}
 \item for every $g\in\calC^\beta$ with $\beta<0$ and $\alpha+\beta\in (-2,1)$  
\begin{equation} 
\Big\| \Pi^{(b)}_g(f)\Big\|_{\calC^{\alpha+\beta}} \lesssim \|g\|_{\calC^\beta} \|f\|_{\calC^\alpha}. 
\label{eq:para2} \end{equation}
\end{itemize}  }
\end{proposition}  

\medskip

The range $(-2,1)$ for the regularity exponent can appear as unusual since in the standard Euclidean theory such continuity properties hold for every $\alpha\in{\mathbb R}$. However, as explained in section \ref{SubsectionHeatSemiGroup}, the restriction $\alpha<1$ comes from our optimal / minimal setting where we only assume a gradient estimate on the heat kernel. The restriction $\alpha>-2$ can be explained as follows. In the Euclidean theory, nice Fourier multipliers can be used to have a 'perfect' frequency decomposition and the study of paraproducts mainly relies on the following rule: {\emph{the spectrum of the product of two functions is included into the sum of the two spectrums}}; this comes from the group structure through the Fourier representation of the convolution. In our setting, our frequency decomposition involving the heat semigroup is not so perfect and the previous rule on the spectrum does not hold, at least not in such a 'perfect' sense. That is why the new limitation $\alpha>-2$ appears; it is inherent to the semigroup approach developed here. No such limitation holds in the more restricted setting developed in \cite{BBF16}.

\medskip

\begin{Dem}
Recall that
$$  
\Pi^{(b)}_g(f) = \frac{1}{\gamma_b} \int_0^1  (tL)P_t^{(b)} \left(Q_t^{(b-1)} f\cdot P_t^{(b)} g\right) + Q_t^{(b-1)} \left((tL)P_t^{(b)} f\cdot P_t^{(b)} g\right) \, \frac{dt}{t}.
$$
Given $s\in(0,1]$, consider $Q_s^{(b-1)} \Pi_g(f)$. For $s\leq t$, we use that 
$$ 
Q_s^{(b-1)} (tL) P_t^{(b)} =  \left(\frac{s}{t}\right)^{b-1} (tL)^b P_t^{(b)} e^{-sL} \quad \textrm{and} \quad Q_s^{(b-1)} Q_t^{(b-1)} =  \left(\frac{s}{t}\right)^{b-1} Q_{t}^{2(b-1)} e^{-sL}, 
$$
and for $t\leq s$ that
$$ 
Q_s^{(b-1)} (tL)P_t^{(b)} =  \frac{t}{s} Q_{s}^{(b)} P_t^{(b)} \quad \textrm{and} \quad Q_s^{(b-1)} Q_t^{(b-1)} =  \frac{t}{s} Q_{s}^{(b)} Q_t^{(b-2)}. 
$$
Hence, with the uniform $L^\infty$-boundedness of $Q_t,P_t$ operators, we have
\begin{align*}
\Big\| Q_s^{(b-1)}&\Pi^{(b)}_g(f) \Big\|_\infty \lesssim \int_0^s \frac{t}{s} \Big\| Q_t^{(b-1)} f\Big\|_\infty \Big\|P_t^{(b)} g\Big\|_\infty + \frac{t}{s} \Big\| (tL) P_t^{(b)} f\Big\|_\infty \Big\|P_t^{(b)} g\Big\|_\infty \, \frac{dt}{t} \\
  & + \int_s^1 \left(\frac{s}{t}\right)^{b-1} \Big\| Q_t^{(b-1)} f\Big\|_\infty \Big\|P_t^{(b)} g\Big\|_\infty + \left(\frac{s}{t}\right)^{b-1} \Big\| (tL) P_t^{(b)} f\Big\|_\infty \Big\|P_t^{(b)} g\Big\|_\infty \, \frac{dt}{t}.
\end{align*}
Since $f\in \calC^\alpha$ we have
$$ 
\Big\| Q_t^{(b-1)} f\Big\|_\infty + \Big\| (tL) P_t^{(b)} f\Big\|_\infty \lesssim t^{\alpha \over 2}\,\|f\|_{\calC^\alpha}.
$$
Moreover, if $g\in L^\infty$ then
$$ 
\Big\| P_1^{(b)} g\Big\|_\infty \lesssim \|g\|_\infty
$$
and if $g\in \calC^\beta$ with $\beta<0$ then
\begin{align*}
\Big\| P_t^{(b)} g \Big\|_\infty & \leq \int_t^1 \Big\| Q^{(b)}_u g\Big\|_\infty \, \frac{du}{u} + \Big\| P^{(b)}_1(g) \Big\|_\infty \\
& \lesssim \left(\int_t^1  \frac{du}{u^{1-\frac{\beta}{2}}} + 1 \right) \|g\|_{\calC^\beta} \lesssim t^\frac{\beta}{2} \,\|g\|_{\calC^\beta}.
\end{align*}
We deduce the following bounds as a consequence.

\ssk

\begin{itemize}
\item If $g\in L^\infty$ then
\begin{align*}
\Big\| Q_s^{(b-1)}\Pi^{(b)}_g(f) \Big\|_\infty  & \lesssim \left(\int_0^s \left(\frac{t}{s}\right) t^{\alpha \over 2} \, \frac{dt}{t} + \int_s^1 \left(\frac{s}{t} \right)^c t^{\alpha \over 2} \, \frac{dt}{t}\right) \|f\|_{\calC^\alpha} \|g\|_\infty \\
& \lesssim s^{\alpha \over 2}\|f\|_{\calC^\alpha} \|g\|_\infty,
\end{align*}
since $\alpha\in(-2,1)$ and $c\geq 1$. This holds for every $s>0$ which yields \eqref{eq:para1}.
\item If $g\in \calC^\beta$ with $\alpha+\beta\in(-2,1)$ then
\begin{align*}
\Big\| Q_s^{(b-1)} \Pi^{(b)}_g(f) \Big\|_\infty  & \lesssim \left(\int_0^s \left(\frac{t}{s}\right) t^{\frac{\alpha+\beta}{2}} \, \frac{dt}{t} + \int_s^1 \left(\frac{s}{t}\right)^{b-1} t^{\frac{\alpha + \beta}{2}}\, \frac{dt}{t}\right) \|f\|_{\calC^\alpha} \|g\|_{\calC^\beta} \\
& \lesssim s^{\frac{\alpha+\beta}{2}}\|f\|_{\calC^\alpha} \|g\|_\infty,
\end{align*}
since $2(b-1)\geq 1> \alpha+\beta >-2$. This holds for every $s>0$, which yields \eqref{eq:para2}.
\end{itemize}
\end{Dem}

\medskip

\begin{proposition} \label{prop:diag}   {\sf 
Fix an integer $b>2$. 
For any $\alpha,\beta\in (-\infty,1)$ with $\alpha+\beta>0$, for every $f\in\calC^\alpha$ and $g\in \calC^\beta$, we have the continuity estimate
$$ 
\Big\| \Pi^{(b)}(f,g) \Big\|_{\calC^{\alpha+\beta}} \lesssim \|f\|_{\calC^\alpha} \|g\|_{\calC^\beta}.
$$   }
\end{proposition}   

\medskip

\begin{Dem}
We recall that
\begin{align*} 
\Pi^{(b)}(f,g)  & = \frac{1}{\gamma_b} \int_0^1 -P_t^{(b)} \left(Q_t^{(b-1)} f\cdot (tL) P_t^{(b)} g\right) + 2P_t^{(b)} \Gamma\left(\sqrt{t} Q_t^{(b-1)} f, \sqrt{t} P_t^{(b)} g\right) \, \frac{dt}{t} \\
  & + \frac{1}{\gamma_b} \int_0^1 -P_t^{(b)} \left((tL)P_t^{(b)} f\cdot Q_t^{(b-1)} g\right) + 2P_t^{(b)} \Gamma\left(\sqrt{t} P_t^{(b-1)} f, \sqrt{t}  Q_t^{(b)} g\right) \, \frac{dt}{t}\\
  & + \frac{1}{\gamma_b} \int_0^1 2 Q_t^{(b-1)} \Gamma\left(\sqrt{t}P_t^{(b)} f, \sqrt{t} P_t^{(b)} g\right) \, \frac{dt}{t}.
\end{align*}
Consider the function $Q_s^{(b-1)} \Pi^{(b)}(f,g)$, for every $s\in(0,1]$. It is given by an integral over $(0,1)$, which we split into $(I)$ an integral over $(0,s)$, and $(II)$ an integral over $(s,1)$. Since $f\in \calC^\alpha$, the use of Proposition \ref{prop:gradient}, with $\alpha<1$,  yields the estimate
\begin{align*}
 & \Big\| Q_t^{(b-1)} f \Big\|_\infty + \Big\| \sqrt{t}\Gamma\big(Q_t^{(b-1)} f\big) \Big\|_\infty + \Big\| (tL)P_t^{(b)} f\Big\|_\infty + \Big\| \sqrt{t}\Gamma\big(Q_t^{(b-1)} f\big)\Big\|_\infty \\
 & \qquad + \Big\|\sqrt{t}\Gamma\big(P_t^{(b)} f\big)\Big\|_\infty \lesssim t^{\alpha\over 2} \, \|f\|_{\calC^\alpha}; 
\end{align*}
a similar estimate holds with $g$ in place of $f$, and $\beta$ in place of $\alpha$. Using the uniform $L^\infty$-boundedness of the different approximation operators, we get for the first part
\begin{align*}
\big\| Q_s^{(b-1)} (I) \big\|_\infty & \lesssim \left(\int_0^s t^{\frac{\alpha+\beta}{2}} \, \frac{dt}{t}\right) \|f\|_{\calC^\alpha} \|g\|_{\calC^\beta} \\
   &   \lesssim s^{\frac{\alpha+\beta}{2}} \, \|f\|_{\calC^\alpha} \|g\|_{\calC^\beta}  ,
 \end{align*}  
where we used the strict inequality $\alpha+\beta>0$. For the second part, we observe that for $t>s$ then 
$$
Q_s^{(b-1)} P_t = \left(\frac{s}{t}\right)^c e^{-sL} (tL)^{b-1} P_t \quad \textrm{and} \qquad Q_s^{(b-1)} Q_t^{(b-1)} =  \left(\frac{s}{t}\right)^{b-1} Q_t^{2(b-1)}e^{-sL}.
$$ 
So we get for the second part
\begin{align*}
\big\| Q_s^{(b-1)} (II) \big\|_\infty & \lesssim \left(\int_s^1 t^{\frac{\alpha+\beta}{2}} \left(\frac{s}{t}\right)^{b-1} \, \frac{dt}{t}\right) \|f\|_{\calC^\alpha} \|g\|_{\calC^\beta}  \\
& \lesssim s^{\frac{\alpha+\beta}{2}}\|f\|_{\calC^\alpha} \|g\|_{\calC^\beta},
\end{align*}
using the fact that $2(b-1)\geq 2 >\alpha+\beta$.
\end{Dem}

\bigskip

\subsection[\hspace{-1cm} {\sf Commutator estimates}]{Commutator estimates}
\label{SubsectionCommutatorEstimates}

Recall the discussion of the paracontrolled calculus approach to the study of the parabolic Anderson model equation given in the introductory section \ref{SubsectionParacontrolledCalculus}. We have at that point the paraproduct and diagonal operators in hands; we deal in this section with the fundamental commutator estimate introduced in \cite{GIP}. Readers familiar with the basics of stochastic analysis will notice the similarity of this continuity result and the rule satisfied by the bracket operator in It\^o's theory
$$
d\left\langle \int_0^\cdot MdN\,,P\right\rangle = P\,d\langle N, P\rangle;
$$
this is obviously not a coincidence.

\medskip

\begin{proposition}  \label{prop:commutator}   {\sf
Consider the a priori unbounded trilinear operator
$$ 
C(f,g,h):=\Pi^{(b)}\Big(\Pi^{(b)}_g(f),h\Big) - g \, \Pi^{(b)}(f,h), 
$$
on $\SSS$. Let $\alpha,\beta,\gamma$ be H\"older regularity exponents with $\alpha\in(-1,1), \beta\in(0,1)$ and $\gamma\in (-\infty,1]$. If
$$
0 < \alpha + \beta + \gamma  \qquad \textrm{ and }\qquad \alpha + \gamma < 0
$$
then, setting $\delta := (\alpha+\beta)\wedge 1 +\gamma$, we have 
\begin{equation}
\label{EqCommutator}
\big\| C(f,g,h) \big\|_{\calC^{\delta}} \lesssim \|f\|_{\calC^\alpha} \, \|g\|_{\calC^\beta} \, \|h\|_{\calC^\gamma},
\end{equation}
for every $f\in \calC^\alpha$ ,$g\in\calC^\beta$ and $h\in \calC^\gamma$; so the commutator defines a continuous trilinear map from $\calC^\alpha\times\calC^\beta\times\calC^\gamma$ to $\calC^\delta$.  }
\end{proposition}

\medskip

\begin{Dem}
Note first that the paraproduct $\Pi_g(f)$ is given, up to a multiplicative constant, by the sum of two terms of the form 
$$ 
\mcA(f,g) = \int_0^1 \mcQ^1_t\Big(\mcQ^2_t f \cdot \mcP_t g \Big) \, \frac{dt}{t},
$$
and the resonant part $\Pi(f,g)$ by the sum of five terms of the following forms
\begin{equation}
\label{EqModelTermR}
\mcR(f,g) = \int_0^1 \mcP^1_t \Gamma\left( \sqrt{t} \mcP^2_t f\, ,\, \sqrt{t} \mcP_t^3 g \right) \, \frac{dt}{t},
\end{equation}
or
$$ 
\mcR(f,g) = \int_0^1 \mcP_t \left( (tL) \mcP_t f  \cdot \mcQ_t g \right) \, \frac{dt}{t},
$$
or 
$$ 
\mcR(f,g) = \int_0^1 \mcP^1_t \left(\mcQ_t f  \cdot (tL)\mcP_t^2 g \right) \, \frac{dt}{t},
$$
where the operators \vspace{0.1cm}
\begin{itemize}
   \item $\mcQ_t, \mcQ^j_t$ are of the form $(tL)^{b-1} p(tL)\,e^{-tL}$ with a polynomial function $p$,   \vspace{0.1cm}

   \item $\mcP_t, \mcP^j_t$ are of the form $p(tL)\,e^{-tL}$ with a polynomial function $p$.   \vspace{0.1cm}
\end{itemize}
Note also that terms of the for $\psi(tL)$ are a posteriori of the form $\phi(tL)$. So it suffices to focus on a generic term of the form
$$ 
D(f,g,h) := \mcR\big(\mcA(f,g),h\big) - g\,\mcR(f,h)
$$
and prove the continuity estimate \eqref{EqCommutator} for it. We focus on the case where $\mcR$ has form \eqref{EqModelTermR}, the treatment of the other cases being similar and somewhat easier. We split the proof of the commutator estimate \eqref{EqCommutator} for $D$ in two steps, and introduce an intermediate quantity 
$$
\mcS(f,g,h) := \int_0^1 \mcP^1_t\Big\{ \Gamma\big(\sqrt{t} \mcP^2_t f , \sqrt{t} \mcP_t^3 h\big) \cdot P_t g\Big\} \, \frac{dt}{t}
$$
for which we shall prove that we have both 
\begin{equation}
\label{EqEstimateD1}
\big\| g\,\mcR(f,h) - \mcS(f,g,h)\big\|_{\calC^\delta} \lesssim \|f\|_{\calC^\alpha} \, \|g\|_{\calC^\beta} \, \|h\|_{\calC^\gamma}
\end{equation}
and 
\begin{equation}
\label{EqEstimateD2}
\big\| D(f,g,h) - \mcS(f,g,h)\big\|_{\calC^\delta} \lesssim \|f\|_{\calC^\alpha} \, \|g\|_{\calC^\beta} \, \|h\|_{\calC^\gamma}.
\end{equation}

\medskip

\noindent {\bf Step 1 -- proof of \eqref{EqEstimateD1}.} We first prove a weaker version of the continuity estimate \eqref{EqEstimateD1}, under the form of the inequality
\begin{equation} 
\label{eq:eq1}
\big\| g\,\mcR(f,h) - \mcS(f,g,h)\big\|_\infty \lesssim \|f\|_{\calC^\alpha} \, \|g\|_{\calC^\beta} \, \|h\|_{\calC^\gamma}.
\end{equation}
As a start, remark that we have
\begin{equation} 
\label{eq:eq2} 
\Big(g \mcR(f,h) - \mcS(f,g,h)\Big)(x) = \int_0^1 \mcP^1_t\Big\{ \Gamma\big(\sqrt{t} \mcP^2_t f , \sqrt{t} \mcP_t^3 h\big) \cdot \big(g(x) - P_t g\big)\Big\}(x) \, \frac{dt}{t},
\end{equation}
for $\mu$-almost every $x\in M$. Since $g\in \calC^\beta$, with $\beta\in (0,1)$, we have
\begin{align*}
\| P_t g -g \|_\infty & \lesssim \int_0^t \| Q_s g \|_\infty \, \frac{ds}{s} \\
 &\lesssim \left(\int_0^t  s^{\frac{\beta}{2}}\, \frac{ds}{s}\right) \|g\|_{\calC^\beta} \lesssim t^{\frac{\beta}{2}}\,\|g\|_{\calC^\beta},
\end{align*} 
so we have
\begin{align*}
\big|P_t g(y) - g(x)\big| &\leq \big| P_t g(y) - g(y)\big| + \big| g(y) - g(x)\big| \\
& \lesssim t^{\frac{\beta}{2}}\,\|g\|_{\calC^\beta} + d(x,y)^\beta \,\|g\|_{\calC^\beta}\\
& \lesssim \left(t^{\frac{\beta}{2}} + d(x,y)^\beta\right)  \|g\|_{\calC^\beta}, 
\end{align*}
for every $x,y\in M$. Coming back to equation \eqref{eq:eq2} and using Gaussian pointwise estimates for the kernel of $\mcP^1_t$, together with Proposition \ref{prop:gradient} and \eqref{eq:CS}, we have
$$
\Big| \mcP^1_t \Big(  \Gamma\big(\sqrt{t} \mcP^2_t f , \sqrt{t} \mcP_t^3 h\big) \cdot \big(g(x) - P_t g\big)\Big)(x) \Big| 
$$
no greater than
\begin{equation*}
\begin{split}
&\left\{\int_M \frac{1}{V(x,\sqrt{t})}\exp\left(-c \frac{d(x,y)^2}{t}\right) \big|g(x) - P_t g(y)\big| \, d\mu(y) \right\}  \big\|\sqrt{t} \Gamma\big(\mcP^2_t f\big)\big\|_\infty \big\|\sqrt{t} \Gamma\big(\mcP^3_t h\big)\big\|_\infty  \\
&\lesssim \left\{\int_M \frac{1}{V(x,\sqrt{t})}\exp\left(-c \frac{d(x,y)^2}{t}\right) \left(t^{\frac{\beta}{2}} + d(x,y)^\beta\right) \, d\mu(y) \right\} \,t^\frac{\alpha}{2}\, t^\frac{\gamma}{2}\,\|f\|_{\calC^\alpha}\|g\|_{\calC^\beta}\|h\|_{\calC^\gamma}
\end{split}
\end{equation*}
The continuity estimate \eqref{eq:eq1} comes from integrating with respect to time, taking into account the fact that $\alpha+\beta+\gamma > 0$.

\medskip

Let then estimate the regularity of  $g \mcR(f,h) - \mcS(f,g,h)$. For $x,y\in M$, with $d(x,y)\leq 1$, write
$$ 
\Big(g(x) \mcR(f,h) - \mcS(f,g,h)\Big)(x) -  \Big(g(y) \mcR(f,h) - \mcS(f,g,h)(y)\Big)(y) =: U+ V
$$
with $U$ defined by the formula
\begin{equation*}
\begin{split}
\bigintssss_0^{d(x,y)^2} \Bigg\{\mcP^1_t\Big(  \Gamma\big(\sqrt{t} \mcP^2_t f , \sqrt{t} \mcP_t^3 h\big) &\cdot \big(g(x) - P_t g\big)\Big)(x)  \\
&- \mcP^1_t\Big(  \Gamma\big(\sqrt{t} \mcP^2_t f , \sqrt{t} \mcP_t^3 h\big) \cdot \big(g(y) - P_t g\big)\Big)(y)\Bigg\}\,\frac{dt}{t},
\end{split}
\end{equation*}
and $V$ is defined by the formula
\begin{equation*}
\begin{split}
\bigintssss_{d(x,y)^2}^1 \Bigg\{\mcP^1_t\Big( \Gamma\big(\sqrt{t} \mcP^2_t f , \sqrt{t} \mcP_t^3 h\big)  & \cdot \big(g(x) - P_t g\big)\Big)(x)  \\ 
&- \mcP^1_t\Big( \Gamma\big(\sqrt{t} \mcP^2_t f , \sqrt{t} \mcP_t^3 h\big) \cdot \big(g(y) - P_t g\big)\Big)(y)\Bigg\}\,\frac{dt}{t}.
\end{split}
\end{equation*}
By repeating, the argument used in the proof of \eqref{eq:eq1}, we easily bound $U$ by the quantity
\begin{align*}
U & \lesssim \left(\int_0^{d(x,y)^2} t^{\frac{\alpha+\beta+\gamma}{2}} \, \frac{dt}{t}\right) \|f\|_{\calC^\alpha}\|g\|_{\calC^\beta}\|h\|_{\calC^\gamma} \\
   & \lesssim d(x,y)^\delta \|f\|_{\calC^\alpha}\|g\|_{\calC^\beta}\|h\|_{\calC^\gamma}.
\end{align*}
For the second part, we use the inequality
$$ 
|V| \leq A+B,
$$
with $A$ equal to
\begin{equation*}
\begin{split}
\left| \bigintssss_{d(x,y)^2}^1 \Bigg\{ \mcP^1_t\Big( \Gamma\Big(\sqrt{t} \mcP^2_t f , \sqrt{t} \mcP_t^3 h\Big) \right.&\cdot \big(g(x) - P_t g\big)\Big)(x) \\ 
&- \left.\mcP^1_t\Big(\Gamma\Big(\sqrt{t} \mcP^2_t f , \sqrt{t} \mcP_t^3 h\Big) \cdot \big(g(x) - P_t g\big)\Big)(y)\Bigg\}\,\frac{dt}{t}\right|
\end{split}
\end{equation*}
and 
$$
B:= \int_{d(x,y)^2}^1 \big|g(x)-g(y)\big|\cdot \Big| \mcP^1_t\Big(\Gamma\big(\sqrt{t} \mcP^2_t f , \sqrt{t} \mcP_t^3 h\big)\Big)(y)\Big| \, \frac{dt}{t}.
$$
The last quantity is bounded by
\begin{align*}
B & \lesssim d(x,y)^\beta\,\|g\|_{\calC^\beta} \int_{d(x,y)^2}^1 \Big\|\sqrt{t} \Gamma\big(\mcP^2_t f\big)\Big\|_\infty \Big\|\sqrt{t} \Gamma\big(\mcP_t^3 h\big)\Big\|_\infty  \, \frac{dt}{t} \\
& \lesssim d(x,y)^\beta\,\|f\|_{\calC^\alpha}\|g\|_{\calC^\beta}\|h\|_{\calC^\gamma}\int_{d(x,y)^2}^1 t^{\frac{\alpha+\gamma}{2}}  \, 
\frac{dt}{t} \\
& \lesssim d(x,y)^{\delta}\,\|f\|_{\calC^\alpha}\|g\|_{\calC^\beta}\|h\|_{\calC^\gamma}.
\end{align*}
For the quantity $A$, we use the Lipschitz regularity \eqref{Lipschitz} of the heat kernel to get the upper bound, 
{\small \begin{align*}
&\int_{d(x,y)^2}^1 \left\{ \int_M \frac{d(x,y)}{\sqrt{t}V(x,\sqrt{t})}\exp\left(-c \frac{d(x,z)^2}{t}\right)  \left|g(x) - P_t g(z)\right|\mu(dz) \right\} \\
&\hspace{1.2cm}\times \Big\|\sqrt{t} \Gamma\big(\mcP^2_t f\big) \Big\|_\infty \Big\|\sqrt{t} \Gamma\big(\mcP_t^3 h\big)\Big\|_\infty \, \frac{dt}{t} \\
 & \lesssim  \left\{ \int_{d(x,y)^2}^1 \int_M \frac{d(x,y)}{\sqrt{t}V\big(x,\sqrt{t}\big)}\,e^{-c \frac{d(x,z)^2}{t}}  \left(d(x,z)^\beta +t^{\frac{\beta}{2}}\right)\, \mu(dz)\,t^{\frac{\alpha+\gamma}{2}} \, \frac{dt}{t} \right\} \|f\|_{\calC^\alpha}\|g\|_{\calC^\beta}\|h\|_{\calC^\gamma} \\
 & \lesssim  \left( \int_{d(x,y)^2}^1 \frac{d(x,y)}{\sqrt{t}} t^{\frac{\alpha+\beta+\gamma}{2}} \, \frac{dt}{t} \right) \|f\|_{\calC^\alpha}\|g\|_{\calC^\beta}\|h\|_{\calC^\gamma} \\
  & \lesssim  d(x,y)^{\alpha+\beta+\gamma}\,\|f\|_{\calC^\alpha}\|g\|_{\calC^\beta}\|h\|_{\calC^\gamma},
\end{align*}}
\noindent where we have used the fact that $\alpha+\beta+\gamma\in(0,1)$. The combination of all the previous estimates yields
\begin{equation*}
\begin{split} 
\Big| \Big\{g\mcR(f,h) - &\mcS(f,g,h)\Big\}(x) -  \Big\{g\mcR(f,h) - \mcS(f,g,h)\Big\}(y)\Big| \\ 
&\leq |U|+A+B \lesssim d(x,y)^{\delta}  \|f\|_{\calC^\alpha}\|g\|_{\calC^\beta}\|h\|_{\calC^\gamma},
\end{split}
\end{equation*}
which concludes the proof of the continuity estimate \eqref{EqEstimateD1}.

\medskip

\noindent {\bf Step 2 -- proof of \eqref{EqEstimateD2}}. Given the collection $\Big(Q_r := Q_r^{(1)}\Big)_{r\in(0,1]}$ of operators, we need to prove that we have
\begin{equation}
\label{EqSecondBitCommutator}
\Big\| Q_r\Big(\mcR\big(\mcA(f,g),h\big) - \mcS(f,g,h)\Big)\Big\|_\infty \lesssim r^\frac{\delta}{2},
\end{equation}
for every $r\in(0,1]$, and where $\mcR\big(\mcA(f,g),h\big) \,- \, \mcS(f,g,h)$ is equal to 
\begin{equation}
\label{eq:eq4}
\int_0^1 \mcP^1_t \Gamma \left( \sqrt{t}\left\{\int_0^1 \mcP^2_t  \mcQ^1_s\left(\mcQ^2_s f \cdot \mcP_s^3 g\right) \, \frac{ds}{s} - P_t g \cdot\mcP^2_t f\right\} \,,\,\sqrt{t} \mcP_t^3 h\right) \, \frac{dt}{t}. 
\end{equation}
The notation may be confusing and we need to be careful: when $\Gamma$ acts on the product $P_t g \cdot\mcP^2_t f$, it is thought to act  only on the variable of $\mcP^2_tf$, with the variable of $P_t g$ frozen. We shall bound above the absolute value of the $\Gamma$ term in the integral, which is of the form $\Gamma(p,q)$, by $\Gamma(p)\Gamma(q)$ -- recall we write $\Gamma(p)$ for $\sqrt{\Gamma(p,p)}$. Set for that purpose
$$ 
A_{t}(f,g) := \sqrt{t} \Gamma\left(\int_0^1 \mcP^2_t  \mcQ^1_s\left(\mcQ^2_s f \cdot \mcP_s^3 g \right) \, \frac{ds}{s} - P_t g \mcP^2_t f\right).
$$
We have for almost every $x\in M$
\begin{align*}
 A_t (f,g)(x) & \leq \sqrt{t} \Gamma \mcP^2_t \left(\int_0^1 \mcQ^1_s\left(\mcQ^2_s f \cdot \mcP_s^3 g \right)(x) \, \frac{ds}{s} - P_t g(x) \cdot f\right)(x) \\
& \leq \int_0^1  \sqrt{t} \Gamma \mcP^2_t  \mcQ^1_s\Big(\mcQ^2_s f  \big(\mcP_s^3 g - P_t g(x)\big)\Big)(x) \, \frac{ds}{s} + \big|P_t g(x)\big|\, \sqrt{t}\Gamma[\mcP^2_t\mcP_1f](x) ,
\end{align*}
where we used the property 
$$ 
\int_0^1 \mcQ^1_s \mcQ^2_s\,\frac{ds}{s} = \textrm{Id} - \mcP_1,
$$
for some $\mcP_1$ operator.
As in Step 1, the fact that $\beta$ is positive, which implies the inequality
\begin{align*}
\Big| \mcP_s^3 g(y) - P_t g(x) \Big|  &\leq  \Big| \mcP_s^3 g(y) - g(y) \Big| +\big| g(y) -  g(x)\big| + \big| g(x) - P_t g(x) \big| \\
& \lesssim \left(s^{\frac{\beta}{2}} + t^{\frac{\beta}{2}} + d(x,y)^\beta\right)  \|g\|_{\calC^\beta} \lesssim \left( \max(s,t)^{\frac{\beta}{2}} + d(x,y)^\beta\right)  \|g\|_{\calC^\beta}.
\end{align*}
Moreover, it follows from Lemma \ref{lem:gauss} in Appendix A -- about the composition of Gaussian pointwise estimates, that the operator $\sqrt{t}\Gamma\Big(\mcP^2_t  Q^1_s \Big)$ has pointwise Gaussian estimates at the scale $\max(s,t)$ with an extra factor $\left(\frac{\min(s,t)}{\max(s,t)}\right)^{\frac{1}{2}}$; so if one sets $\tau:=\max(s,t)$, we have
{\small \begin{align*}
& \sqrt{t} \Gamma\mcP^2_t  \mcQ^1_s\Big(\mcQ^2_s f  \big(\mcP_s^3 g - P_t g(x)\big) \Big)(x) \\
& \qquad\lesssim \left(\frac{\min(s,t)}{\max(s,t)}\right)^{\frac{1}{2}} 
\left\{\int_M \frac{1}{V(x,\sqrt{\tau})}e^{-c \frac{d(x,y)^2}{\tau}} \left(\tau^{\frac{\beta}{2}} + d(x,y)^\beta\right) \, \mu(dy) \right\}\big\|Q^2_s f \big\|_\infty \|g\|_{\calC^\beta} \\
& \qquad \lesssim \left(\frac{\min(s,t)}{\max(s,t)}\right)^{\frac{1}{2}}  \tau^{\frac{\beta}{2}} \,s^{\frac{\alpha}{2}}\, \|f \|_{\calC^\alpha} \|g\|_{\calC^\beta}  \\
& \qquad \lesssim  \left(\frac{\min(s,t)}{\max(s,t)}\right)^{\frac{1}{2}} \,\max(s,t)^{\frac{\beta}{2}}  \,s^{\frac{\alpha}{2}} \,\|f \|_{\calC^\alpha} \|g\|_{\calC^\beta}.
\end{align*}}
Integrating in $s$, and taking into account the fact that $\alpha>-1$ and $\alpha+\beta< 1$, we obtain for $A_t(f,g)$ the estimate
\begin{align*}
\big\| A_t(f,g)\big\|_\infty & \lesssim \left\{\int_0^t \left(\frac{s}{t}\right)^{\frac{1}{2}}  t^{\frac{\beta}{2}}  s^{\frac{\alpha}{2}} \, \frac{ds}{s} + \int_t^1 \left(\frac{t}{s}\right)^{\frac{1}{2}}  s^{\frac{\beta}{2}}  s^{\frac{\alpha}{2}} \, \frac{ds}{s} +\sqrt{t} \right\}  \|f \|_{\calC^\alpha} \|g\|_{\calC^\beta} \\
 & \lesssim t^{\frac{\alpha+\beta}{2}}\,\|f \|_{\calC^\alpha} \|g\|_{\calC^\beta}.
\end{align*}
Observe that in the case where $\alpha+\beta\geq 1$, we get
\begin{align*}
\big\| A_t(f,g)\big\|_\infty \lesssim t^{\frac{1}{2}}\|f \|_{\calC^\alpha} \|g\|_{\calC^\beta},
\end{align*}
Coming back to identity \eqref{eq:eq4}, we have
$$ 
\Big| \mcR\big(\mcA(f,g),h\big) - \mcS(f,g,h)\Big| \leq \int_0^1 \mcP^1_t\Big( A_t(f,g)\cdot \sqrt{t} \Gamma\big(\mcP_t^4 h\big)\Big)\,\frac{dt}{t},
$$
and since $\alpha+\beta+\gamma>0$, it follows that
\begin{align*}
\Big\| \mcR(\mcA(f,g),h) - \mcS(f,g,h)\Big\|_\infty & \lesssim \left(\int_0^1 t^{\frac{\alpha+\beta+\gamma}{2}} \, \frac{dt}{t} \right) \|f \|_{\calC^\alpha} \|g\|_{\calC^\beta} \|h\|_{\calC^\gamma} \\
 & \lesssim \|f \|_{\calC^\alpha} \|g\|_{\calC^\beta} \|h\|_{\calC^\gamma}.
 \end{align*}
Moreover, taking into account that we have $Q_r^{(1)} \mcP_t^1 = \frac{r}{t} \mcP_r^1 \mcQ_t^{(1)}$  for $t\geq r$, and $\alpha+\beta+\gamma<1$, we see that the estimate \eqref{EqSecondBitCommutator} holds true
\begin{align*}
\Big\|Q_r^{(1)}\Big(\mcR\big(\mcA&(f,g),h\big) - \mcS(f,g,h) \Big)\Big\|_\infty \\
& \lesssim \int_0^r \big\| A_t(f,g)\big\|_\infty \Big\|\sqrt{t} \Gamma\big(\mcP_t^4 h\big)\Big\|_\infty \, \frac{dt}{t} +  \int_r^1 \frac{r}{t}
\,\big\| A_t(f,g)\big\|_\infty \Big\|\sqrt{t} \Gamma\big(\mcP_t^4 h\big)\Big\|_\infty \, \frac{dt}{t}  \\
& \lesssim \left(\int_0^r t^{\frac{\alpha+\beta+\gamma}{2}} \, \frac{dt}{t} + \int_r^1 r t^{\frac{\alpha+\beta+\gamma-2}{2}} \, \frac{dt}{t} \right) \| f\|_{\calC^\alpha} \|g\|_{\calC^\beta} \|h\|_{\calC^\gamma} \\
& \lesssim r^{\frac{\delta}{2}}\,\|f \|_{\calC^\alpha} \|g\|_{\calC^\beta} \|h\|_{\calC^\gamma}.
\end{align*}
\end{Dem}

\subsection[\hspace{-1cm} {\sf Paralinearization and composition estimates}]{Paralinearization and composition estimates}
\label{SubsectionParalinearization}

Two ingredients are needed to turn the machinery of paraproducts into an efficient tool. To understand how nonlinear functions act on H\"older functions $\calC^\alpha$, with $0<\alpha<1$, and to understand how one can compose two paraproducts. The first point is the object of the following analogue of Bony's classical result on paralinearization \cite{Bony}, while the second point is dealt with by theorem \ref{thm:composition} below.

\medskip

\begin{theorem}  \label{thm:paralinearization}    {\sf
Let fix an integer $b\geq2, \alpha\in(0,1)$, and consider a nonlinearity $\textrm{F}\in C_b^{3}$. Then for every $f\in \calC^\alpha$, we have $\textrm{F}(f)\in\calC^\alpha$ and
$$ 
R_\textrm{F}(f) := \textrm{F}(f) - \Pi^{(b)}_{\textrm{F}'(f)}(f) \in \calC^{2\alpha}.
$$
More precisely
$$ 
\Big\|\textrm{F}(f) - \Pi^{(b)}_{\textrm{F}'(f)}(f)\Big\|_{\calC^{2\alpha}} \lesssim \|\textrm{F}\|_{C^3_b} \Big(1+\|f\|_{\calC^\alpha}^2\Big).
$$
If $F\in C^4_b$ then the remainder term $R_\textrm{F}(f)$ is locally Lipschitz with respect to $f$, in so far far as we have
$$ 
\big\|R_\textrm{F}(f)-R_\textrm{F}(g)\big\|_{\calC^{2\alpha}} \lesssim \|\textrm{F}\|_{C^4_b} \big(1+\|f\|_{\calC^\alpha} + \|g\|_{\calC^\alpha} \big)^2\,\|f-g\|_{\calC^\alpha}.
$$  }
\end{theorem}

\medskip

\begin{Dem} 
First using the Leibniz rule for the operator $L$, we know that for $h\in \calC^\alpha$ then
$$ 
L\big(\textrm{F}(h)\big) = \textrm{F}'(h) L(h) + \textrm{F}''(h) \Gamma(h)^2.
$$
Now, since the semigroup is continuous at $t=0$, we have
\begin{align*}
\textrm{F}(f) & = \lim_{t\to 0}\,P_t^{(b)} \textrm{F}\big(P_t^{(b)} f\big),
\end{align*}
so we can write
\begin{align*}
\textrm{F}(f)  & = -\int_0^1 \frac{d}{dt} P_t^{(b)} \textrm{F}\big(P_t^{(b)} f\big) \, dt + P_1^{(b)} \textrm{F}\big(P_1^{(b)} f\big) \\
        & = \frac{1}{\gamma_b}  \int_0^1  \Big\{Q_t^{(b)}\big(\textrm{F}(P_t^{(b)} f)\big) + P_t^{(b)}\Big(Q_t^{(b)} f \cdot  \textrm{F}'\big(P_t^{(b)} f\big)\Big) \Big\}\, \frac{dt}{t}+  P_1^{(b)} \textrm{F}\big(P_1^{(b)} f\big).
\end{align*}
Using the relation $Q_t^{(b)}=Q_t^{(b-1)}(tL)$, together  with the chain rule
$$ 
L \Big(\textrm{F}\big(P_t^{(b)} f\big)\Big) = \textrm{F}'\big(P_t^{(b)}f\big) LP_t^{(b)}f + \textrm{F}''\big(P_t^{(b)}f\big) \Gamma\big(P_t^{(b)}f\big)^2,
$$
we get
$$  
Q_t^{(b)}\Big(\textrm{F}\big(P_t^{(b)} f\big)\Big) = Q_t^{(b-1)}\Big( (tL)P_t^{(b)}f \cdot \textrm{F}'\big(P_t^{(b)}\big)\Big) +  Q_t^{(b-1)}\Big(\textrm{F}''\big(P_t^{(b)}f\big) \cdot t\Gamma\big(P_t^{(b)}f\big)^2\Big).
$$
Note here the identity
\begin{align*}
P_t^{(b)}\Big(Q_t^{(b)} f \cdot  \textrm{F}'\big(P_t^{(b)} f\big)\Big) = & (tL) P_t^{(b)}\Big(Q_t^{(b-1)} f \cdot  \textrm{F}'\big(P_t^{(b)} f\big)\Big) - P_t^{(b)}\Big(Q_t^{(b-1)} f \cdot  tL \textrm{F}'\big(P_t^{(b)} f\big)\Big) \\
& - 2P_t^{(b)} t\Gamma\Big(Q_t^{(b-1)} f ,  \textrm{F}'\big(P_t^{(b)} f\big)\Big).
\end{align*}
So we have 
$$ 
\textrm{F}(f) - \Pi^{(b)}_{\textrm{F}'(f)}(f) =: (a) + (b) + (c) + (d) + (e) + (f)
$$
with
\begin{align*} 
(a) & := P_1^{(b)} \textrm{F}\big(P_1^{(b)} f\big), \\
(b) & := \frac{1}{\gamma_b} \int_0^1  Q_t^{(b-1)}\Big\{(tL)P_t^{(b)}f \cdot \Big(\textrm{F}'\big(P_t^{(b)} f\big) - P_t^{(b)}\big(\textrm{F}'(f)\big)\Big)\Big\}  \, \frac{dt}{t} \\   
(c) & := \frac{1}{\gamma_b} \int_0^1  (tL) P_t^{(b)}\Big\{ Q_t^{(b-1)} f \cdot  \Big(\textrm{F}'\big(P_t^{(b)} f\big) - P_t^{(b)}\big(\textrm{F}'(f)\big)\Big)\Big\}\, \frac{dt}{t}\\   
(d) & := \frac{1}{\gamma_b} \int_0^1  Q_t^{(b-1)}\Big\{\textrm{F}''\big(P_t^{(b)}f\big) \cdot t\Gamma\big(P_t^{(b)}f\big)^2\Big\} \, \frac{dt}{t} \\   
(e) & := -\frac{1}{\gamma_b} \int_0^1  P_t^{(b)}\Big\{Q_t^{(b-1)} f \cdot  tL \textrm{F}'\big(P_t^{(b)} f\big)\Big\} \, \frac{dt}{t}\\   
(f) & := -\frac{2}{\gamma_b} \int_0^1  P_t^{(b)} t\Gamma\Big(Q_t^{(b-1)} f ,  \textrm{F}'\big(P_t^{(b)} f\big) \Big)\, \frac{dt}{t}.
\end{align*}
We are now going to control each of these terms in the H\"older space $\calC^{2\alpha}$.

\medskip

\noindent {\bf Step 1 -- term $(a)$.} Since $f\in \calC^\alpha$, we know that $P_1^{(b)} f\in L^\infty$, so $\textrm{F}\big(P_1^{(b)} f\big)$ is also bounded. From Proposition \ref{prop:low}, we get
$$ 
\big\|(a)\big\|_{\calC^{2\alpha}} \lesssim \Big\|\textrm{F}\big(P_1^{(b)} f\big)\Big\|_\infty \lesssim \|f\|_{\calC^\alpha}.
$$

\medskip

\noindent {\bf Step 2 -- terms $(b), \,(c)$.} The following quantity appears in these two terms
\begin{align}
\Big| \textrm{F}'\big(P_t^{(b)} f\big) - P_t^{(b)}\big(\textrm{F}'(f)\big) \Big|  &\lesssim  \Big\| \textrm{F}'(f) - \textrm{F}'\big(P_t^{(b)} f\big) \Big\|_\infty + \Big\| F'(f) - P_t^{(b)}\big(\textrm{F}'(f)\big)\Big\|_\infty    \nonumber \\
 & \lesssim \big\|\textrm{F}''\big\|_\infty \, \Big\| f - P_t^{(b)}f \Big\|_\infty + \Big\| \textrm{F}'(f) - P_t^{(b)}\big(\textrm{F}'(f)\big)\Big\|_\infty    \nonumber \\
 & \lesssim \big\|\textrm{F}''\big\|_\infty \, \int_0^t \big\|Q_s^{(b)} f\big\|_\infty \, \frac{ds}{s} + \int_0^t \Big\|Q_s^{(b)}\big(\textrm{F}'(f)\big)\Big\|_\infty \, \frac{ds}{s}    \nonumber \\
 & \lesssim \big\|\textrm{F}''\big\|_\infty \, \left(\int_0^t s^{\frac{\alpha}{2}} \, \frac{ds}{s}\right) \|f\|_{\calC^\alpha} + \left(\int_0^t s^{\frac{\alpha}{2}} \, \frac{ds}{s}\right) \big\|\textrm{F}'(f)\big\|_{\calC^\alpha} \nonumber \\
 & \lesssim t^{\frac{\alpha}{2}}\,\big\|\textrm{F}''\big\|_\infty \, \|f\|_{\calC^\alpha}; \label{eq:qi}
\end{align}
we used along the way the characterization of H\"older space, for $0\alpha<1$, given by Proposition \ref{prop:caraholder}, to see that
$$ 
\big\|\textrm{F}'(f)\big\|_{\calC^\alpha} \lesssim \big\|\textrm{F}''\big\|_\infty\,\|f\|_{\calC^\alpha}.
$$
Using this estimate \eqref{eq:qi}, we deduce the following bound. Uniformly for every $s\in(0,1)$, we have
\begin{align*}
\Big\|Q_s^{(1)} (b)\Big\|_\infty & \lesssim \left(\int_0^1  \Big\|Q_s^{(1)}Q_t^{(b-1)}\Big\|_{\infty\to \infty} \,t^{\frac{\alpha}{2}}  \left\|\textrm{F}'(P_t^{(b)} f)-P_t^{(b)}\big(\textrm{F}'(f)\big)\right\|_\infty  \, \frac{dt}{t}\right)\|f\|_{\calC^\alpha}  \\   
   & \lesssim \left(\int_0^s t^{\alpha} \, \frac{dt}{t} + \int_s^1 \frac{s}{t} t^\alpha \, \frac{dt}{t}\right) \big\|\textrm{F}''\big\|_\infty \|f\|_{\calC^\alpha}^2 \\
   & \lesssim t^\alpha \big\|\textrm{F}''\big\|_\infty \|f\|_{\calC^\alpha}^2,
\end{align*}
where we used that $\big\|Q_s^{(1)}Q_t^{(b-1)}\big\|_{L^\infty\to L^\infty} \lesssim \frac{\min(s,t)}{\max(s,t)}$. That yields
$$
 \big\|(b)\big\|_{\calC^{2\alpha}} \lesssim \big\|F''\big\|_\infty \|f\|_{\calC^\alpha}^2,
 $$
and a similar inequality holds also for the third term $(c)$.

\medskip

\noindent {\bf Step 3 -- terms $(d),(e)$ and $(f)$.} We quickly sketch the boundedness of each of these three terms. Using Proposition \ref{prop:gradient}, we get a bound uniform in $s\in(0,1)$, of the form
\begin{align*}
 \big\|Q_s^{(1)} (d)\big\|_\infty & \lesssim \left(\int_0^1  \Big\|Q_s^{(1)}Q_t^{(b-1)}\Big\|_{\infty \to \infty} \Big\| \sqrt{t}\Gamma\big(P_t^{(b)}f\big)\Big\|_\infty ^2  \, \frac{dt}{t}\right) \big\|\textrm{F}''\big\|_\infty \\   
   & \lesssim \left(\int_0^s t^{\alpha} \, \frac{dt}{t} + \int_s^1 \frac{s}{t} t^\alpha \, \frac{dt}{t}\right) \big\|\textrm{F}''\big\|_\infty \|f\|_{\calC^\alpha}^2 \\
   & \lesssim t^\alpha \big\|\textrm{F}''\big\|_\infty \|f\|_{\calC^\alpha}^2.
\end{align*}
Similarly
\begin{align*}
 \big\|Q_s^{(1)} (f)\big\|_\infty & \lesssim \left(\int_0^1  \Big\|Q_s^{(1)}P_t^{(b)}\Big\|_{\infty \to \infty} \Big\| \sqrt{t}\Gamma\big(Q_t^{(b-1)}f\big)\Big\|_\infty \Big\| \sqrt{t}\Gamma\big(F'(P_t^{(b)}f)\big)\Big\|_\infty  \, \frac{dt}{t}\right)  \\   
   & \lesssim \left(\int_0^s t^{\alpha} \, \frac{dt}{t} + \int_s^1 \frac{s}{t} t^\alpha \, \frac{dt}{t}\right) \big\|\textrm{F}''\big\|_\infty \|f\|_{\calC^\alpha}^2 \\
   & \lesssim t^\alpha \big\|\textrm{F}''\big\|_\infty \|f\|_{\calC^\alpha}^2,
\end{align*}
where we used the Leibniz rule 
$$
\Big\| \sqrt{t}\Gamma\Big(\textrm{F}'\big(P_t^{(b)}f\big)\Big)\Big\|_\infty \leq \big\|\textrm{F}''\big\|_\infty \Big\| \sqrt{t}\Gamma\big((P_t^{(b)}f)\big)\Big\|_\infty\lesssim t^{\frac{\alpha}{2}}\,\big\|\textrm{F}''\big\|_\infty \|f\|_{\calC^\alpha}.
$$
For the remaining last term $(e)$, we can still using the Leibniz rule and get
$$
\Big\| tL \textrm{F}'\big(P_t^{(b)} f\big)\Big\|_\infty \leq \Big(\big\|\textrm{F}''\big\|_\infty+\big\|\textrm{F}'''\big\|_\infty\Big) \left(\Big\|(tL)P_t^{(b)} f\Big\|_\infty + \Big\|\sqrt{t}\Gamma\big(P_t^{(b)}f\big)\Big\|_\infty^2\right)
 $$
which then yields
\begin{align*}
 \|Q_s^{(1)} (e)\|_\infty & \lesssim \int_0^1  \Big\|Q_s^{(1)}P_t^{(b)}\Big\|_{\infty \to \infty} \, \Big\| Q_t^{(b-1)}f\Big\|_\infty \, \Big\| tL\, \textrm{F}'\big(P_t^{(b)} f\big)\Big\|_\infty  \, \frac{dt}{t}  \\   
   & \lesssim \left(\int_0^s t^{\alpha} \, \frac{dt}{t} + \int_s^1 \frac{s}{t} t^\alpha \, \frac{dt}{t}\right) \|\textrm{F}\|_{C^3_b} \|f\|_{\calC^\alpha}^2 \\
   & \lesssim t^\alpha \, \|\textrm{F}\|_{C^3_b} \|f\|_{\calC^\alpha}^2.
\end{align*}
By combining the previous estimates, we conclude that we have
$$ 
\big\|(d)\big\|_{\calC^{2\alpha}} + \big\|(e)\big\|_{\calC^{2\alpha}} + \big\|(f)\big\|_{\calC^{2\alpha}} \lesssim \|\textrm{F}\|_{C^3_b} \, \|f\|_{\calC^\alpha}\,(1+\|f\|_{\calC^\alpha}),
$$
which ends the proof of the estimate of the remainder. The Lipschitz regularity of the remainder term is proved by very similar arguments which we leave to the reader.
\end{Dem}

\medskip

Let us now examine the composition of two paraproducts. Note that for $u\in \calC^\alpha$ and $v\in \calC^\beta$, with $\alpha\in(0,1)$, $\beta\in(0,\alpha]$, we have $uv\in\calC^\beta$.

\medskip

\begin{theorem}  \label{thm:composition}    {\sf
Fix an integer $b\geq2$, $\alpha\in(0,1)$, $\beta\in(0,\alpha]$ and consider $u\in \calC^\alpha$ and $v\in \calC^\beta$. Then for every $f\in \calC^\alpha$, we have
$$
\Pi^{(b)}_{u} \Big( \Pi^{(b)}_{v}(f) \Big) - \Pi^{(b)}_{uv}(f) \in \calC^{\alpha+\beta}
$$
with
$$ 
\left\| \Pi^{(b)}_{u}\Big( \Pi^{(b)}_{v}(f) \Big) - \Pi^{(b)}_{uv}(f) \right\|_{\calC^{\alpha+\beta}} \lesssim \|f\|_{\calC^\alpha} \, \|u\|_{\calC^\alpha} \|v\|_{\calC^\beta}.
$$   }
\end{theorem}

\medskip

\begin{Dem} 
We leave a detailed proof to the reader and we just sketch it, since it is similar and easier than the proof of Theorem \ref{thm:paralinearization}. Following Proposition \ref{prop:para}, we know that the two terms $\Pi^{(b)}_{u} \Big( \Pi^{(b)}_{v}(f) \Big)$ and $\Pi^{(b)}_{uv}(f)$ belong to $\calC^\alpha$. The idea is to use the $\calC^\beta$-regularity of $v$ to gain the same regularity in the difference.

\ssk

Indeed, adopting the notations used above, the paraproduct $\Pi^{(b)}_g(f)$ is given, up to a multiplicative constant, by two terms with the form
$$ 
I(f,g) = \int_0^1 \mcQ^1_t\Big(\mcQ^2_t f \cdot \mcP_t^1 g \Big) \, \frac{dt}{t},
$$
where in $\mcQ^1_t$ and $\mcQ^2_t$ we have at least a term $(tL)$ to the power $1$. Let us focus on this form. Then we have
$$  
I\big(I(f,v),u\big) =    \int_0^1\int_0^1  \mcQ^1_t\Big(\mcQ^2_t \mcQ^1_s\left(\mcQ^2_s f \cdot \mcP_s^1 v \right)  \cdot \mcP_t^1 u \Big) \, \frac{ds\,dt}{st}
$$
and
$$ 
I(f,vu) =   \int_0^1\int_0^1  \mcQ^1_t\Big(\mcQ^2_t \mcQ^1_s\left(\mcQ^2_s f\right)  \cdot\mcP_t^1 (uv) \Big) \, \frac{ds\,dt}{st},
$$
where we have used the normalization $\Pi^{(b)}_{1} = \textrm{Id}$, which means here that $ I(f,vu) = I\big(I(f,1),vu\big)$. Then using the $\calC^\beta$-regularity of $v$ and the fact that $\mcQ^i_t$ involves at least a power $1$ of $(tL)$, one can check that uniformly in $s,t\in(0,1)^2$, we have
\begin{equation*}
\begin{split}
\Big\|\mcQ^1_t\Big(\mcQ^2_t \mcQ^1_s\left(\mcQ^2_s f \cdot \mcP_s^1 v \right)  \cdot P_t^1 u \Big) - \mcQ^1_t &\Big(\mcQ^2_t \mcQ^1_s\left(\mcQ^2_s f\right)  \cdot \mcP_t^1 (uv) \Big)\Big\|_{\infty} \\
&\lesssim \frac{\min(s,t)}{\max(s,t)} \,s^{\frac{\alpha}{2}} \,(s+t)^{\frac{\beta}{2}}\, \|f\|_{\calC^\alpha} \|v\|_{\calC^\beta} \|u\|_{\calC^\alpha}.
\end{split}
\end{equation*}
So integrating in $s\in(0,1)$ yields for $\alpha+\beta<2$
\begin{equation*}
\begin{split}
\int_0^1 \Big\|\mcQ^1_t\Big(\mcQ^2_t \mcQ^1_s\left(\mcQ^2_s f \cdot \mcP_s^1 g \right)  \cdot \mcP_t^1 u \Big) - \mcQ^1_t &\Big(\mcQ^2_t \mcQ^1_s\left(\mcQ^2_s f\right)  \cdot \mcP_t^1 (uv) \Big)\Big\|_{\infty} \, \frac{ds}{s} \\
&\lesssim t^{\frac{\alpha+\beta}{2}} \, \|f\|_{\calC^\alpha} \|v\|_{\calC^\beta} \|u\|_{\calC^\alpha}.
\end{split}
\end{equation*}
Then as previously we check that for every $\tau\in(0,1)$ we have
\begin{align*}
& \Big\| \mcQ^{1}_{\tau}\Big(I\big(I(f,v),u\big) - I(f,vu)\Big) \Big\|_{\infty} \\
& \qquad \lesssim \left(\int_0^1 \frac{\min(\tau,t)}{\max(\tau,t)}\, t^{(\alpha+\beta)/2} \, \frac{dt}{t} \right) \|f\|_{\calC^\alpha} \|v\|_{\calC^\beta} \|u\|_{\calC^\alpha} \\
& \qquad  \lesssim \tau^{\frac{\alpha+\beta}{2}} \, \|f\|_{\calC^\alpha} \|v\|_{\calC^\beta} \|u\|_{\calC^\alpha},
\end{align*}
since $\alpha+\beta<2$. That allows us to conclude that 
\begin{align*}
 \Big\| I\big(I(f,v),u\big) - I(f,vu) \Big\|_{\calC^{\alpha+\beta} } \lesssim  \|f\|_{\calC^\alpha} \|v\|_{\calC^\beta} \|u\|_{\calC^\alpha}.
\end{align*}
\end{Dem}

\medskip

\subsection[\hspace{-1cm} {\sf Schauder estimates}]{Schauder estimates}
\label{SubsectionSchauder}

Proposition \ref{prop:schauder1} gives an elementary proof in our setting of a Schauder-type estimate about the regularizing character of the convolution operation with the operators $P^{(b)}_s$. Its paracontrolled analogue, given in section \ref{SubsectionParacontrolledSchauder} provides a crucial ingredient in the study of parabolic singular PDEs, from the point of view of paracontolled distributions.

\medskip

\begin{definition}  \label{def:space}    {\sf
For $\alpha\in(0,2)$ and $T>0$, we set 
$$ 
C_T \calC^\alpha := \left\{ f\in L^\infty(\SSS),\  \|f\|_{C_T \calC^\alpha}:=\sup_{t\in[0,T]} \big\|f(t)\big\|_{\calC^\alpha}<\infty\right\}
$$
and 
$$ 
C_T^{\frac{\alpha}{2}} L^\infty:= \left\{f\in L^\infty([0,T] \times M),\ \|f\|_{C_T^{\frac{\alpha}{2}} L^\infty} :=\sup_{\genfrac{}{}{0pt}{}{s\neq t}{0\leq s,t \leq T}} \frac{\big\|f(t)-f(s)\big\|_\infty}{|t-s|^{\frac{\alpha}{2}}} <\infty \right\}.
$$
We then define the space 
$$
\mcE^{(\alpha)}_{[0,T]} := C_T\calC^\alpha \cap C_T^{\frac{\alpha}{2}} L^\infty.
$$   }
\end{definition}

\medskip

The resolution operator $\mcR$ is formally defined by the formula
$$
\mcR(v)_t := \int_0^t P^{(b)}_{t-s} v(s)\, ds;
$$
it implicitly depends on the parameter $b$.

\medskip

\begin{proposition}  \label{prop:schauder1}    {\sf
Consider an integer $b\geq 0$ and a regularity exponent $\beta\in{\mathbb R}$. For every positive finite time horizon $T$, and every $v\in C_T\calC^\beta$, then $\mcR(v)$ belongs to $C_T\calC^{\beta+2}$, and 
$$ 
\big\|\mcR(v)_t\big\|_{\calC^{\beta+2}} \lesssim (1+T) \sup_{s\in[0,t]} \big\|v(s)\big\|_{\calC^\beta},
$$
for every $t \in [0,T]$. Moreover if $-2<\beta<0$ then we also have
$$ 
\big\|\mcR(v)\big\|_{C_T^{\frac{\beta+2}{2}}L^\infty} \lesssim \|v\|_{C_T \calC^\beta}.
$$   }
\end{proposition}

\medskip

\begin{Dem}
We consider another integer $c \geq \frac{|\beta|}{2} + 1$ and a parameter $\tau\in(0,1]$. Then
$$ 
Q_\tau^{(c)}\big(\mcR(v)_t\big) = \int_0^t Q_\tau^{(c)} P^{(b)}_{t-s} v(s) \, ds.
$$
We have
\begin{equation}
\Big\|Q_\tau^{(c)} P^{(b)}_{t-s} v(s)\Big\|_{\infty} \lesssim \left(\frac{\tau}{\tau + t-s}\right)^{c} (\tau + t-s)^{\frac{\beta}{2}} \|v(s)\|_{\calC^\beta}.
\label{eq:eq:} \end{equation}
Indeed, if $t-s\leq \tau$ then we only use that $Q_\tau^{(c)}$ and $P^{(b)}_{t-s}$ commute with the $L^\infty$-boundedness of $P^{(b)}_{t-s}$ to have
$$
\Big\|Q_\tau^{(c)} P^{(b)}_{t-s} v(s)\Big\|_{\infty} \lesssim \Big\|Q_\tau^{(c)} v(s)\Big\|_{\infty} \lesssim \tau^{\frac{\beta}{2}} \|v(s)\|_{\calC^\beta}.
$$
In the other hand, if $t-s\geq \tau$ then we use that
$$ Q_\tau^{(c)}P_{t-s}^{(b)}=\left(\frac{\tau}{t-s}\right)^c e^{-\tau L} \big((t-s)\big)^cP_{t-s}^{(b)}$$
and we conclude similarly with the $L^\infty$-boundedness of $e^{-\tau L}$ and the property that 
\begin{align*}
\left\|\big((t-s)\big)^cP_{t-s}^{(b)} v(s)\right\|_{\infty} & \lesssim \int_{t-s}^\infty \left\|\big((t-s)\big)^c Q_r ^{(b)} v(s)\right\|_{\infty} \,\frac{dr}{r} \\
& \lesssim \int_{t-s}^\infty \left(\frac{t-s}{r}\right)^c\left\| Q_r^{(b+c)} v(s)\right\|_{\infty} \,\frac{dr}{r} \\
& \lesssim \int_{t-s}^\infty \left(\frac{t-s}{r}\right)^c r^{\beta/2} \left\| v(s)\right\|_{\calC^\beta} \,\frac{dr}{r} \\
& \lesssim (t-s)^{\beta/2}\left\| v(s)\right\|_{\calC^\beta}.
\end{align*}
That concludes the proof of \eqref{eq:eq:}. So by integrating, it comes
\begin{align*}
\Big\| Q_\tau^{(c)}\big(\mcR(v)_t\big)\Big\|_\infty & \lesssim \left\{\int_0^t \left(\frac{\tau}{\tau+t-s}\right)^{c} (\tau + t-s)^{\frac{\beta}{2}} \, ds\right\} \sup_{s\in[0,t]}  \big\|v(s)\big\|_{\calC^\beta} \\
& \lesssim \tau^{\frac{\beta}{2}+1} \sup_{s\in[0,t]} \big\|v(s)\big\|_{\calC^\beta}.
\end{align*}
This holds uniformly in $\tau\in(0,1]$ and so one concludes the proof of the first statement with the global inequality
$$ 
\big\|\mcR(v)_t\big\|_\infty \lesssim \left\{\int_0^t (\cdots) ds\right\} \|v\|_{C_t\calC^\beta} \lesssim T \|v\|_{C_T\calC^\beta}.
$$
For the second statement, we note that for $s<t\leq T$ we have
\begin{align*}
\mcR(v)_t - \mcR(v)_s  &= \Big(P_{t-s}^{(b)}-\textrm{Id}\Big)\big(\mcR(v)_s\big) + \int_s^t P_{t-r}^{(b)}\big(v(r)\big) \, dr \\ 
							           &= \frac{1}{\gamma_a}\,\int_0^{t-s} Q^{(a)}_{r} \mcR(v)_s \, \frac{dr}{r} + \int_s^t P_{t-r}^{(b)}\big(v(r)\big) \, dr.
 \end{align*}
We have
\begin{align*} 
\left\| \int_0^{t-s} Q^{(a)}_{r} \mcR(v)_s \, \frac{dr}{r} \right\|_{\infty} & \lesssim \left(\int_0^{t-s} r^{\frac{\beta}{2}+1} \, \frac{dr}{r} \right) \big\|\mcR(v)_s\big\|_{\calC^{\beta+2}} \\ 
 & \lesssim (t-s)^{\frac{\beta}{2}+1} \big\|\mcR(v)_s\big\|_{\calC^{\beta+2}}
\end{align*}
and since $\beta<0$, we also have
\begin{align*}
\left\| \int_s^t P_{t-r}^{(b)}\big(v(r)\big) \, dr \right\|_{L^\infty} & \lesssim \int_s^t \left(\int_{t-r}^1 \Big\|Q_{\tau}^{(b)} v(r)\Big\|_\infty \, \frac{d\tau}{\tau} + \Big\|P_{1}^{(b)}\big(v(r)\big)\Big\|_{\infty}\right) dr \\
& \lesssim \int_s^t  \left(\big\|v(r)\big\|_{\calC^\beta} \int_{t-r}^1 \tau^{\frac{\beta}{2}} \, \frac{d\tau}{\tau} + \Big\|P_{1}^{(b)}\big(v(r)\big)\Big\|_{\infty}\right) dr \\
& \lesssim (t-s)^{\frac{\beta}{2}+1} \sup_{r\in[0,t]} \big\|v(r)\big\|_{\calC^\beta},
\end{align*}
where we used $\frac{\beta}{2}+1\in (0,1)$.
\end{Dem}

\medskip

\begin{coro} \label{cor:schauder1}    {\sf
For a fixed integer $b\geq 0$ and $\alpha\in(0,2)$, we gave
$$
\big\|\mcR(v)\big\|_{\mcE^{(\alpha)}_{[0,T]}} \lesssim (1+T)\,\|v\|_{C_T \calC^{\alpha-2}},
$$
uniformly in $T>0$.   }
\end{coro}

\medskip

\begin{rem} \label{rem:T} \textit{\textrm{
Observe that in Proposition \ref{prop:schauder1} the weight $(1+T)$ can be weakened, up to a little loss on the regularity exponent. Indeed, the exact same proof allows us to prove
$$ 
\big\|\mcR(v)_t\big\|_{\calC^{\beta+2-2\epsilon}} \lesssim T^\epsilon \sup_{s\in[0,t]}\, \big\|v(s)\big\|_{\calC^\beta}
$$
and 
$$
\big\|\mcR(v)\big\|_{C_T^{\frac{1}{2}(\beta+2-2\epsilon)}L^\infty} \lesssim T^\epsilon\, \|v\|_{C_T \calC^\beta},
$$
for any $\epsilon\in(0,1)$; so we have 
$$ 
\big\|\mcR(v)\big\|_{\mcE^{(\alpha)}_{[0,T]}} \lesssim T^\epsilon \,\|v\|_{C_T \calC^{\alpha-2+2\epsilon}}.
$$ 
We refer the reader to Proposition \ref{prop:schauder2} for a detailed proof of a more difficult statement, where we show how we can improve the bound $(1+T)$ up to a small loss on the regularity.   }}
\end{rem}

\bigskip

\section[\hspace{0.6cm} {\sf Paracontrolled calculus}]{Paracontrolled calculus} 
\label{sec:para}

The idea of paracontrolled calculus, such as introduced in \cite{GIP}, has its roots in Gubinelli's notion of controlled path \cite{GControlled}. The latter provides an alternative formulation of Lyons' rough paths theory \cite{Lyons98,LyonsStFlour} that offers a simple approach to the core of the theory, while rephrasing it in a very useful Banach setting. Let us have here a glimpse at this field, as a guide for what we shall be doing in this section and the next one. We refer the reader to \cite{FH} for a very nice and pedagogical introduction to the subject, assuming only here that she/he knows only the very definition of a (weak geometric) $\alpha$-H\"older rough path, for some $\frac{1}{3} < \alpha \leq \frac{1}{2}$; see also \cite{BM2Course}. Let just mention that these objects are nothing else than objects that play the role of the collection 
$$
\left(h_t-h_s\,,\,\int_s^t\int_s^r \dot h_u\otimes \dot h_r\,dudr\right)_{0\leq s\leq t\leq T}
$$
of the increments of an $\RR^\ell$-valued control $h$ in a controlled ordinary differential equation $\dot x_t = V_i(x_t)\,\dot h^i_t$, together with its second order iterated integral. These quantities are precisely what appears when making a second order Euler-Taylor expansion of the solution to the equation. Rough paths are such kind of objects for which $h$ is too irregular to make sense of the product $dh_u\otimes dh_r$; typical irregularity is $\alpha$-H\"older, with $\alpha<\frac{1}{2}$, like for Brownian motion.

\medskip

Assume we are given an $\RR^\ell$-valued (weak geometric) $\alpha$-H\"older rough path 
$$
{\bfX} = \big((X_{ts},\bbX_{ts})\big)_{0\leq s\leq t\leq T},
$$ 
with $X_{ts}\in\RR^\ell$ and $\bbX_{ts}\in \RR^\ell\otimes\RR^\ell$; recall that $X$ is $\alpha$-H\"older and $\bbX$ is $(2\alpha)$-H\"older.
Let also $\sigma = (V_1,\dots,V_\ell) \in C^3\big(\RR^d,\textrm{L}\big(\RR^\ell,\RR^d\big)\big)$ be given, with each column $V_i$ seen as a vector field on $\RR^d$. Following Lyons, an $\RR^d$-valued path $x_\bullet$ is said to solve the rough differential equation 
\begin{equation}
\label{EqRDE}
dx_t = \sigma(x_t)\,{\bfX}(dt)
\end{equation}
if one has, for every smooth observable $f$, the second order Euler-Taylor expansion
\begin{equation}
\label{EqRDEDefn}
f(x_t) - f(x_s) = X^i_{ts}(V_if)(x_s) + \bbX^{jk}_{ts}(V_jV_kf)(x_s) + O\big(|t-s|^a\big)
\end{equation}
for all $0\leq s\leq t\leq T$, for some constant $a>1$. Note that if $X_{ts} = h_t-h_s$, and $\bbX_{ts} = \int_s^t (h_r-h_s)\otimes dh_r$, for some $\RR^\ell$-valued $C^1$ control $h$, equation \eqref{EqRDEDefn} is nothing but a second order Taylor expansion for the solution to the controlled differential equation $\dot x_t = \sigma(x_t)\,\dot h_t$. Gubinelli's crucial remark in \cite{GControlled} was to notice that for a path $x_\bullet$ to satisfy equation \eqref{EqRDEDefn}, it needs to be controlled by $X$ in the sense that one has
\begin{equation}
\label{EqTaylorOneRDE}
x_t-x_s = x'_s\,X_{ts} + O\big(|t-s|^{2\alpha}\big),
\end{equation}
for some $\textrm{L}(\RR^\ell,\RR^d)$-valued $\alpha$-H\"older path $x'_\bullet$, here $x'_s = \sigma(x_s)$. This set of paths has a natural Banach topology. The point of this remark is that, somewhat conversely, if we are given an $\textrm{L}(\RR^\ell,\RR^d)$-valued $\alpha$-H\"older path $(z,z')$ controlled by $X$, then there exists a unique $\RR^d$-valued path $y_\bullet$ whose increments satisfy 
$$
y_t-y_s = z_s\,X_{ts} + z'_s\,\bbX_{ts} + O\big(|t-s|^a\big),
$$
for some exponent $a>1$. Note that we indeed need the full rough path $\bfX$ to define that path, and not just $X$. With a little bit of abuse, we write $\int_0^\bullet z_s\,{\bfX}(ds)$ for that path $y_\bullet$ -- this path depends not only on $z$ but rather on  $(z,z')$. Given an $\RR^d$-valued path $(x,x')$ controlled by $X$, and $\sigma$ sufficiently regular, the $\textrm{L}(\RR^\ell,\RR^d)$-valued path $z_s := \sigma(x_s)$ is controlled by $X$, with derivative $z'_s = \sigma'_{x_s}(x'_s)$. With $z_s=\sigma(x_s)$ and $x_\bullet$ satisfying the first order Euler-Taylor expansion \eqref{EqTaylorOneRDE}, the above second order Euler-Taylor expansion for $y$ is nothing else than the right hand side of \eqref{EqRDEDefn}, with the identity in the role of $f$. It can be proved that the rough integral $\int_0^\bullet z_s\,{\bfX}(ds)$ depends continuously on $(z,z')$ and $\bfX$ in the right topologies -- this is the main selling point of rough paths theory. So, for a path $x_\bullet$ to solve the rough differential equation \eqref{EqRDE}, it is necessary and sufficient that it satisfies
$$
x_t-x_s = \int_s^t\sigma(x_r)\,{\bfX}(dr),
$$
for all $0\leq s\leq t\leq T$, that is, $x_\bullet$ is a fixed point of the continuous map
$$
x_\bullet \mapsto \int_0^\bullet \sigma(x_r)\,{\bfX}(dr),
$$
from the space of paths controlled by $\bfX$ to itself. The well-posed character of equation \eqref{EqRDE} is then shown by proving that this map is a contraction if one works on a sufficiently small time interval.

\medskip

Our present setting will not differ much from the above description. We aim in the sequel at solving equations of the form 
$$
\big(\partial_t+\Delta\big) u = \textrm{F}(u)\,\zeta,
$$
for some distribution $\zeta$. Comparing this equation with the rough differential equation \eqref{EqRDE}, the role of the rough path $\bfX$ will be played in that setting by a rough distribution $\widehat{\zeta} = \big(\zeta,\zeta^{(2)}\big)$, with $\zeta$ in the role of $dX_t$ -- or $dh$, with $\zeta^{(2)}$ somehow in the role of $d\bbX_t$ -- or $dh\otimes dh$, and the operator $\big(\partial_t+\Delta\big)$ in the role of the time derivative operator $\frac{d}{dt}$. The elementary insight that the/a solution $u$ should behave at small space scales like the solution $Z$ to the equation $(\partial_t+L)Z = \zeta$, is turned into the definition of a distribution "controlled by" $Z$, such as given below, using the paraproduct as a means of comparison, for writing a first order Taylor expansion of $u$ similar to identity \eqref{EqTaylorOneRDE} -- compare this ansatz with the fact that a solution to a rough differential equation should be controlled by $X = \int dX$. The crucial point of this definition is that one can make sense of the product $\textrm{F}(u)\,\zeta$, in that controlled setting, which provides an analogue of the definition of the product $\sigma(x_s)\,{\bfX}(ds)$, given by the right hand side of identity \eqref{EqRDEDefn}  -- see theorem \ref{thm:prod}. To run formally the above argument, we shall need to see how controlled distributions are transformed by a nonlinear map; this is the content of theorem \ref{thm:nolinear} on paralinearisation. Some problems intrinsically linked with the multidimensional setting of the problem are dealt with in section \ref{SubsectionParacontrolledSchauder}, where a version of Schauder theorem is proved for paracontrolled distributions.

\medskip

A last look at section \ref{SubsectionParacontrolledCalculus} may provide a helpful guide for this section, before proceeding.

\bigskip

\subsection[\hspace{-1cm} {\sf Paracontrolled distributions}]{Paracontrolled distributions}
\label{SubsectionParaControlledDistributions}

We fix throughout that section and the next one an integer $b\geq 2$.

\medskip

\begin{definition*} {\sf
Let $\alpha\in (-2,1)$ and $\beta>0$ be given, together with a reference distribution $Z\in \calC^\alpha$. A pair of \textbf{distributions} $(f,f')\in \calC^\alpha \times \calC^\beta$ is said to be \textbf{paracontrolled by $Z$} if 
$$
(f,f')^\sharp:=f-\Pi^{(b)}_{f'}(X) \in \calC^{\alpha+\beta}.
$$ 
In such a case, we write $(f,f')\in {\sf C}^\beta_\alpha(Z)$ and define the norm
$$ 
\big\|(f,f')\big\|_{{\sf C}^\beta_\alpha}:= \big\|(f,f')^\sharp\big\|_{\calC^{\alpha+\beta}} + \|f'\|_{\calC^\beta}.
$$
If $\overline{Z}\in \calC^\alpha$ is a possibly different reference distribution and $(h,h')\in {\sf C}^\beta_\alpha(\overline{Z})$, we set
$$ 
d_{{\sf C}^\beta_\alpha} \big( (f,f'),(h,h') \big):= \|f'-h'\|_{\calC^\beta} + \big\| (f,f')^\sharp - (h,h')^\sharp \big\|_{\calC^{\alpha+\beta}}.
$$    }
\end{definition*}

\medskip

\noindent Note that this choice of norm allows to compare paracontrolled distributions associated with \textit{different} model distributions $X$ and $\overline{Z}$. Following the terminology of \cite{GIP}, the function $f'$ is called the \textbf{derivative of $f$}, and the term $(f,f')^\sharp$, the \textbf{remainder}; one should think of the decomposition 
$$ 
f = \Pi^{(b)}_{f'}(X) + (f,f')^\sharp
$$
as a kind of first order Taylor formula for $f$, in terms of regularity properties. The notion of derivative depends of course on which model distribution is used. As a first step towards completing the above program, the following statement gives an analogue in our setting of the right hand side of identity \eqref{EqRDEDefn} defining $\sigma(x_s)\,{\bfX}(ds)$ in the rough paths context. It is motivated by the following simple regularity analysis based on propositions \ref{prop:para} and \ref{prop:diag}, giving regularity conditions for the well-posed character of paraproducts and resonnant terms. Given $f\in\calC^\alpha$ and $\zeta\in\calC^\gamma$, with $\gamma<0<\alpha<1$, we have from Calderon's identity the formal identity
$$
f\zeta = \Pi^{(b)}_f(\zeta) + \Pi^{(b)}_\zeta(f) + \Pi^{(b)}(f,\zeta) + \Delta_{-1}(f,\zeta),
$$
where the only term that is potentially undefined is the diagonal term $\Pi^{(b)}(f,\zeta)$. If however, $f$ is controlled by $Z$, with derivative $f'\in\calC^\beta$, we can write
$$
\Pi^{(b)}(f,\zeta) = \Pi^{(b)}\Big(\Pi^{(b)}_{f'}(Z),\zeta\Big) + \Pi^{(b)}\big((f,f')^\sharp,\zeta\big),
$$
with $\Pi^{(b)}\Big((f,f')^\sharp , \zeta\Big)$ well-defined if $(\alpha+\beta)+\gamma > 0$. So, writing
$$
\Pi^{(b)}\Big(\Pi^{(b)}_{f'}(Z) , \zeta\Big) = C(Z,f',g) + f'\,\Pi^{(b)}(Z,\zeta),
$$
we finally see that, in the paracontrolled setting, the only undefined term in the above a priori decomposition of $fg$ is the term $\Pi^{(b)}(Z,\zeta)$. The following theorem turns that elementary regularity analysis into a constructive recipe for defining $f\zeta$.

\medskip

\begin{theorem}[\textbf{Product operation in a paracontrolled space}]  \label{thm:prod}    {\sf
Let $\gamma<0<\beta<\alpha<1$ be regularity exponents such that 
$$
\alpha+\gamma<0,   \qquad   \alpha+\beta+\gamma \in (0,1).
$$ 
Let $(Z_n)_{n\geq 0}$ be a sequence of smooth functions converging to $Z$ in $\calC^\alpha$, and $(\zeta_n)_{n\geq 0}$ be a sequence of smooth functions converging to $\zeta$ in $\calC^\gamma$. Assume that $\Pi^{(b)}(Z_n,\zeta_n)$ converges in $\calC^{\alpha+\gamma}$ to some limit distribution $\zeta^{(2)}$; write $\widehat{\zeta}$ for the pair $\big(\zeta,\zeta^{(2)}\big)$ -- call it an \textbf{enhanced, or rough, distribution}.
\begin{itemize}
   \item Given a paracontrolled function $(f,f')\in {\sf C}_\alpha^\beta(Z)$, the formula
\begin{equation}
\label{EqDefnProduct}
(f,f')\cdot\widehat{\zeta} := \Pi^{(b)}_f(\zeta) + \Pi^{(b)}_\zeta(f) + \Pi^{(b)}\big( (f,f')^\sharp , \zeta\big) + C(Z,f',\zeta) + f'\,\zeta^{(2)}
\end{equation}
defines an element of $\calC^\gamma$ which satisfies the estimate
\begin{equation} 
\label{eq:DD}
\big\| (f,f')\cdot\widehat{\zeta} - \Pi^{b}_f(\zeta) \big\|_{\calC^{\alpha+\gamma}} \lesssim \big\|(f,f')\big\|_{{\sf C}^\beta_\alpha}\,\Big\{\|\zeta\|_{\calC^\gamma} + \|Z\|_{\calC^\alpha}\|\zeta\|_{\calC^\gamma} + \big\|\zeta^{(2)}\big\|_{\calC^{\alpha+\gamma}}\Big\}; 
\end{equation}
so $\big( (f,f')\cdot\zeta , f\big)\in {\sf C}^\alpha_\gamma(\zeta)$.   \vspace{0.15cm}

   \item Furthermore, this operation is locally Lipschitz in the sense that if $\big(\overline{Z},\overline{\zeta},\overline{\zeta}^{(2)}\big)$ is another set of objects similar to $\big(Z,\zeta,\zeta^{(2)}\big)$, and if $(h,h')\in{\sf C}_\alpha^\beta(\overline{Z})$, then we have the estimate
\begin{equation}
\label{EqLipschitzEstimateProduct}
\begin{split}
&\Big\| \big( (f,f')\cdot\widehat{\zeta} , f\big)^\sharp - \big( (h, h')\cdot\widehat{\overline{\zeta}}, h\big)^\sharp  \Big\|_{\calC^{\alpha+\gamma}} \\ 
&\lesssim C_M \left\{d_{{\sf C}^\beta_\alpha}\big((f,f'),(h,h')\big) + \big\|Z-\overline{Z}\big\|_{\calC^\alpha} + \big\|\zeta - \overline{\zeta}\big\|_{\calC^\gamma} + \big\|\zeta^{(2)} - \overline{\zeta}^{(2)}\big\|_{\calC^{\alpha+\gamma}}\right\},
\end{split}
\end{equation}
where $C_M$ is a positive constant with polynomial growth in 
{\small   $$
M := \max\Big\{\|Z\|_{\calC^\alpha}, \|\zeta\|_{\calC^\gamma}, \big\|\zeta^{(2)}\big\|_{\calC^{\alpha+\gamma}}, \|\overline{Z}\|_{\calC^\alpha}, \big\|\overline{\zeta}\big\|_{\calC\gamma}, \big\| \overline{\zeta}^{(2)}\big\|_{\calC^{\alpha+\gamma}}, \big\|(f,f')\big\|_{{\sf C}^\beta_\alpha}, \big\|(h,h')\big\|_{{\sf C}^\beta_\alpha} \Big\}.
$$   }
\end{itemize}   }
\end{theorem}

\medskip

Running backward the computations preceeding theorem \ref{thm:prod} one sees that $(f,f')\cdot\zeta$ coincides with $f\zeta$ if $f$ and $\zeta$ are both smooth -- in which case one can choose $f'=0$.

\medskip

\begin{Dem}
We examine the regularity of each terms of the defining identity \eqref{EqDefnProduct}. By Proposition \ref{prop:para}, we have $\Pi^{(b)}_f(\zeta)\in \calC^{\gamma}$ and $\Pi^{(b)}_\zeta(f) \in \calC^{\alpha+\gamma}$. Proposition \ref{prop:diag} yields that $ \Pi^{(b)} \big( (f,f')^\sharp,\zeta\big) \in \calC^{\alpha+\beta+\gamma}$. Applying Proposition \ref{prop:commutator} on the continuity properties of the commutator operator, with $\alpha+\gamma<0$, we see that $C(Z,f',\zeta)\in \calC^{\delta}$, with $\delta=\min(\alpha+\beta,1)+\gamma$. Since $\alpha+\gamma<\beta$, then $f'\zeta^{(2)} \in \calC^{\alpha+\gamma}$. Each term in formula \eqref{EqDefnProduct} is then an element of $\calC^\gamma$, and \eqref{eq:DD} holds since $\alpha<1$. The proof of the Lipschitz estimate \eqref{EqLipschitzEstimateProduct} is left to the reader. \end{Dem}

\medskip

Let insist here on the fact that $Z$ is not sufficient by itself to define a product operation, and that different choices of $\zeta^{(2)}$ provide different definitions of the product operation. In another direction, using the paralinearization formula, we are able to study the action of a nonlinearity on paracontrolled distributions, giving us the equivalent of the elementary fact that, in the above classical controlled setting for rough differential equations, the image by some map  $\sigma$ of a path $(x,x')$ controlled by some reference path $X$.

\medskip

\begin{theorem}  \label{thm:nolinear}    {\sf
Let $0<\beta<\alpha<1$, a reference function $Z\in\calC^\alpha$ and $(f,f')\in{\sf C}_\alpha^\beta(Z)$ be given, together witha function  $\textrm{F}\in C^4_b(\RR,\RR)$. Then $\big(\textrm{F}(f),\textrm{F}'(f)f' \big)$ belongs to ${\sf C}^\beta_\alpha(Z)$, and 
$$ 
\Big\| \big( \textrm{F}(f), \textrm{F}'(f)f' \big) \Big\|_{{\sf C}^\beta_\alpha} \lesssim \|\textrm{F}\|_{C^3_b} \big(1+\big\|(f,f')\big\|_{{\sf C}^\beta_\alpha}^2\big) \big(1+\|Z\|_{\calC^\alpha}^2\big).
$$
Moreover, this operation is locally Lipschitz in the sense that we have, with the same notations as in theorem \ref{thm:prod} and $(h,h')\in{\sf C}_\alpha^\beta(\overline{Z})$,
$$ 
d_{{\sf C}^\beta_\alpha} \Big( \big(\textrm{F}(f),\textrm{F}'(f)f'\big), \big(\textrm{F}(h),\textrm{F}'(h)h'\big) \Big) \lesssim C_M \|\textrm{F}\|_{C^4_b}^3 \left(d_{{\sf C}^\beta_\alpha} \big((f,g),(h,k)\big) + \big\|Z-\overline{Z}\big\|_{\calC^\alpha} \right),
$$
where $C_M$ is a constant with a polynomial growth in 
$$
M := \max\left\{\|Z\|_{\calC^\alpha}, \big\|\overline{Z}\big\|_{\calC^\alpha}, \big\|(f,f')\big\|_{{\sf C}^\beta_\alpha}, \big\|(h,h')\big\|_{{\sf C}^\beta_\alpha} \right\}.
$$   }
\end{theorem}

\medskip

\begin{Dem}
Given $f\in \calC^\alpha$, we have $\textrm{F}(f)\in \calC^\alpha$, since $F$ is Lipschitz. We know that $\textrm{F}'(f) f' \in \calC^\beta$, since $\textrm{F}'(f)\in \calC^\alpha$ and $f'\in \calC^\beta$. Using the notations of Theorem \ref{thm:paralinearization}, we have
$$ 
\textrm{F}(f) - \Pi^{(b)}_{\textrm{F}'(f) f'}(Z)  = \Pi^{(b)}_{\textrm{F}'(f)}(f) - \Pi^{(b)}_{\textrm{F}'(f) f'}(Z) + R_\textrm{F}(f)
$$
with a remainder $R_\textrm{F}(f)\in \calC^{2\alpha} \subset \calC^{\alpha+\beta}$. Since $f = \Pi^{(b)}_{g}(Z) + (f,f')^\sharp$, we have
$$  
\Pi^{(b)}_{\textrm{F}'(f)}(f) = \Pi^{(b)}_{\textrm{F}'(f)}\Pi^{(b)}_{f'}(Z) + \Pi^{(b)}_{\textrm{F}'(f)} \big((f,f')^\sharp\big)
$$
with $ \Pi^{(b)}_{\textrm{F}'(f)}\big((f,f')^\sharp\big) \in \calC^{\alpha+\beta}$, after Proposition \ref{prop:para}. So
$$ 
\textrm{F}(f) - \Pi^{(b)}_{\textrm{F}'(f) f'}(Z) \in \Pi^{(b)}_{\textrm{F}'(f)}\Pi^{(b)}_{f'}(Z) - \Pi^{(b)}_{\textrm{F}'(f) f'}(Z) + \calC^{\alpha+\beta}.
$$
Using Theorem \ref{thm:composition}, we deduce that 
$$ 
\textrm{F}(f) - \Pi^{(b)}_{\textrm{F}'(f) f'}(Z) \in \calC^{\alpha+\beta},
$$
which indeed shows that $\big(\textrm{F}(f),\textrm{F}'(f)f'\big) \in {\sf C}^\beta_\alpha(Z)$. We let the reader check the Lipschitz inequality for this operation.
\end{Dem}

\medskip

\subsection[\hspace{-1cm} {\sf Schauder estimates for paracontrolled distributions}]{Schauder estimates for paracontrolled distributions}
\label{SubsectionParacontrolledSchauder}

The above definition of a paracontrolled distribution is adapted to a time-independent setting. To deal with the time-dependent setting needed to handle the parabolic equations considered in practical examples, we use an adapted notion. Recall the definition of the space $\mcE^{(\alpha)}_{[0,T]}$ given in definition \ref{def:space}.

\medskip

\begin{definition*}    {\sf
Let $\alpha\in\R$ and $\beta>0$ be such that $\alpha+\beta\in(0,2)$; fix a reference distribution $Z\in \mcE^{(\alpha)}_{[0,T]}$, for some finite positive horizon $T$. A pair of distributions $(f,f')\in \mcE^{(\alpha)}_{[0,T]} \times \mcE^{(\beta)}_{[0,T]}$ is said to be \textbf{paracontrolled by $Z$} if 
$$
(f,f')^\sharp:=f-\Pi^{(b)}_{f'}(Z) \in C_T\calC^{\alpha+\beta} \cap C_T^{\frac{\beta}{2}} L^\infty.
$$ 
In such a case, we write $(f,f')\in {\sf C}^\beta_{\alpha,[0,T]}(Z)$ and define the norm
$$ 
\big\|(f,f')\big\|_{{\sf C}^\beta_{\alpha,[0,T]}} := \|f'\|_{\mcE_T^{(\beta)}} + \big\|(f,f')^\sharp\big\|_{C_T\calC^{\alpha+\beta}} + \big\|(f,f')^\sharp\big\|_{C_T^{\frac{\beta}{2}} L^\infty}.
$$
If $\overline{Z}\in \mcE_T^\alpha$ is another reference distribution, and $(h,h')\in {\sf C}^\beta_{\alpha,T}(\overline{Z})$, we set
$$ 
d_{{\sf C}^\beta_{\alpha,[0,T]}}\big((f,f'),(h,h')\big) := \|f'-h'\|_{\mcE_T^{(\beta)}} + \big\|(f,f')^\sharp - (h,h')^\sharp \big\|_{C_T\calC^{\alpha+\beta} \cap C_T^{\frac{\beta}{2}} L^\infty}.
$$   }
\end{definition*}

\medskip

Let us point out that the remainder $(f,f')^\sharp$ is not a priori in $\mcE_{[0,T]}^{(\alpha+\beta)} = C_T\calC^{\alpha+\beta} \cap C_T^{\frac{\alpha+\beta}{2}} L^\infty$, as the time-regularity of elements in the latter space is stronger than what is assumed in the above definition.

\medskip

\begin{theorem}  \label{thm:az}    {\sf
Let $\beta\in(0,1), \alpha\in(0,2-\beta)$, and a fixed positive finite time horizon $T$ be given. Given $\zeta\in C_T\calC^{\alpha-2}$, let $Z$ be the solution on $[0,T)$ of the equation
$$ 
\LL X := (\partial_t+L) Z = \zeta, \qquad Z_{\big| t=0}=0.
$$
Given $f'\in \mcE_T^\beta$ and $h\in C_T\calC^{\alpha+\beta-2}$, denote by $f$ the solution to the initial value problem
$$ 
\LL f = \Pi^{(b)}_{f'}(\zeta) + h, \qquad f_{\big| t=0}=f_0\in \calC^{\alpha+\beta}.
$$
Then $(f,f') \in {\sf C}^\beta_{\alpha,[0,T]}(Z)$ and 
$$ 
\big\|(f,f')\big\|_{{{\sf C}^\beta_{\alpha,[0,T]}}}\lesssim \|f_0\|_{\calC^{\alpha+\beta}} + (1+T)\left\{ \|f'\|_{\mcE^\beta_T} \big(1+\|\zeta\|_{C_T\calC^{\alpha-2}} \big) + \|h\|_{C_T\calC^{\alpha+\beta-2}} \right\}.
$$
Moreover, the map which associates $(f,f')$ to $(\zeta,Z,f',h,f_0)$ is locally Lipschitz.   }
\end{theorem}

\medskip

\begin{Dem} 
Since $\zeta\in C_T\calC^{\alpha-2}$ the Schauder estimates, Corollary \ref{cor:schauder1}, yield that $Z\in \mcE^{(\alpha)}_{[0,T]}$ and $f\in\mcE_T^{(\alpha)}$; so we are left with checking that
\begin{equation} 
\label{eq:tc}
(f,f')^\sharp:= f-\Pi^{(b)}_{f'}(Z) \in C_T\calC^{\alpha+\beta} \cap C_T^{\frac{\beta}{2}} L^\infty. 
\end{equation}
Let us derive an equation for this quantity
\begin{align*} 
\LL (f,f')^\sharp & = \LL f - \LL \Pi^{(b)}_{f'}(Z) = \Pi^{(b)}_{f'}(\zeta) - \LL \Pi^{(b)}_{\zeta}(Z) + h   \\
 & = \Pi^{(b)}_{f'}(\LL Z) - \LL \Pi^{(b)}_{f'}(Z) + h \\
 & = \Big[\LL, \Pi^{(b)}_{f'}\Big] (Z) + h.
\end{align*}
We have $h\in C_T\calC^{\alpha+\beta-2}$, and we have seen that $\Pi^{(b)}_{f'}(\zeta)\in C_T\calC^{\alpha-2}$ and $\Pi^{(b)}_{f'}(Z) \in C_T\calC^{\alpha}$, so that $\LL \Pi^{(b)}_{f'}(Z) \in C_T\calC^{\alpha-2}$. By studying the difference (which consist to commute the paraproduct $\Pi_{f'}^{(b)}$ with $\LL$) with introducing an intermediate time-space paraproduct\footnote{We refer the reader to a more recent work by the authors \cite{BBF16} where such time-space paraproducts are more crucial to deal with higher dimensional space. More details are also provided there. Indeed, they consist of tensorial product between the semigroup and some standard $1$-dimensional approximations of the unity to localize the time variable.}, such as done in \cite[Lemma 5.1]{GIP} -- whose proof can easily be extended to our setting, we obtain that
$$ 
\mcR\Big[\LL, \Pi^{(b)}_{f'}\Big] (Z) \in C_T\calC^{\alpha+\beta} \cap C_T^{\frac{\beta}{2}} L^\infty,
$$
where we recall that $\mcR$ is the resolution operator of heat equation. We invite the reader to check the Lipschitz inequality for this operation, in terms of $\zeta,Z,f',h$ and $f_0$.
\end{Dem}

\bigskip

The fact that $\mcR\Big[\LL, \Pi^{(b)}_{f'}\Big] (Z)$ takes values in $C_T\calC^{\alpha+\beta} \cap C_T^{\frac{\beta}{2}} L^\infty$ is the reason why we define the space of paracontrolled distributions as above rather than with the more natural norm $\mcE_{[0,T]}^{\alpha+\beta}$ to measure the size of the remainder.

\ssk 

With this result in hands, we now have all the theoretical apparatus needed to study some examples of singular parabolic PDEs. We have chosen to illustrate our machinery on what may be one of the simplest examples of such an equation, the generalized parabolic Anderson equation, (gPAM), that was already handled in the 2-dimensional torus both by Hairer in \cite{H} using his theory of regularity structures, and by Gubinelli, Imkeller and Perkowski in \cite{GIP}, using their Fourier-based paracontrolled approach. This choice is motivated by the fact that only one (probabilistic) renormalization is needed to implement the paracontrolled machinery, while further renormalizations are needed in the stochastic quantization or KPZ equations. So the reader can see in the next section the machinery at work without being overwhelmed by side probabilistic matters. It makes sense now to make the following definition, in the present setting, where $\alpha$ stands for a real number in $(0,1)$.

\medskip

\begin{definition*}
An \textbf{enhanced, or rough, distribution} $\widehat{\zeta}$ is a pair $\big(\zeta,\zeta^{(2)}\big)$, with $\zeta\in\calC^{\alpha-2}$ and $\zeta^{(2)}\in C_T\calC^{2\alpha-2}$. 
\end{definition*}

\bigskip

\section[\hspace{0.6cm} {\sf The (generalized) parabolic Anderson Model in dimension 2}]{The (generalized) parabolic Anderson Model in dimension 2} 
\label{sec:pam}

This section is devoted to the study in our abstract setting of the (generalized) parabolic Anderson Model, in dimension $2$. The analytical/geometrical \textbf{setting} is described  in Section \ref{sec:preli}. The space $(M,d,\mu)$ is a space of homogeneous type, equipped with a semigroup $\big(e^{-tL}\big)_{t>0}$ satisfying the regularity assumptions \eqref{UE} and \eqref{Lipschitz}. Let us insist here on the fact that even in this modest setting, the above semigroup approach offers some results that seem to be beyond the present scope of the theory of regularity structures, in so far as we are for instance  allowed to work in various underlying spaces and even in the Euclidean space with operators $L$ of the form $\textrm{div}\big(A\nabla\big)$, with $A$ H\"older continuous -- see example 2 in section \ref{SubsectionHeatSemiGroup}. The first two subsections are dedicated to proving some local and global in time well-posedness results, for the deterministic (gPAM) and (PAM) equations respectively. To turn that machinery into an efficient tool for investigating stochastic PDEs in which the singular term involves a Gaussian noise, we need to lift this noise into an enhanced/rough distribution; this step requires a probabilistic limit procedure generically called a \textit{renormalization step}. It is performed in section \ref{SubsectionRenormalization}, in the geometric framework of a potentially unbounded manifold, when working with a weighted noise.

\bigskip

\subsection[\hspace{-1cm} {\sf Local well-posedness result for generalized PAM}]{Local well-posedness result for generalized PAM}
\label{SubsectionPAM}

A big enough parameter $b$ is fixed for all the previous results to hold.

\medskip

\begin{theorem} \label{thm:gpam}     {\sf
Let $\alpha\in\big(\frac{2}{3},1\big)$ be given, and $\alpha'<\alpha$ be close enough to $\alpha$ to have $2\alpha+\alpha'-2>0$; let also a finite positive time horizon $T$ be given. Fix an initial data $u_0\in \calC^{2\alpha}$, and a nonlinearity $\textrm{F}\in C^3_b(\RR,\RR)$. Given $\zeta \in \calC^{\alpha-2}$, set $Z := \mcR(\zeta)$, and assume there is a sequence of smooth functions $\zeta_n$ converging to $\zeta$ in $\calC^{\alpha-2}$, and a sequence of smooths functions $Z_n$ converging to $Z$ in $\mcE^{(\alpha)}_{[0,T]}$, such that $\Pi^{(b)}(Z_n,\zeta_n)$ converges in $C_T\calC^{2\alpha-2}$ to some limit element $\zeta^{(2)}$. Write $\widehat{\zeta}=\big(\zeta,\zeta^{(2)}\big)$ for the associated enhanced distribution. Then the generalized parabolic Anderson model equation (gPAM)  
$$ 
\partial_t u + Lu = \big(\textrm{F}(u),\textrm{F}'(u)u'\big)\cdot\widehat{\zeta},\qquad u(0) = u_0 
$$
has a unique paracontrolled solution $\big(u,u'\big) \in {\sf C}^{\alpha'}_{\alpha,[0,T]}(Z)$, with $u'=F(u)$, provided $T$ is small enough.   }
\end{theorem}

\medskip

Recall that the very notion of product in the right hand side of the (gPAM)equation depends on $\widehat{\zeta}$. Since, we have established in the previous sections the main analytic estimates of paracontrolled calculus, we can prove this theorem \ref{thm:gpam} following the pattern of proof devised by Gubinelli, Imkeller and Perkowski in their seminal work \cite{GIP}, such as extended here to our more abstract setting. 

\medskip 

\begin{Dem}
Given a singular spatial perturbation $\zeta \in \calC^{\alpha-2}$, we know, by Proposition \ref{prop:schauder1}, that $Z := \mcR(\zeta)\in \mcE^{(\alpha)}_{[0,T]}$. Let $\frak{A}_{[0,T]}(K)$ be the set of controlled distributions $(u,u')\in {\sf C}^{\alpha'}_{\alpha,[0,T]}(Z)$ such that 
$$
\big\|(u,u')\big\|_{{\sf C}^{\alpha'}_{\alpha,[0,T]}} \leq K,\quad  u'(0) = \textrm{F}(u_0) \ \textrm{and} \ (u,u')^\sharp(0)=u_0.
$$
Denoting by $v$ the solution of the equation 
$$ 
\partial_t v + Lv = \big(\textrm{F}(u),\textrm{F}'(u)u'\big) \cdot \widehat\zeta, \quad v(0)=u_0,
$$
we define on $\frak{A}_{[0,T]}(K)$ a map $\Gamma$ setting
$$
\Phi\big((u,u')\big) = \big(v,\textrm{F}(u)\big);
$$
it also takes its values in ${\sf C}_{\alpha,[0,T]}^{\alpha'}(Z)$, by Theorems \ref{thm:prod}, \ref{thm:nolinear} and \ref{thm:az}, and has a ${\sf C}_{\alpha,[0,T]}^{\alpha'}(Z)$-norm bounded above by 
{\small   \begin{equation*}
\begin{split} 
&\big\|u_0\big\|_{\calC^{2\alpha}} + T^{\frac{\alpha-\alpha'}{2}}\left( \big\|\textrm{F}(u)\big\|_{\mcE^{(\alpha)}_{[0,T]}} \big(1+\|\zeta\|_{\calC^{\alpha-2}}\big) + \Big\|\big(\textrm{F}(u) , g-\textrm{F}'(u) u'\big)\cdot \zeta - \Pi^{(b)}_{\textrm{F}(u)} (\zeta)\Big\|_{C_T\calC^{2\alpha-2}} \right)   \\
&\leq \big\|u_0\big\|_{\calC^{2\alpha}} \\
& + T^{\frac{\alpha-\alpha'}{2}} \left(\|\textrm{F}\|_{C^3_b}\|u\|_{\mcE^{(\alpha)}_{[0,T]}}\big(1+\|\zeta \|_{\calC^{\alpha-2}}\big) + \|\textrm{F}\|_{C^3_b} \big(1+\big\|(u,u')\big\|_{{\sf C}^{\alpha'}_{\alpha,[0,T]}}^2\big)\big(1+\|Z\|_{C_T\calC^\alpha}^2\big)(\star)\right), 
\end{split} 
\end{equation*}   } 
with 
$$
(\star) := \|\zeta\|_{\calC^{\alpha-2}}+\|\textrm{F}\|_{C^3_b}\|u\|_{C_T\calC^\alpha}\|\zeta\|_{\calC^{\alpha-2}} + \big\|\zeta^{(2)}\big\|_{C_T\calC^{2\alpha-2}}.
$$
from the estimates provided in the above-mentioned theorems, together with Remark \ref{rem:T} describing how we can obtain the extra factor in terms of $T^{\frac{\alpha-\alpha'}{2}}$. The map $\Phi$ is then a contraction from $\frak{A}_{[0,T]}(K)$ to itself if $K$ is chosen large enough and $T$ small enough. Moreover, the Lipchitz estimates in Theorems \ref{thm:prod} and \ref{thm:nolinear} give us the estimate
\begin{align*}
d_{{\sf C}^{\alpha'}_{\alpha,[0,T]}}\Big(\Phi(u,u'),\Phi(v,v')\Big) \lesssim T^{\frac{\alpha-\alpha'}{2}}\left(\|u-v\|_{\mcE^{(\alpha')}_T} + d _{{\sf C}^{\alpha'}_{\alpha,[0,T]}}\big((u,u'),(v,v')\big)\right) 
\end{align*}
for all $(u,u'), (v,v')$ in $\frak{A}_{[0,T]}(K)$, with implicit constants depending only on $K$, $Z$ and $\zeta$. So $\Phi$ happens to be a contraction of $\frak{A}_{[0,T]}(K)$, provided the time horizon $T$ is small enough, from which it follows that $\Phi$ has a unique fixed point in $\frak{A}_{[0,T]}(K)$. IFinally, it is easy to proceed as in Gubinelli, Imkeller, Perkowski's work \cite{GIP} and check that a solution of (gPAM) has to be in $\frak{A}_{[0,T]}(K)$, at least for a small enough $T$.
\end{Dem} 
 
\bigskip

\subsection[\hspace{-1cm} {\sf Global well-posedness result for linear PAM}]{Global well-posedness result for linear PAM}
\label{SubsectionGlobalPAM}

We focus in this subsection on the linear (PAM) equation and prove a global in time well-posedness result in that setting. With that aim in mind, we define a weighted (in time) version of the previous functional spaces of paracontrolled distributions.

\medskip

\begin{definition*}    {\sf 
Given $\lambda\geq 1$ and $\alpha\in(-2,2)$, set 
$$ 
\mathcal{W}_\lambda \calC^\alpha := \left\{ f\in L^\infty_{\textrm{loc}}(\SSS),\  \|f\|_{\mathcal{W}_\lambda \calC^\alpha}:=\sup_{t\geq 0}\ e^{-\lambda t}\,\big\|f(t)\big\|_{\calC^\alpha}<\infty\right\}
$$
and  
$$ 
\mathcal{W}_\lambda^{\frac{\alpha}{2}} L^\infty:= \left\{f\in L^\infty_{\textrm{loc}}(\SSS),\ \|f\|_{\mathcal{W}_\lambda^{\frac{\alpha}{2}} L^\infty} :=\sup_{\genfrac{}{}{0pt}{}{s\neq t}{0\leq s,t \leq 1}} \ e^{-\lambda s}\,\frac{\big\|f(t)-f(s)\big\|_\infty}{|t-s|^{\frac{\alpha}{2}}} <\infty \right\}.
$$
We then define the space
$$
\mathcal{W}^{(\alpha)}_\lambda := \mathcal{W}_\lambda\calC^\alpha \cap \mathcal{W}_\lambda^{\frac{\alpha}{2}} L^\infty,
$$ 
and, given a reference distribution $Z\in\mathcal{W}^{(\alpha)}_\lambda$, we define accordingly the space ${\sf C}^{\beta}_{\alpha\,;\lambda}(Z)$ of pairs of distributions $(f,f')$ in $\mathcal{W}^{(\alpha)}_\lambda\times\mathcal{W}^{(\beta)}_\lambda$ such that 
$$
(f,f')^\sharp := f - \Pi^{(b)}_{f'}(Z) \in C_T \big(\mathcal{W}_\lambda\calC^{\alpha+\beta}\big) \cap C^{\frac{\beta}{2}}L^\infty.
$$   
} \end{definition*}

\medskip

Following the reasoning of Theorem \ref{thm:gpam}, we prove in this section a global in time well-posedness result. One of the main ingredients used in the proof of Theorem \ref{thm:gpam} was the Schauder estimates, through Proposition \ref{prop:schauder1} or Corollary \ref{cor:schauder1}. We now give an extension of these estimates to the setting provided by the above exponentially weighted spaces.

\medskip

\begin{proposition}  \label{prop:schauder2}     {\sf
Consider an integer $b\geq 0$, $\beta\in(-2,0)$ and $\lambda\geq 1$. For every $\epsilon\in(0,1)$ and $v\in \mathcal{W}_\lambda \calC^\beta$ the function $\mcR(v)_t := \int_0^t P^{(b)}_{t-s} v(s)\, ds$ belongs to $\mathcal{W}_\lambda \calC^{\beta+2-2\epsilon}$ and satisfies the $\lambda$-uniform bounds
$$ 
\big\|\mcR(v)\big\|_{\mathcal{W}_\lambda \calC^{\beta+2-2\epsilon}} \lesssim \lambda^{-\epsilon} \, \|v\|_{\mathcal{W}_\lambda \calC^\beta}.
$$
and
$$ 
\big\|\mcR(v)\big\|_{\mathcal{W}_\lambda^{\frac{\beta+2-2\epsilon}{2}}L^\infty} \lesssim \lambda^{-\epsilon} \, \|v\|_{\mathcal{W}_\lambda \calC^\beta}.
$$
Consequently,
$$ 
\big\|\mcR(v)\big\|_{\mathcal{W}^{\frac{\beta+2-2\epsilon}{2}}_\lambda} \lesssim \lambda^{-\epsilon} \, \|v\|_{\mathcal{W}_\lambda \calC^\beta}.
$$   }
\end{proposition}

\medskip

\begin{Dem} 
We adapt the proof of Proposition \ref{prop:schauder1} and add an extra new argument to consider the exponential weight in time. So consider another integer $a \geq \frac{|\beta|}{2} + 1$, and a parameter $\tau\in(0,1]$. Then
 $$  
 Q_\tau^{(a)}\big(\mcR(v)_t\big) = \int_0^t Q_\tau^{(a)} P^{(b)}_{t-s} v(s) \, ds.
 $$
Hence,
\begin{align*}
\Big\|Q_\tau^{(a)} P^{(b)}_{t-s} v(s)\Big\|_{\infty} & \lesssim \left(\frac{\tau}{\tau + t-s} \right)^a \Big\|Q_{\tau+t-s}^{(a)} v(s)\Big\|_{\infty} \\
 & \lesssim \left(\frac{\tau}{\tau + t-s}\right)^a (\tau + t-s)^{\frac{\beta}{2}} \big\|v(s)\big\|_{\calC^\beta}.
\end{align*}
So by integrating, it comes
\begin{align*}
 e^{-\lambda t} \Big\| Q_\tau^{(a)}\big(\mcR(v)_t\big) \Big\|_\infty & \lesssim \left(\int_0^t e^{-\lambda(t-s)} \left(\frac{\tau}{\tau+t-s}\right)^a (\tau + t-s)^{\frac{\beta}{2}} \, ds\right) \sup_{s\in[0,t]} e^{-\lambda s} \big\|v(s)\big\|_{\calC^\beta} \\
  & \lesssim \left(\int_0^t e^{-\lambda(t-s)} \left(\frac{\tau}{\tau+t-s}\right)^a (\tau + t-s)^{\frac{\beta}{2}} \, ds\right) \|v\|_{\mathcal{W}_\lambda \calC^\beta}.
\end{align*}
Let us just consider the integral term, temporarily denoted by $I$. If $t\leq \tau$, then 
\begin{align*} 
I & \leq \tau^{\frac{\beta}{2}} \int_0^t e^{-\lambda(t-s)} \, ds \lesssim \tau^{\frac{\beta}{2}} \frac{1-e^{-\lambda t}}{\lambda} \lesssim \tau^{\frac{\beta}{2}} (\lambda t)^{(1-\epsilon)} \lambda^{-1} \\
  & \lesssim \tau^{\frac{\beta}{2}+1-\epsilon} \lambda^{-\epsilon}.
\end{align*}
If $t\geq \tau$ then
\begin{align*} 
I & \leq \tau^{\frac{\beta}{2}} \int_{t-\tau}^t e^{-\lambda(t-s)} \, ds + \int_0^{t-\tau} e^{-\lambda(t-s)} \left(\frac{\tau}{t-s}\right)^a (t-s)^{\frac{\beta}{2}} \, ds \\
 & \lesssim \tau^{\frac{\beta}{2}} \frac{1-e^{-\lambda \tau}}{\lambda} +  \tau^{\frac{\beta}{2}+1-\epsilon} \int_0^{t-\tau} e^{-\lambda(t-s)} (t-s)^{\epsilon-1}  \, ds \\
 & \lesssim \tau^{\frac{\beta}{2}} \frac{1-e^{-\lambda \tau}}{\lambda} +  \tau^{\frac{\beta}{2}+1-\epsilon} \lambda^{-\epsilon} \left(\int_0^{\infty} e^{-x} x^{\epsilon-1}  \, dx\right) \\
 & \lesssim \tau^{\frac{\beta}{2}+1-\epsilon} \lambda^{-\epsilon} . 
\end{align*}
So in both situations, we deduce that uniformly in $\lambda\geq 1$ and $t>0$, it comes
$$ 
e^{-\lambda t} \Big\| Q_\tau^{(a)}\big(\mcR(v)_t\big) \Big\|_\infty \lesssim \tau^{\frac{\beta}{2}+1-\epsilon} \lambda^{-\epsilon} \|v\|_{\mathcal{W}_\lambda \calC^\beta} 
$$
and similarly
$$ 
e^{-\lambda t} \Big\| P_1\big(\mcR(v)_t\big) \Big\|_\infty \lesssim \lambda^{-\epsilon} \|v\|_{\mathcal{W}_\lambda \calC^\beta}.
$$
Consequently, we deduce that for every $t\geq 0$
$$ 
e^{-\lambda t} \big\|\mcR(v)_t\big\|_{\calC^{\beta+2-2\epsilon}} \lesssim \lambda^{-\epsilon} \|v\|_{\mathcal{W}_\lambda \calC^\beta},
$$
which yields
$$ 
\big\|\mcR(v)\big\|_{\mathcal{W}_\lambda \calC^{\beta+2-2\epsilon}} \lesssim \lambda^{-\epsilon} \|v\|_{\mathcal{W}_\lambda \calC^\beta}.
$$
For the second statement, for $s<t$ we have
\begin{align*}
\mcR(v)_t - \mcR(v)_s   &=  \big(P_{t-s}^{(b)}-\textrm{Id}\big)\big(\mcR(v)_s\big) + \int_s^t P_{t-r}^{(b)}\big(v(r)\big) \, dr \\ 
      									&=  \int_0^{t-s} Q^{(b)}_{r} \mcR(v)_s \, \frac{dr}{r} + \int_s^t P_{t-r}^{(b)}\big(v(r)\big) \, dr.
 \end{align*}
So
\begin{align*} 
e^{-\lambda t} \left\| \int_0^{t-s} Q^{(b)}_{r} \mcR(v)_s \, \frac{dr}{r} \right\|_{\infty} & \lesssim \left(\int_0^{t-s} r^{\frac{\beta}{2}+1} \, \frac{dr}{r} \right) e^{-\lambda t} \big\|\mcR(v)_s\big\|_{\calC^{\beta+2}} \\ 
 & \lesssim (t-s)^{\frac{\beta}{2}+1} e^{-\lambda (t-s)} e^{-\lambda (t-s)} \big\|\mcR(v)\big\|_{\mathcal{W}_\lambda \calC^{\beta+2}} \\
  & \lesssim (t-s)^{\frac{\beta}{2}+1} \big(\lambda(t-s)\big)^{-\epsilon}  \big\|\mcR(v)\big\|_{\mathcal{W}_\lambda \calC^{\beta+2}} \\
  & \lesssim (t-s)^{\frac{\beta}{2}+1-\epsilon} \lambda^{-\epsilon}  \big\|\mcR(v)\big\|_{\mathcal{W}_\lambda \calC^{\beta+2}} 
\end{align*}
and also (since $\beta<0$)
\begin{align*}
e^{-\lambda t} \left\| \int_s^t P_{t-r}^{(b)}\big(v(r)\big) \, dr \right\|_{L^\infty} & \lesssim \int_s^t e^{-\lambda t} \left(\int_{t-r}^1  \Big\|Q_{\tau}^{(b)} v(r)\Big\|_\infty \, \frac{d\tau}{\tau} + \Big\|P_{1}^{(b)}\big(v(r)\big)\Big\|_{\infty}\right) dr \\
& \lesssim \int_s^t  e^{-\lambda t} \left(\big\|v(r)\big\|_{\calC^\beta} \int_{t-r}^1 \tau^{\frac{\beta}{2}} \, \frac{d\tau}{\tau} + \Big\|P_{1}^{(b)}\big(v(r)\big)\Big\|_{L^\infty}\right) dr \\
& \lesssim \|v\|_{\mathcal{W}_\lambda \calC^\beta} \int_s^t e^{-\lambda (t-r)} \left( \int_{t-r}^1 \tau^{\frac{\beta}{2}} \, \frac{d\tau}{\tau} + 1\right) dr \\
& \lesssim \|v\|_{\mathcal{W}_\lambda \calC^\beta} \int_s^t e^{-\lambda (t-r)} \left((t-r)^{\frac{\beta}{2}} +1 \right) dr \\
& \lesssim \lambda^{-\epsilon} \|v\|_{\mathcal{W}_\lambda \calC^\beta} (t-s)^{\frac{\beta}{2}+1-\epsilon},
\end{align*}
where we used the fact that $\frac{\beta}{2}+1\in (0,1)$. So we conclude that
$$ 
\big\|\mcR(v)\big\|_{\mathcal{W}^{\frac{\beta}{2}+1-\epsilon}_\lambda L^\infty} \lesssim \lambda^{-\epsilon} \|v\|_{\mathcal{W}_\lambda \calC^\beta}.
$$
\end{Dem}

\medskip

\begin{theorem} \label{thm:gpam-2}   {\sf
Let $\alpha\in\big(\frac{2}{3},1\big)$ be given, and $\alpha'<\alpha$ be close enough to $\alpha$ to have $2\alpha+\alpha'-2>0$; fix also an initial data $u_0\in \calC^{2\alpha}$, together with some exponent $\lambda\geq 1$. Given $\zeta \in \calC^{\alpha-2}$, set $Z := \mcR(\zeta)$, and assume there is a sequence of smooth functions $\zeta_n$ converging to $\zeta$ in $\calC^{\alpha-2}$, and a sequence of smooths functions $Z_n$ converging to $Z$ in $\mathcal{W}^{(\alpha)}_\lambda$, such that $\Pi^{(b)}(Z_n,\zeta_n)$ converges in $\mathcal{W}_\lambda\calC^{2\alpha-2}$ to some limit element $\zeta^{(2)}$. Write $\widehat{\zeta} = \big(\zeta,\zeta^{(2)}\big)$ for the associated enhanced distribution. Then, for a large enough choice of parameter $\lambda$, the linear (PAM) equation 
$$ 
\partial_t u + Lu = (u,u')\cdot\widehat{\zeta} ,\qquad u(0)=u_0,
$$
has a unique global in time solution $\big(u,u'\big) \in {\sf C}^{\alpha'}_{\alpha\,;\lambda}(Z)$, with $u'=u$.   }
\end{theorem}

\medskip

\begin{Dem}
Since $\zeta  \in \calC^{\alpha-2}$, Proposition \ref{prop:schauder2} implies that $Z = \mcR(\zeta )\in \mathcal{W}^{(\alpha)}_\lambda$. Define  $\frak{A}_\lambda(K)$ as the set of distributions $(u,u')\in {\sf C}^{\alpha'}_{\alpha\,;\lambda}(Z)$ with 
$$
\big\|(u,u')\big\|_{{\sf C}^{\alpha'}_{\alpha;\lambda}} \leq K, \qquad u'(0)=u_0, \ (u,u')^\sharp(0)=u_0,
$$
and define $\Phi_\lambda\big((u,u')\big) := (v,u)$, for $(u,u')\in\frak{A}_\lambda(K)$, as the solution of $v$ of the equation
$$ 
\partial_t v +Lv = (u,u') \cdot \widehat\zeta, \qquad v(0)=u_0.
$$
Then by the same considerations such as those detailled in the proof of Theorem \ref{thm:gpam}, with some minor simplifications getting around the non-necessary paralinearization step, and using Proposition \ref{prop:schauder2} instead of Proposition \ref{prop:schauder1}, we get the following result. The pair $\Phi_\lambda(u,u')$ belongs to ${\sf C}^{\alpha'}_{\alpha\,;\lambda}(Z)$ and satisfies 
\begin{align*}
\big\|\Phi_\lambda(u,u')\big\|_{{\sf C}^{\alpha'}_{\alpha\,;\lambda}} \lesssim \|u_0\|_{\calC^{2\alpha}} + \lambda^{-(\alpha-\alpha')}\left(\|u\|_{\mathcal{W}^{(\alpha)}_\lambda}\big(1+\|\zeta\|_{\calC^{\alpha-2}}\big) + \|(u,u')\|_{{\sf C}^{\alpha'}_{\alpha;\lambda}}(\star)\right), 
\end{align*}
with
$$
(\star) := \left(\|\zeta\|_{\calC^{\alpha-2}}+\|Z\|_{\mathcal{E}_\lambda\calC^\alpha}\|\zeta \|_{\calC^{\alpha-2}} + \big\|\zeta^{(2)}\big\|_{C_T\calC^{2\alpha-2}}\right),
$$
uniformly in $\lambda\geq 1$. A large enough choice of constants $K,\lambda$ ensures as a consequence that $\Phi_\lambda$ sends $\frak{A}_\lambda(K)$ into itself. Moreover, we also have for $(u,u')$ and $(v,v')$ in $\frak{A}_\lambda(K)$ the Lipschitz estimate
\begin{align*}
d_{{\sf C}^{\alpha'}_{\alpha;\lambda}}\Big(\Phi_\lambda(u,u'),\Phi_\lambda(v,v')\Big) \lesssim \lambda^{-(\alpha-\alpha')}\left(\|u-v\|_{\mathcal{W}^{(\alpha')}_\lambda} + d _{{\sf C}^{\alpha'}_{\alpha;\lambda}}\big((u,u'),(v,v')\big)\right) 
\end{align*}
with implicit constants depending only on $K$, $Z$ and $\zeta$. So the result follows from the fact that $\Phi_\lambda$ happens to be a contraction of $\frak{A}_\lambda(K)$ if $\lambda$ is chosen large enough.  
\end{Dem}

\bigskip

\subsection[\hspace{-1cm} {\sf Renormalization for a weighted noise}]{Renormalization for a weighted noise}
\label{SubsectionRenormalization}

We cannot expect to work in the Besov spaces used above when working  in unbounded ambiant spaces and with a spatial white noise; so weights need to be introduced, with a choice to be made. We can either put the weight in the H\"older spaces and still consider a uniform white noise, or we can put the weight on the noise and consider a weighted noise with values in unweighted H\"older spaces. The first approach has been recently implemented by Hairer and Labb\'e in a forthcoming work on the linear (PAM) equation in ${\mathbb R}^3$; see \cite{HL}. We choose to work with the second option here, partly motivated by exploring this unexplored question, partly because it seems to us that spatial white noise in an unbounded space has more something of a mathematical abstraction than of a model for real-life phenomena. Refer to our follow up work \cite{BBF16} for the use of weighted H\"older spaces in a paracontrolled setting.

\medskip

\begin{definition*}  {\sf
Let $\omega$ be an $L^2(\mu)$ weight on $M$; the {\bf noise with weight $\omega$} is the centered Gaussian process $\xi$ indexed by $L^2(\omega\mu)$, such that for every continuous function $f\in L^2(\omega\mu)$ we have
\begin{equation} 
\E\big[\xi(f)^2\big] = \int f^2(x)\,\omega(x)\mu(dx). 
\label{eq:wn} \end{equation}   }
\end{definition*}

\medskip

Let us define the following notation. For $t>0$, we denote by $\G_t$ the Gaussian kernel
$$ 
\G_t(x,y):=\frac{1}{V(x,\sqrt{t})} e^{-c\,\frac{d(x,y)^2}{t}};
$$
it also depends on the positive constant $c$, although we do not mention it in the notation for convenience. Since we shall need to commute in some sense the Gaussian kernels with the weight $\omega$, it seems natural to make the following assumption. We \textbf{assume} the existence of some implicit constants such that we have, for every $t\in (0,1]$ and every $x,y\in M$, 
\begin{equation} 
\omega(x)\,\G_{t}(x,y)\lesssim \omega(y)\,\G_{t}(y,x),
\label{ass:weight} 
\end{equation}
where we allow the implicit constants (in the exponential in $\G_t$) to be different from the left to the right hand side, but but that they are independent of $t,x,y$. 

\medskip

Recall the definition of \textbf{Ahlfors regularity} of a measure $\mu$ on a metric space $(M,d)$, given in section \ref{SubsectionHolder} before the Sobolev embedding theorem \ref{LemBesovEmbedding}, and quantified in equation \eqref{EqDoublingProperty}. In that setting, it is relatively elementary to use the latter and prove by classical means that a weighted noise, as defined above, has a realization that takes almost-surely its values in some H\"older space.   

\medskip

\begin{proposition} \label{prop:whitenoise} {\sf 
Assume that $(M,d,\mu)$ is Ahlfors regular and let $\xi$ be a noise on $M$, with weight $\omega \in L^1 \cap L^\infty$ satisfying the assumption \eqref{ass:weight}. Then, for every $\sigma < -\frac{\nu}{2}$, there exists a version of $\xi$, still denoted by the same symbol, which takes almost surely its values in $\calC^\sigma$.
}
\end{proposition}

\medskip

\begin{Dem}
It suffices from general principles and lemma \ref{LemBesovEmbedding} to check that the two expectations
$$
\EE\left(\int_{\RR^d}\big|e^{-L}f\big|^p(x)\mu(dx)\right)
$$
and 
$$
(\star) := \EE\left(\int_0^1 t^{-p\frac{\sigma}{2}} \Big\| Q^{(a)}_tf\Big\|_{p}^p\,\frac{dt}{t}\right)
$$
are finite for every $p>2$. We show how to deal with the second expectation, the first one being easier to treat with similar arguments. Starting from the fact that $\big(Q^{(a)}_tf\big)(x)$ is, for every $x\in M$, a Gaussian random variable with covariance the $L^2(\omega\mu)$-norm of $K_{Q^{(a)}_t}(x,\cdot)$, the equivalence of Gaussian moments with \eqref{eq:wn} give the upper bound 
\begin{equation*}
\begin{split}
\EE\Big[\Big| Q^{(a)}_tf \Big|^p (x) \Big] &\lesssim \EE\left[ \Big| Q^{(a)}_tf\Big|^2(x) \right]^\frac{p}{2} \\
&\lesssim \left|\int K_{Q^{(a)}_t}(x,z)^2\omega(z)\mu(dz)\right|^\frac{p}{2}.
\end{split}
\end{equation*}
Using the Gaussian bounds for the kernel of $Q^{(a)}_t$ with property \eqref{ass:weight} and Ahlfors regularity, this implies for $t\in(0,1)$
\begin{equation*}
\begin{split}
\EE\Big[\Big| Q^{(a)}_tf \Big|^p (x) \Big] &\lesssim \omega(x)^{\frac{p}{2}}\, t^{-\nu\frac{p}{4}}.
\end{split}
\end{equation*}
Hence, it follows that
\begin{equation*}
(\star) \lesssim \int_0^1\int_M \omega(x)^{\frac{p}{2}} \,t^{-p\frac{\nu}{4}} \,t^{-p\frac{\sigma}{2}} \,\mu(dx)\frac{dt}{t} \lesssim \|\omega\|_{\frac{p}{2}}
\end{equation*}
if $\sigma < -\frac{\nu}{2}$; the conclusion follows since $\omega\in L^1 \cap L^\infty \subset L^{\frac{p}{2}}$. 
\end{Dem}
 
\medskip

Let $\xi$ be a weighted noise, with weight $\omega$, and define for every $s>0$, a function $g_s:M \rightarrow \R$, by the formula
$$ 
g_s(x):=\E\Big[ \Pi \big(e^{-sL}\xi,\xi \big) (x)\Big];
$$
so that we formally have
$$ 
\int_0^\infty g_s(x) \,ds = \E\Big[ \Pi \big(L^{-1}\xi,\xi \big) (x)\Big].
$$
An explicit computation can be used in the case of the torus and the white noise to show that this integral diverges; see \cite{GIP}. A similar computation can be done in our setting with the help of a highly non-trivial estimate on the kernel of the operators $Q^{(1)}_s$, showing that the above integral also diverges at almost all points $x$ of $M$. These facts justifies that we consider the modified integral \eqref{EqRenormalizedIntegral} below. Even though we shall only use here theorem \ref{thm:Renormalization} in a $2$-dimensional setting, we prove it in the optimal range of homogeneous dimensions $d\in [2,4)$, for use in forthcoming works. Denote by $\Xi$ the function $\mcR(\xi)$ solution to the linear equation $\big(\partial_t+\Delta\big) \Xi = \xi$, with null initial condition.

\medskip

\begin{theorem}[Renormalization] \label{thm:Renormalization}  {\sf 
Assume that $(M,d,\mu)$ is locally Ahlfors regular, with homogeneous dimension $d=\nu\in[2,4)$. Consider $\xi$ a weighted noise with weight $\omega \in L^1\cap L^\infty$ satisfying assumptions \eqref{ass:weight}.
Set 
\begin{equation} 
\label{EqRenormalizedIntegral}
\big(\Xi \diamondsuit \xi\big)(t) := \int_0^t \Big\{\Pi^{(b)}\big(e^{-sL}\xi,\xi\big) (x) - g_s(x)\Big\} \, ds,
\end{equation}
where we recall that $g_s(x) := \E\Big[ \Pi^{(b)}\big(e^{-sL}\xi,\xi \big) (x)\Big]$. Consider one of the following time functional space $\mathcal{F}=C_T$, for an arbitrary finite time horizon $T$, or $\mathcal{F}=\mathcal{W}_\lambda$, for some arbitrary $\lambda\geq 1$. Then for every $\alpha\in\Big(1-\frac{d}{4},2-\frac{d}{2}\Big)$ and $p\in(1,\infty)$, we have
$$ 
\E\Big[ \big\|\,\Xi \diamondsuit \xi\,\big\|_{\mathcal{F}\calC^{2\alpha-2}}^p \Big]<\infty.
$$
Moreover, by considering for $\epsilon\in(0,1)$, the regularized versions $\xi^\epsilon:=e^{-\epsilon L} \xi$, and  $\Xi^\epsilon := e^{-\epsilon L} \Xi$, and 
$$ 
C^\epsilon:=\int_0^\infty \E\Big[ \Pi^{(b)}\big(e^{-sL}\xi^\epsilon,\xi^\epsilon\big)\Big]\, ds,
$$ 
then for every $p\in[1,\infty)$, we have
$$ 
\lim_{\epsilon\to 0} \E\Big[ \Big\|\,\Xi \diamondsuit \xi - \Big(\Pi^{(b)}(\Xi^\epsilon,\xi^\epsilon)-C^\epsilon\Big) \Big\|_{\mathcal{F} \calC^{2\alpha-2}}^p \Big]=0.
$$   }
\end{theorem}

\medskip

\noindent Note the following points before proceeding to the proof. 

\begin{itemize}
   \item  In particular, if the ambiant space $M$ has finite measure then the constant weight $\omega\equiv 1$ satisfies \eqref{ass:weight} and belongs to $L^1 \cap L^\infty$. So the previous results can be applied to white noise.  \vspace{0.1cm}
   
   \item In Proposition \ref{prop:whitenoise} as well as in Theorem \ref{thm:Renormalization}, we do not really need $\omega \in L^1 \cap L^\infty$; it suffices that $\omega\in L^p$ for a sufficiently large and finite exponent $p$.   \vspace{0.1cm}
   
   \item Given a point $o\in M$, any weight of the form $\omega(x) = \big(1+d(x,o)\big)^{-M}$, satisfies the assumption \eqref{ass:weight} and belongs to $L^\infty \cap L^2$, provided for $M>\frac{d}{2}$.
\end{itemize} 

\medskip 

\begin{Dem} 
By definition of white noise with weight $\omega$, we know that if $T,T'$ are two self-adjoint operators then for every $y,z\in M$
\begin{equation} 
\E\Big[ \big(T\xi\big)(y)\big(T'\xi\big)(z)\Big] =   \int_M K_T(u,y) K_{T'}(u,z)\,\omega(u) \mu(du). 
\label{eq:rule-bis} \end{equation}
Moreover if $T$ and $T'$ are self-adjoint operators, with a kernel pointwisely bounded by Gaussian kernels at scale $t,t'\in(0,1]$, meaning that we have 
$$
\big| K_T(x,y) \big| \lesssim \G_{t}(x,y) \quad \textrm{and} \quad  \big| K_{T'}(x,y) \big| \lesssim \G_{t'}(x,y)
$$
for almost every $x,y\in M$, then it follows from identity \eqref{eq:rule-bis} and Assumption \eqref{ass:weight} that we have
\begin{align}
\E\big[T\xi(y)T'\xi(z)\big] & = \int  K_T(u,y) K_{T'}(u,z) \omega(u) \, \mu(du) \lesssim \int \G_t(u,y) \, \G_{t'}(u,z) \omega(u) \, \mu(du) \nonumber \\
& \lesssim \omega(y) \int \G_t(y,u) \, \G_{t'}(u,z) \, \mu(du) \lesssim \omega(y)\,\G_{t+t'}(y,z); \label{eq:rule-ter}
\end{align}
we used Lemma \ref{lem:gauss}. Fix now an integer $b\geq 2+{d \over 2}$, and for $r\in(0,1]$ and $s>0$, define the quantity
\begin{align*}
A(r,s):= & \int_0^1 \int_0^1 \left(\frac{r}{r+t_1}\right)^b \left(\frac{r}{r+t_2}\right)^b \ldots \\
& \hspace{2cm} \ldots \left(\frac{t_1t_2}{(t_1+s)(t_2+s)}\right)^{\frac{1}{2}} (s+t_1+t_2)^{-{d \over 2}} (r+t_1+t_2)^{-{d \over 2}} \, \frac{dt_1}{t_1} \frac{dt_2}{t_2},
 \end{align*}
Set
$$
\Theta_s:= \Pi^{(b)} \big(e^{-sL}\xi,\xi\big) - g_s;
$$
we claim that for every $r\in(0,1]$, $s>0$ and every $x\in M$ then we have
\begin{equation}
\E\left[\big|Q_r^{(b)} \Theta_s (x) \big|^2 \right] \lesssim A(r,s) \omega(x)^2.
\label{eq:moment-bis}
\end{equation}

\medskip

\noindent {\bf Step 1 -- Proof of \eqref{eq:moment-bis}.} The resonant, or diagonal, part of the paraproduct $\Pi^{(a)}$ is given by five terms, of the form
$$ 
\mcR^1(f,g) = \int_0^1 \mcP_t \Big( (tL) \mcP_t^1 f  \cdot \mcQ_t g \Big) \, \frac{dt}{t} \quad \textrm{or} \quad \mcR^2(f,g) = \int_0^1 \mcP_t \Big(\mcQ_t f  \cdot (tL)\mcP_t^1 g \Big) \, \frac{dt}{t},
$$
or
$$ 
\mcR^3(f,g) = \int_0^1 \mcP_t \Gamma\left( \sqrt{t} \mcP^1_t f\, ,\, \sqrt{t} \mcP_t^2 g \right) \, \frac{dt}{t}
$$
where 
\begin{itemize}
\item $\mcP_t,\mcP^1_t$ and $\mcP^2_t$ are operators of the form $p(tL)e^{-tL}$ with $p$ a polynomial function;
\item $\mcQ_t$ is of the form $(tL)^{a-1} p(tL)e^{-tL}$ with a polynomial function $p$.
\end{itemize}
So both of these operators have a kernel with Gaussian pointwise estimates and we only need to deal with these three generic quantities.

\ssk

Let us focus on a term of the first form and study
$$
\Theta^1_s:= \mcR^1\big(e^{-sL}\xi,\xi\big) (x) - g^1_s(x) \qquad \textrm{with} \qquad g^1_s(x):=\E\Big[\mcR^1\big(e^{-sL}\xi,\xi\big) (x)\Big].
$$
Due to the covariance rule of Gaussian variables, we have for $T,U,T',U'$ self-adjoint operators (using \eqref{eq:rule-bis}) and every $y,z\in M$
\begin{equation*}
\begin{split}
& \E \Big[ T\xi(y)U\xi(y) T'\xi(z) U'\xi(z) - \E\big[T\xi(y)U\xi(y)\big]\E\big[T'\xi(z)U'\xi(z)\big] \Big] \\
&  \qquad = \E\big[T\xi(y)T'\xi(z)\big]\E\big[U\xi(y)U'\xi(z)\big] + \E\big[T\xi(y)U'\xi(z)]\E[U\xi(y)T'\xi(z)\big].
\end{split}
\end{equation*}
Hence $ \E \left[\big|Q_r^{(b)} \Theta^1_s(x)\big|^2\right]$ is equal to
{\footnotesize 
\begin{align*} 
 &\bigintssss_{0}^1\bigintssss_0^1 \Bigg\{ \bigg(Q_r^{(b)} \mcP_{t_1} \otimes Q_r^{(b)} \mcP_{t_2} \bigg)\bigg(\E\Big[ (t_1 L) \mcP^1_{t_1} e^{-sL}\xi(\bullet) (t_2 L) \mcP^1_{t_2} e^{-sL} \xi(\bullet)\Big] \E\Big[\mcQ_{t_1}\xi(\bullet)\mcQ_{t_2}\xi(\bullet)\Big] \bigg) (x,x)  \\ 
  & + \bigg(Q_r^{(b)} \mcP_{t_1} \otimes Q_r^{(b)} \mcP_{t_2} \bigg)\,\bigg(\E\Big[ (t_1 L) \mcP^1_{t_1} e^{-sL}\xi(\bullet) \mcQ_{t_2}\xi(\bullet)\Big]\, \E\Big[\mcQ_{t_1}\xi(\bullet) (t_2 L) \mcP^1_{t_2} e^{-sL} \xi(\bullet)\Big]\bigg)(x,x)\Bigg\} \, \frac{dt_1}{t_1} \frac{dt_2}{t_2},
\end{align*}} 
\noindent where we use the notation $f(\bullet,\bullet)$ for a function of two variables, with $(fg)(\bullet,\bullet)$ standing for the map $(y,z)\mapsto f(y,z)\,g(y,z)$. Moreover, to shorten notations, we shall use below the notation $dm$ for the measure $\mu(dy) \mu(dz)\frac{dt_1}{t_1} \frac{dt_2}{t_2}$.  By applying \eqref{eq:rule-ter}, it follows that 
$$
\E \left[\big|Q_r^{(b)} \Theta^1_s(x)\big|^2\right] \lesssim J_1+J_2
$$
with 
{\small\begin{align*}
J_1:=\int \Big| K_{Q_r^{(b)} \mcP_{t_1}}(x,y)\Big|\, \Big|K_{Q_r^{(b)} \mcP_{t_2}}(x,z)\Big| \,\omega(y) \,\omega(z) \,\frac{t_1}{t_1+s}\, \frac{t_2}{t_2+s}\, \G_{t_1+t_2+s}(y,z)\, \G_{t_1+t_2}(y,z)\, dm
\end{align*}}
and 
{\small 
\begin{align*} 
J_2:=\int \Big|K_{Q_r^{(b)} \mcP_{t_1}}(x,y)\Big|\,\Big| K_{Q_r^{(b)} \mcP_{t_2}}(x,z)\Big| \,\omega(y) \,\omega(z) \,\frac{t_1}{t_1+s} \,\frac{t_2}{t_2+s} \,\G_{t_1+t_2+s}(y,z)\, \G_{t_1+t_2+s}(y,z)\, dm.
 \end{align*}} 
Let us first explain how we can estimate the kernel of $Q_r^{(b)} \mcP_{t_1}$. Using the notation $\mcP_{t_1}=p(t_1L)e^{-t_1 L}$ for some polynomial function $p$, it comes
\begin{align*}
Q_r^{(b)}\mcP_{t_1} &= \left(\frac{r}{r+\fr{1}{2} t_1}\right)^b \Big((r+\fr{1}{2}t_1)L\Big)^b e^{-r L} p(t_1 L) e^{-t_1L} \\
&= \left(\frac{r}{r+\fr{1}{2}t_1}\right)^b Q_{r+\fr{1}{2}t_1}^{(b)} p(t_1L)e^{-\fr{1}{2}t_1L},
\end{align*}
so since $r+\fr{1}{2}t_1 \simeq r+t_1$, $Q_{r+\fr{1}{2}t_1}^{(b)} $ has a kernel with Gaussian bounds at the scale $r+\fr{1}{2}t_1$ and $p(t_1L)e^{-\fr{1}{2}t_1L}$ at the scale $t_1$, it follows by Lemma \ref{lem:gauss} that $Q_r^{(b)}\mcP_{t_1}$ has a kernel pointwise bounded by $\G_{r+t_1}$ with an extra factor $\left(\frac{r}{r+t_1}\right)^b$. Coming back to estimate the first term $J_1$. We have the upper bound for $J_1$
{\small
\begin{align*} 
&\int \left|\frac{r^2}{(r+t_1)(r+t_2)}\right|^b \G_{r+t_1}(x,y) \,\G_{r+t_2}(x,z) \,\omega(y)\, \omega(z)\, \frac{t_1}{t_1+s} \,\frac{t_2}{t_2+s} \,\G_{t_1+t_2+s}(y,z) \,\G_{t_1+t_2}(y,z)\, dm \\
& \lesssim \omega(x)^2 \int \left|\frac{r^2}{(r+t_1)(r+t_2)}\right|^b \frac{t_1}{t_1+s} \,\frac{t_2}{t_2+s} \,\G_{r+t_1}(x,y) \,\G_{r+t_2}(x,z)\, \G_{t_1+t_2+s}(y,z) \,\G_{t_1+t_2}(y,z)\, dm,
\end{align*}   }
\noindent where we used Assumption \ref{ass:weight}. Due to Lemma \ref{lem:gauss} with Ahlfors regularity \eqref{EqDoublingProperty}, we have
\begin{align*} 
\iint \G_{r+t_1}(x,y) & \G_{r+t_2}(x,z) \,\G_{t_1+t_2+s}(y,z) \,\G_{t_1+t_2}(y,z) \, \mu(dy) \mu(dz) \\ 
&\lesssim (t_1+t_2+s)^{-{d\over 2}} \iint \G_{r+t_1}(x,y) \,\G_{r+t_2}(x,z) \,\G_{t_1+t_2}(y,z) \, \mu(dy) \mu(dz) \\
& \lesssim (t_1+t_2+s)^{-{d \over 2}}\, (r+t_1+t_2)^{-{d \over 2}}.
\end{align*}
Hence,
\begin{align*}
J_1 & \lesssim \omega(x)^2 \int_{0}^1\int_0^1 \left|\frac{r^2}{(r+t_1)(r+t_2)}\right|^b \frac{t_1}{t_1+s} \,\frac{t_2}{t_2+s} \,(s+t_1+t_2)^{-{d \over 2}}\, (r+t_1+t_2)^{-{d \over 2}} \, \frac{dt_1}{t_1} \frac{dt_2}{t_2} \\
       & \lesssim \omega(x)^2 A(r,s). 
\end{align*} 
\noindent The second term $J_2$ can be similarly bounded, which concludes the proof of \eqref{eq:moment-bis} for $\Theta^1$. The corresponding term $\Theta^2$ with $\mcR^2$ can be estimated in the same way. So it remains us now to focus on the last and third term with
$$ 
\mcR^3(f,g) = \int_0^1 \mcP_t\, t\Gamma\left(\mcP^1_t f,\mcP_t^2 g \right) \, \frac{dt}{t}
$$
and $ \Theta^3_s:= \mcR^3\big(e^{-sL}\xi,\xi\big) (x) - g^3_s(x)$. Following the exact same reasoning we have
$$
\E \left(\big|Q_r^{(b)} \Theta^3_s(x)\big|^2\right) \lesssim K_1+K_2
$$
with $K_1$ equal to
{\small
\begin{align*}
\int \Big|K_{Q_r^{(b)} P_{t_1}}(x,y)\Big|\, \Big|K_{Q_r^{(b)} P_{t_2}}(x,z)\Big|\, \omega(y)\, \omega(z) \left|\frac{t_1t_2}{(t_1+s)(t_2+s)}\right|^{\frac{1}{2}} \G_{t_1+t_2+s}(y,z)\, \G_{t_1+t_2}(y,z)\, d\nu
\end{align*}}
and $K_2$ equal to
{\small
\begin{align*} 
\int \Big|K_{Q_r^{(b)} P_{t_1}}(x,y)\Big|\, \Big|K_{Q_r^{(b)} P_{t_2}}(x,z)\Big| \,\omega(y) \omega(z) \left|\frac{t_1t_2}{(t_1+s)(t_2+s)}\right|^{\frac{1}{2}} \G_{t_1+t_2+s}(y,z) \,\G_{t_1+t_2+s}(y,z)\, d\nu.  
 \end{align*}}
Following the same computations, gives us that both $K_1$ and $K_2$ are bounded as follows
$$ 
K_1+K_2 \lesssim \omega(x)^2 A(r,s),
$$
which concludes the proof of \eqref{eq:moment-bis}.
 
\medskip

\noindent {\bf Step 2 -- Conclusion.} We refer the reader to Lemma \ref{lemma:A} for a precise control of quantity $A$. Combining \eqref{eq:moment-bis} with Lemma \ref{lemma:A} gives
\begin{equation}  
\E\left[\big|Q_r^{(b)} \Theta_s(x)\big|\right] \lesssim \E\left[ \big|Q_r^{(b)} \Theta_s(x)\big|^2\right]^{\frac{1}{2}} \lesssim \omega(x) \left(\frac{r}{s+r}\right)^{\frac{1}{2}} (rs)^{-{d \over 4}} \left(1+\log\big(\frac{s+r}{s}\big)\right)^{{1 \over 2}}. 
\label{eq:eq5-bis} 
\end{equation}
We then consider 
$$ 
\big(\Xi \diamondsuit \xi \big)(t) := \int_0^t \Big( \Pi\big(e^{-sL}\xi,\xi\big) (x) - g_s(x)\Big) ds = \int_0^t\Theta_s\,ds.
$$
We refer the reader to Definition \ref{def:besov} for the definition of Besov spaces. For all $0\leq s<t$, it comes (where we forget the low-frequency part in the Besov norm for simplicity, since it is easier than the high-frequency part which we deal with)
\begin{align*}
 & \E \left[\big\| \Xi \diamondsuit \xi (t) - \Xi \diamondsuit \xi (s) \big\|^{2p}_{B^{2\alpha-2}_{2p,2p}}\right] = \int_0^1 r^{-p(2\alpha-2)} \int_M \E \left(\Big| Q_r^{(b)} \big[\Xi \diamondsuit \xi (t) - \Xi \diamondsuit \xi (s)\big]\Big|^{2p} \right) \, \mu(dx)\frac{dr}{r} \\ 
 & \qquad \lesssim  \int_0^1 r^{-p(2\alpha-2)} \int_M  \left(\int_s^t \E\Big(\big| Q_r^{(b)} \Theta(\tau)(x)\big| \Big) \, d\tau\right)^{2p} \, \mu(dx)\frac{dr}{r} \\ 
 & \qquad \lesssim  \|\omega\|_{2p}^{2p} \int_0^1 r^{-p(2\alpha-2)}  \left(\int_s^t \left(\frac{r}{r+\tau}\right)^{\frac{1}{2}} (r\tau)^{-{d \over 4}} \left(1+\log\big(\frac{\tau+r}{\tau}\big)\right)^{1 \over 2} \, d\tau\right)^{2p} \, \frac{dr}{r},
 \end{align*}
where we have used Gaussian hypercontractivity, that is the basic fact that the $L^p$ norms are all equivalent in a fixed Gaussian chaos. So it comes, by Minkowski inequality,
\begin{align*}
& \E \left[\big\| \Xi \diamondsuit \xi (t) - \Xi \diamondsuit \xi (s)\big\|^{2p}_{B^{2\alpha-2}_{2p,2p}}\right] \\
& \qquad \lesssim  \|\omega\|_{2p}^{2p} \left( \int_s^t \left( \int_0^1  r^{-2p(\alpha-1)} \left(\frac{r}{r+\tau}\right)^{p} (r \tau)^{-{{dp} \over 2}} \left(1+\log\big(\frac{\tau+r}{\tau}\big)\right)^{p} \, \frac{dr}{r} \right)^{1 \over {2p}} \, d\tau \right)^{2p}.
\end{align*}
We have
{\small
 \begin{align*}
 & \int_0^1  r^{-2p(\alpha-1)} \left(\frac{r}{r+\tau}\right)^{p} \big(r\tau\big)^{-{{dp}\over 2}} \left(1+\log\big(\frac{\tau+r}{\tau}\big)\right)^{p} \, \frac{dr}{r} \\
 & \qquad \lesssim \int_0^{\tau}  r^{-2p(\alpha-1)} \left(\frac{r}{\tau}\right)^{p} (r \tau)^{-{{dp}\over 2}}  \, \frac{dr}{r} + \int_{\tau}^1  r^{-2p(\alpha-1)} (r\tau)^{-{{dp}\over 2}} \left(1+\log\big(\frac{r}{\tau}\big)\right)^{p} \, \frac{dr}{r}  \\
  & \qquad \lesssim \tau^{-2p(\alpha-1)-dp} + \int_{\tau}^1  r^{-2p(\alpha-1)} (r\tau)^{-{{dp}\over 2}} \left(\frac{r}{\tau}\right)^{p\epsilon} \, \frac{dr}{r} \\
  & \qquad \lesssim \tau^{-2p(\alpha-1)-dp},
 \end{align*}} 
where we chose an arbitrary small parameter $\epsilon$ with
$$ -2(\alpha-1)-\frac{d}{2}+\epsilon<0<-2(\alpha-1)+1-\frac{d}{2}$$ which is equivalent to
$$ 1-\frac{d}{4} +\frac{\epsilon}{2} < \alpha <\frac{3}{2}-\frac{d}{4}$$
and such $\epsilon>0$ exists.
Observe that the latter condition is satisfied since $d\in[2,4)$, so $1-\frac{d}{4} < \alpha< 2-\frac{d}{2} \leq \frac{3}{2}-\frac{d}{4}$. Then because $\alpha<2-\frac{d}{2}$, so $-(\alpha-1)-\frac{d}{2}>-1$, we have in the end
\begin{align*}
\E \left[\big\| \Xi \diamondsuit \xi (t) - \Xi \diamondsuit \xi (s)\big\|^{2p}_{B^{2\alpha-2}_{2p,2p}}\right] & \lesssim \|\omega\|_{2p}^{2p} \left( \int_s^t \tau^{-(\alpha-1)-\frac{d}{2}} \, d\tau \right)^{2p} \\
& \lesssim \|\omega\|_{2p}^{2p} |s-t|^{-2p(\alpha-2)-dp}.
\end{align*}
We can then use Kolmogorov's continuity criterion to deduce that for every $T<\infty$ and $\lambda\geq 1$, we have
$$  
\E \left[\big\| \Xi \diamondsuit \xi \big\|^{2p}_{C_T B^{2\alpha-2}_{2p,2p}}\right] + \E\left[\big\| \Xi \diamondsuit \xi \big\|^{2p}_{\mathcal{E}_\lambda B^{2\alpha-2}_{2p,2p}}\right] < \infty.
$$
As we are in a position to use Besov embedding, as we assume Ahlfors regularity \eqref{EqDoublingProperty}, we know that 
$$ 
B^{2\alpha-2}_{2p,2p} \hookrightarrow B^{2\alpha-2}_{2p,\infty} \hookrightarrow B^{2\alpha-2-\frac{d}{2p}}_{\infty,\infty}=\calC^{2\alpha-2-\frac{d}{2p}}.
$$
So for every $p\geq 1$ and every $\alpha \in \big(1-\frac{d}{2},2-\frac{d}{2}\big)$ 
$$  
\E\left[\big\| \Xi \diamondsuit \xi \big\|^{2p}_{C_T \calC^{2\alpha-2-\frac{d}{2p}}}\right] + \E\left[\big\| \Xi \diamondsuit \xi \big\|^{2p}_{\mathcal{E}_\lambda \calC^{2\alpha-2-\frac{d}{2p}}}\right] < \infty,
$$
which allows us to conclude.
\medskip

\noindent {\bf Step 3 -- Second part of the statement.}  The very same reasonings prove the second part of the statement about the approximation results. We only sketch the arguments. First, observe that estimating the main quantity 
$$
\Xi \diamondsuit \xi - \Big(\Pi^{(b)}(\Xi^\epsilon,\xi^\epsilon)-C^\epsilon\Big)
$$
is almost equivalent to the study of
$$ 
\Xi \diamondsuit \xi - \Xi^\epsilon \diamondsuit \xi^\epsilon.
$$
Then by the bilinear structure, this is equal to
\begin{equation} 
\label{eqaa}
\big(\Xi-\Xi^\epsilon\big) \diamondsuit \xi + \Xi^\epsilon \diamondsuit \big(\xi-\xi^\epsilon\big). 
\end{equation}
We may then repeat the two previous steps to estimate the two terms, independently from one another.
Note that $\Xi-\Xi^\epsilon$ is the resolution of the heat equation associated with $\xi-\xi^\epsilon$, so both of the two terms make appear the quantity $\xi-\xi^\epsilon$.

\smallskip

So to apply the two previous steps to any of the two terms in \eqref{eqaa} is very similar and we only have to include some additional factors coming from
\begin{equation} 
\label{eqaa1}
\xi-\xi^\epsilon = \big(1-e^{-\epsilon L}\big)\xi = \int_0^\epsilon Q_{\sigma}^{(1)} \xi \, \frac{d\sigma}{\sigma}.
\end{equation}
A careful examination shows that if we replace in the previous reasoning (steps 1 and 2) the terms $\xi$ by $Q_{\sigma}^{(1)} \xi$, then all the estimates are identical and we get an extra factor of the forms
$$ 
\frac{s\sigma}{(s+\sigma)^2},\quad \frac{\sigma t_1}{(\sigma+t_1)^2} \qquad \textrm{or} \quad \frac{\sigma t_2}{(\sigma+t_2)^2}.
$$
So in these three situations (by replacing $\xi$ with $Q_{\sigma}^{(1)} \xi$), the same estimates hold with a quantity $A_\sigma(r,s)$ satisfying, for $\eta>0$ as small as we want the inequality
\begin{equation} 
\label{eq:A-bis}
A_\sigma(r,s) \lesssim \left(\frac{r}{s+r}\right) (rs)^{-\frac{d}{2}}\left(1+\log\Big(\frac{s+r}{s}\Big)\right) \left(\frac{\sigma}{\min(r,s)}\right)^\eta,
\end{equation}
instead of estimate \eqref{eq:A}. Such an estimate can be plugged into \eqref{eqaa1} and provides an integrable quantity on an interval $\tau \in(0,\epsilon)$; this gives similar estimates for $\Xi \diamondsuit \xi - \big(\Pi^{(b)}(\Xi^\epsilon,\xi^\epsilon)-C^\epsilon\big)$ as those obtained for $\Xi \diamondsuit \xi$, with an extra factor in terms of $\epsilon^\eta$, and $\eta$, as small as we want, is also appearing in the conditions on $\alpha,p,...$. We then let the reader check that since all the conditions on the exponents are ``open conditions'', then the previous reasoning can be reproduced, up to a small loss of regularity. That means that we can incorporate such a sufficiently small parameter $\eta>0$, and with
$$
F_\epsilon:=  \Xi \diamondsuit \xi - \Big(\Pi^{(b)}(\Xi^\epsilon,\xi^\epsilon) - C^\epsilon\Big)
$$
we can get 
\begin{align*}
\E \left[\big\| F_\epsilon(t) - F_\epsilon(s)\big\|^{2p}_{B^{2\alpha-2}_{2p,2p}}\right] & \lesssim \left( \int_0^\epsilon \int_s^t \left(\frac{\sigma}{\tau}\right)^\eta \tau^{-(\alpha-1) - \frac{d}{2}} \, d\tau \frac{d\sigma}{\sigma}\right)^{2p} \\
& \lesssim  |s-t|^{-2p(\alpha-2)-dp} \left( \frac{\epsilon}{|s-t|} \right)^{2p\eta},
\end{align*}
which allows us to conclude as previously using Besov embedding.
\end{Dem}

\bigskip 

\begin{Dem}[\hspace{0.5cm}{\sf \textbf{Proof of  Theorem \ref{thm:pam} --}}] Fix the weighted white noise $\xi$ and its regularized version $\xi^\epsilon:=e^{-\epsilon L} \xi$. As in Theorem \ref{thm:Renormalization} or Theorem \ref{thm:pam}, set
$$
C^\epsilon :=\int_0^\infty \E\left[ \Pi(e^{-sL}\xi^\epsilon,\xi^\epsilon)\right]\, ds = \E\Big[\Pi(L^{-1}\xi^\epsilon,\xi^\epsilon)\Big].
$$
In order to make appear this term in the equation, we can introduce a suitable correction term in the regularized problems, which leads us to study the following renormalized PDE
$$ 
\partial_t u^\epsilon + L u^\epsilon = F(u^\epsilon)\xi^\epsilon - C^\epsilon F'(u^\epsilon) F(u^\epsilon);
$$
this is elementary, and detailed in \cite{GIP}. We can then follow the exact same approach as for Theorem \ref{thm:gpam}, or Theorem \ref{thm:gpam-2} for the global estimates with the spaces $\mathcal{E}_\lambda$, adapted to this modified PDE. So we only detail the required modification.
We cannot use $\Pi^{(b)}(X^\epsilon,\xi^\epsilon )$, which does not converge, and use instead $\Pi^{(b)}(X^\epsilon,\xi^\epsilon ) - C^\epsilon$, which converges in $C_T\calC^{2\alpha-2}$.

\ssk

Using the arguments of Theorem \ref{thm:prod}, with $\gamma=\alpha-2$, one defines a new ``product'' for $(u,u')\in {\sf C}_\alpha^\alpha(X^ \epsilon)$ 
$$ 
\Big(\big(F(u),F'(u)u'\big) \cdot \xi^\epsilon\Big) - C^\epsilon u'F'(u)\in {\sf C}^{\alpha'}_{\alpha-2,T}\big(\xi^\epsilon\big)
$$
(whose the norms can be estimated by quantities independent with respecto $\epsilon$),
for which we have the uniform estimate 
{\small   \begin{align*} 
&\Big\|\big(F(u),F'(u) u'\big)\cdot \xi^\epsilon - \Pi^{b}_{F(u)} (\xi^\epsilon) - C^\epsilon u'F'(u)\Big\|_{C_T\calC^{2\alpha-2}} \\ 
&\lesssim \Big\|\big(F(u),F'(u)u'\big)\Big\|_{{\sf C}^{\alpha'}_{\alpha,[0,T]}(X)}\left(\|\xi^\epsilon\|_{\calC^{\alpha-2}}+\|X^\epsilon\|_{C_T\calC^\alpha}\|\xi^\epsilon\|_{\calC^{\alpha-2}} + \Big\|\Pi^{(b)}\big(X^\epsilon,\xi^\epsilon\big) - C^\epsilon\Big\|_{C_T\calC^{2\alpha-2}}\right).
\end{align*}   } 
We conclude as in the proofs of Theorems \ref{thm:gpam} and \ref{thm:gpam-2}, by a fixed point argument. 
\end{Dem}


\appendix

\section[\hspace{0.6cm} {\sf Heat kernel and technical estimates}]{Heat kernel and technical estimates} 
\label{Appendix}
 
We gather in this Appendix a number of propositions whose proofs were not given in the course of the paper, so as to keep focused on the most essential aspects of our work. These proofs are given here.

\medskip
 
We start by proving the following pointwise and $L^p$-estimate for the gradient of the heat semigroup. 
 
\medskip 
 
\begin{proposition}   {\sf
Assume that $(M,d,\mu)$ is a doubling space equipped with a semigroup satisfying \eqref{UE} and \eqref{Lipschitz}. Then for every $t>0$, $x_0\in M$ and every function $f\in L^2$ we have 
\begin{equation}  \label{eq:gradientpt}
\Big|\big(\sqrt{t} \Gamma\big)\big(e^{-tL} f\big)(x_0)\Big| \lesssim \bigintsss_M \frac{1}{\sqrt{V\big(x_0,\sqrt{t}\big)V\big(y,\sqrt{t}\big)}}\,\exp\left(-c\,\frac{d(x_0,y)^2}{t}\right) \big|f(y)\big| \, \mu(dy).
\end{equation}
}  
\end{proposition}

\medskip

Let us first introduce the following notation: for a function $f\in L^2_{\textrm{loc}}$ and a ball $B\subset M$, we write $\osc_B(f)$ for the $L^2$ oscillation of $f$ on $B$ defined by 
$$
\osc_B(f) := \left(\aver{B} \left|f - \aver{B} f \, d\mu \right|^2 \, d\mu \right)^\frac{1}{2},
$$
where $\aver{B} f \, d\mu$ stands for the average of $f$ on the ball $B$.

\medskip

\begin{Dem} 
Fix the function $f\in L^2$ and consider $g=e^{-tL}f$. By $L^2$-Caccioppoli inequality (see Lemma below), we have for every $x_0$ and $r>0$ that
\begin{align*}
\left(\aver{B(x_0,r)} \Gamma(g)^2 \, d\mu \right)^{\frac{1}{2}} &\lesssim \frac{1}{r}  \osc_{B(x_0,2r)}(g)  + \left(\aver{B(x_0,2r)} \left|L g\right|^2 \, d\mu \right)^\frac{1}{2}.
\end{align*}
So if $x_0$ is a Lebesgue point of $\Gamma(g)^2$ and $|Lg|^2$ (which is the case for almost every point $x_0 \in M$) then taking the limit for $r\to 0$ yields
\begin{equation} 
\label{eq:i0}
\Gamma(g)(x_0) \lesssim \liminf_{r\to 0} \frac{1}{r}  \osc_{B(x_0,2r)}(g) + \big|Lg(x_0)\big|.  
\end{equation}
Since $(tL)e^{-tL}$ has a kernel satisfying the Gaussian upper estimates \eqref{UE} (by analyticity), we deduce that
\begin{equation} 
\label{eq:i1}
\big|tLg(x_0)\big| \lesssim \int_M \frac{1}{\sqrt{V(x_0,\sqrt{t})V(y,\sqrt{t})}}\exp
\left(-c \,\frac{d(x_0,y)^2}{t}\right) \big|f(y)\big| \, \mu(dy).  
\end{equation}
Using Lipschitz regularity \eqref{Lipschitz} for the heat kernel and doubling property, it comes for $x,z\in B(x_0,2r)$ with $r\leq \sqrt{t}$
\begin{align*}
\big|g(x)-g(z)\big| & = \Big|e^{-tL}f(x)-e^{-tL}f(z)\Big|   \\
& \lesssim \left(\frac{d(x,z)}{\sqrt{t}}\right)  \int_M \frac{1}{\sqrt{V\big(x_0,\sqrt{t}\big)V\big(y,\sqrt{t}\big)}} \, \exp
\left(-c \,\frac{d(x,y)^2}{t}\right) \big|f(y)\big| \, \mu(dy).
\end{align*}
Hence, uniformly with respect to $r\in(0,\sqrt{t})$ we obtain
\begin{align} 
\label{eq:i2}
\frac{1}{r}\osc_{B(x_0,2r)}(g) \lesssim  \frac{1}{\sqrt{t}}  \int_M \frac{1}{\sqrt{V\big(x_0,\sqrt{t}\big)V\big(y,\sqrt{t}\big)}} \, \exp
\left(-c \,\frac{d(x,y)^2}{t}\right) \big|f(y)\big| \, \mu(dy).  
\end{align}
By combining these last inequalities \eqref{eq:i1} and \eqref{eq:i2} into \eqref{eq:i0}, one concludes to \eqref{eq:gradientpt}.
\end{Dem}

\medskip

\begin{lemma}[Cacciopoli inequality] \label{lem:caccio}   {\sf
For every ball $B$ of radius $r>0$ and every function $f\in {\mathcal D}_2(L)$ we have
$$ 
\left(\aver{B} \Gamma(f)^2 \, d\mu \right)^{\frac{1}{2}} \lesssim \frac{1}{r} \osc_{2B}(f) + r \left(\aver{2B} \left| Lf \right|^2 \, d\mu \right)^{\frac{1}{2}}.
$$   }
\end{lemma}

\medskip

Before to check this inequality, let us first recall some consequences of the Gaussian upper estimates \eqref{UE}. Under \eqref{UE}, we know that a scale-invariant local Sobolev inequality holds, more precisely
$$ 
\|f\|_{q}^2 \lesssim |B|^{\frac{2}{q}-1}\left( \|f\|_2^2 + r^2 {\mathcal E}(f,f)\right),
$$
for every ball $B$ of radius $r>0$, every $f\in {\mathcal D}_2(\Gamma)$ supported in $B$ and for some $q > 2$. This inequality was introduced in \cite{S} and was shown, under \eqref{d}, to be equivalent  to \eqref{UE} in the Riemannian setting. The equivalence was stated in  our more general setting in \cite{ST2}. See also \cite{BCS} for many reformulations of local Sobolev inequalities, an alternative proof of the equivalence with \eqref{UE}, and more references.

Such a local Sobolev inequality also implies a following relative Faber-Krahn inequality (see for instance \cite[Theorem 2.5]{HS}, as well as \cite[Section 3.3]{BCS}): for every ball $B$ with a small enough radius $r>0$ , every function $f\in {\mathcal D}_2(\Gamma)$ supported in $B$ then
\begin{equation} 
\label{eq:FK}
\|f\|_{2} \lesssim r \big\| \Gamma(f)\big\|_2.  
\end{equation}

\medskip

\begin{Dem}[\hspace{0.5cm}{\sf \textbf{Proof of Lemma \ref{lem:caccio} --}}]
We refer to \cite[Lemma A.1]{BCF1} for such a result for harmonic function: if $u\in{\mathcal D}_2(L)$ is harmonic on $2B$ (which means $L(u)=0$ on $2B$) then
\begin{equation} \left(\aver{B} \Gamma(u)^2 \, d\mu \right)^{\frac{1}{2}} \lesssim \frac{1}{r} \osc_{2B}(u). \label{eq:caccioh} \end{equation}
Now consider $f\in {\mathcal D}_2(L)$. By \cite[Lemma 4.6]{BCF1}, it is known that there exists $u\in{\mathcal D}_2(L)$ harmonic on $2B$ such that $f-u\in{\mathcal D}(\Gamma)$ is supported on the ball $2B$.
By the support property, it follows
$$ 
\big\| \Gamma(f-u) \big\|_2^2 = \int (f-u) L(f-u) \, d\mu = \int  (f-u) L(f) \, d\mu.
$$
So using  Faber-Krahn inequality \eqref{eq:FK} we obtain
\begin{equation} 
\label{eq:FK2}
\| f-u \|_2 \lesssim  r \big\|\Gamma(f-u)\big\|_{L^2(2B)}  
\end{equation}
and so
$$ 
\big\| \Gamma(f-u) \big\|_2^2 \lesssim \|f-u\|_2 \big\|L(f)\big\|_{L^2(2B)} \lesssim r \big\|\Gamma(f-u)\big\|_2 \|L(f)\|_{L^2(2B)},
$$
which yields
\begin{equation} 
\label{eq:ga} 
\big\| \Gamma(f-u) \big\|_2 \lesssim r\big\|L(f)\big\|_{L^2(2B)}. 
\end{equation}
Then we split
$$
\big\| \Gamma(f) \big\|_2 \leq \big\| \Gamma(f-u)\big\|_2 + \big\| \Gamma(u)\big\|_2 \lesssim r \big\|L(f)\big\|_{L^2(2B)} + \big\|\Gamma(u)\big\|_2
$$
and then use \eqref{eq:caccioh} to get
\begin{align*}
\left( \aver{B}  \Gamma(f)^2 \, d\mu \right)^{\frac{1}{2}} & \lesssim r \left( \aver{2B} |L(f)|^2 \, d\mu \right)^{\frac{1}{2}} + \frac{1}{r} \osc_{2B}(u) \\
 &  \lesssim r \left( \aver{2B} |L(f)|^2 \, d\mu \right)^{\frac{1}{2}} + \frac{1}{r} \osc_{2B}(f) + \frac{1}{r} \left(\aver{2B} | f-u|^2 \, d\mu \right)^{\frac{1}{2}} \\
 & \lesssim r \left( \aver{2B} |L(f)|^2 \, d\mu \right)^{\frac{1}{2}} + \frac{1}{r} \osc_{2B}(f),
\end{align*}
where we used again \eqref{eq:FK2} and \eqref{eq:ga} at the last step.
\end{Dem}
 
\medskip
 
We also give a proof of the following basic important fact about the H\"older spaces $\calC^\sigma$.
 
\medskip 
 
\begin{proposition}   {\sf
For $\sigma<2$, the H\"older spaces $\calC^\sigma$ do not depend on the parameter $a$ used to define them, and the two norms on $\calC^\sigma$ corresponding to two different parameters $a,a'$, are equivalent.   } 
\end{proposition}
 
\medskip
  
\begin{Dem}
Given two positive integers $a$ and $a'$, consider the two spaces $\calC^ \sigma_a$ and $\calC^ \sigma_{a'}$, and their corresponding norms. Fix $t\in(0,1]$.
If $a'\geq a$, then writing
$$ 
Q_t^{(a')} =2^{a'} Q_{\frac{t}{2}}^{(a')} e^{-\,\frac{t}{2}\,L} = 2^{a'}  Q_{\frac{t}{2}}^{(a)} Q_{\frac{t}{2}}^{(a'-a)}
$$
and using the fact that the operators $Q_{\frac{t}{2}}^{(a'-a)}$ are uniformly bounded on $L^\infty$, we get
$$ 
\|\cdot\|_{\calC^ \sigma_{a'}} \lesssim \|\cdot\|_{\calC^ \sigma_a}.
$$
If now $a'< a$, write
$$ 
Q_t^{(a')} =  \frac{1}{\gamma_{a-a'}}\,\int_0^1 Q_t^{(a')} Q_s^{(a-a')} \, \frac{ds}{s} + Q_t^{(a')} P_1^{(a-a')}.
$$
For $s\leq t$, we have
$$ 
Q_t^{(a')} Q_s^{(a-a')} = \left(\frac{s}{t}\right)^{a-a'} Q_{t+s}^{(a)} \left(\frac{t}{t+s}\right)^a
$$
so that for $f\in \calC^\sigma_a$
\begin{align*}
\Big\| Q_t^{(a')} Q_s^{(a-a')} f \Big\|_{\infty } & \lesssim \left(\frac{s}{t}\right)^{a-a'} \Big\| Q_{t+s}^{(a)}f \Big\|_{\infty} \\
  & \lesssim \left(\frac{s}{t}\right)^{a-a'}  t^{\sigma\over 2}\,\|f \|_{\calC^ \sigma_a}.
\end{align*}
For $t\leq s$, we have
$$ 
Q_t^{(a')} Q_s^{(a-a')} =  \left(\frac{s}{t}\right)^{a-a'} Q_{t+s}^{(a)} \left(\frac{t}{t+s}\right)^a
$$
so that 
\begin{align*}
\Big\|Q_t^{(a')} Q_s^{(a-a')} f \Big\|_{\infty } & \lesssim \left(\frac{t}{s}\right)^{a'} \Big\| Q_{t+s}^{(a)} f \Big\|_{\infty} \\
& \lesssim \left(\frac{t}{s}\right)^{a'}  s^{\sigma\over 2}\,\|f \|_{\calC^\sigma_a},
\end{align*}
and similarly
\begin{align*}
\Big\|Q_t^{(a')} P_1^{(a-a')} f \Big\|_{\infty } \lesssim t^{a'} \|f \|_{\calC^\sigma_a}.
\end{align*}
Then by integrating (and since $a', a-a'\geq 1>\frac{\sigma}{2}$) we have
\begin{align*}
\Big\| Q_t^{(a')}f \Big\|_{\infty} & \lesssim \left(\int_0^t \left(\frac{s}{t}\right)^{a-a'} \, \frac{ds}{s}\right)  t^{\sigma\over 2}\,\|f \|_{\calC^ \sigma_a} + \left(\int_t^1 \left(\frac{t}{s}\right)^{a'}  s^{\sigma\over 2} \, \frac{ds}{s}\right) \|f \|_{\calC^ \sigma_a} + t^{a'}\,\|f \|_{\calC^ \sigma_a} \\
& \lesssim t^{\sigma\over 2} \|f \|_{\calC^ \sigma_a},
\end{align*}
which concludes the proof that
$$ 
\|\cdot\|_{\calC^ \sigma_{a'}} \lesssim \|\cdot\|_{\calC^ \sigma_a}.
$$
\end{Dem} 

\medskip

The following lemma provides a useful way of proving that a distribution is H\"older; it was used in sections \ref{SubsectionHolder} and \ref{SubsectionRenormalization} to investigate the almost sure regularity properties of white noise and the renormalized paraproduct dealt with in theorem \ref{thm:Renormalization}. Recall that Besov spaces were defined in Definition \ref{def:besov}.
 
\medskip 
 
\begin{lemma}  \label{lem:Besov}    {\sf
Assume that the metric measure space $(M,d,\mu)$ is Ahlfors regular (see \eqref{EqDoublingProperty}), with exponent $\nu$. Then, given $-\infty<\sigma<2$, and $1<p<\infty$, we have the continuous embeddings
$$ 
B^{\sigma}_{p,p} \hookrightarrow B^{\sigma}_{p,\infty} \hookrightarrow B^{\sigma - \frac{\nu}{p}}_{\infty,\infty}=\calC^{\sigma - \frac{\nu}{p}}.
$$   } 
\end{lemma}

\medskip

\begin{Dem}
The first embedding is a direct application of the following fact. For $s\in(0,1)$ and an integer $a\geq 2$ then
$$
 Q^{(a)}_{s} f = \frac{2}{s} \int_{\frac{s}{2}}^s Q^{(a)}_{t} \left(\frac{s}{t}\right)^a e^{-(s-t)L} f \, dt.
 $$
Since the semigroup is uniformly bounded on $L^p$, we get
$$
\Big\|Q^{(a)}_sf \Big\|_p \lesssim \int_{\frac{s}{2}}^s \Big\| Q^{(a)}_t f\Big\|_p \, \frac{dt}{t} 
$$
and by H\"older inequality
$$
\Big\|Q^{(a)}_sf \Big\|_p \lesssim \left(\int_{\frac{s}{2}}^s \Big\| Q^{(a)}_t f\Big\|_p^p \, \frac{dt}{t}\right)^{\frac{1}{p}} \lesssim s^{\sigma \over 2} \|f\|_{B^{\sigma}_{p,p}}. 
$$
The second embedding comes from the following elementary fact. For $t\in(0,1)$, let $T$ a linear operator with a kernel, pointwisely bounded by a Gaussian kernel $\G_t$ at scale $t$, then with Ahlfors regularity \eqref{EqDoublingProperty}, we have
$$
\| T \|_{L^p\to L^\infty} \lesssim t^{-\frac{\nu}{2p}}.
$$
So for $s\in \big(0,\frac{1}{2}\big)$, applying to $T=Q_s^{(a)}$ we obtain since $Q^{(2a)}_{2s} = 2^{2a} Q^{(a)}_s Q_s^{(a)}$
$$
\Big\| Q^{(2a)}_{2s} f \Big\|_{\infty} \lesssim s ^{-\frac{\nu}{2p}}\,\Big\| Q_s^{(a)} f \Big\|_{L^p} \lesssim s^{-\frac{\nu}{2p}+\frac{\sigma}{2}} \, \|f\|_{B^{\sigma}_{p,\infty}},
$$
which proves the embedding $B^{\sigma}_{p,\infty} \hookrightarrow B^{\sigma - \frac{\nu}{p}}_{\infty,\infty}$.
\end{Dem} 

\bigskip

\noindent The next three lemmas were used in the proof of the renormalization theorem \ref{thm:Renormalization}.

\medskip

\begin{lemma}  \label{lem:gauss}     {\sf
For every $t>0$, set $\G_t$ the Gaussian kernel at scale $t$
$$ 
\G_t(x,y):= \frac{1}{V(x,\sqrt{t})} e^{-c\,\frac{d(x,y)^2}{t}},
$$
where we forget the dependence with respect to the constant $c$ in the notation. Then for $s>t>0$ and every $x,z\in M$, we have
$$ 
\int \G_s(x,y) \G_t(y,z) \, \mu(dy) \lesssim \G_s(x,z)\lesssim \frac{1}{V(x,\sqrt{s})+V(z,\sqrt{s})}.
$$   }
\end{lemma}

\medskip

\begin{Dem} 
By considering that $\G_t(\cdot,z)$ belongs to $L^1$, uniformly in $s$, we directly obtain that
$$ 
\int \G_s(x,y) \G_t(y,z) \, \mu(dy) \lesssim \sup_{y\in M} \ \G_s(x,y) \lesssim \frac{1}{V(x,\sqrt{s})}.
$$
Moreover, 
\begin{align*}
\exp\left(-\frac{d(x,y)^2}{s}\right) \cdot \exp\left(-\frac{d(y,z)^2}{t}\right) & \leq \exp\left(-\frac{d(x,y)^2}{s}\right) \cdot \exp\left(-\frac{d(y,z)^2}{s}\right) \\
& \leq \exp\left(-\frac{d(x,z)^2}{2s}\right). 
\end{align*}
So in the product $\G_s(x,y) \G_t(y,z)$, we may factorize an exponential decay and so for some implicit constants, we have
$$ 
\int \G_s(x,y) \G_t(y,z) \, \mu(dy) \lesssim \G_{s}(x,z). 
$$
\end{Dem}

\medskip

\begin{lemma} \label{lemma:A}  {\sf 
For $r\in(0,1]$, $s>0$ and $d\geq 2$, let us consider the quantity
\begin{align*}
A(r,s):= & \int_0^1 \int_0^1 \left(\frac{r}{r+t_1}\right)^b \left(\frac{r}{r+t_2}\right)^b \ldots \\
 & \hspace{2cm} \ldots \left(\frac{t_1t_2}{(t_1+s)(t_2+s)}\right)^{\frac{1}{2}} (s+t_1+t_2)^{-{d\over 2}} (r+t_1+t_2)^{-{d \over 2}} \, \frac{dt_1}{t_1} \frac{dt_2}{t_2},
 \end{align*}
 where $b\geq 2+{d \over 2}$ is an integer. 
Then we have
\begin{equation} 
\label{eq:A}
A(r,s) \lesssim \left(\frac{r}{s+r}\right) (rs)^{-{d \over 2}} \left(1+\log\Big(\frac{s+r}{s}\Big)\right).  
\end{equation}   }
\end{lemma}

\medskip

\begin{Dem} The two variables $t_1,t_2$ play a symmetric role so we may restrict our attention to the double integral under the condition $t_2\leq t_1$.
The part $A_1$ of the double integral where $t_2\leq t_1 \leq r$ gives
\begin{align*}
A_1(r,s) & = \int_0^r \int_0^{t_1} \left(\frac{t_1 t_2}{(t_1+s)(t_2+s)}\right)^{\frac{1}{2}} r^{-{d \over 2}} (t_1+s)^{-{d \over 2}} \, \frac{dt_2}{t_2} \frac{dt_1}{t_1} \\
&  \lesssim \int_0^r \left(\frac{t_1}{t_1+s}\right)^{\frac{1}{2}} (t_1+s)^{-{d \over 2}} r^{-{d \over 2}} \left(\frac{t_1}{t_1+s}\right)^{\frac{1}{2}} \left(1+\log\big(\frac{s+t_1}{s}\big)\right) \, \frac{dt_1}{t_1} \\
& \lesssim  r^{-{d \over 2}} \int_0^r \left(\frac{t_1}{(t_1+s)^{d/2+1}}\right) \left(1+\log\big(\frac{s+t_1}{s}\big) \right)  \, \frac{dt_1}{t_1} \\
& \lesssim r^{-{d \over 2}} \left\{s^{-d/2} - \left(1+\log\big(\frac{s+r}{s}\big)\right)(s+r)^{-{d \over 2}}\right\} \\
& \lesssim (sr)^{-{d \over 2}} \frac{r}{s+r} \left(1+\log\big(\frac{s+r}{s}\big)\right),
\end{align*}
where we used the basic inequality
$$ 
\int_0^{t_1} \frac{dt_2}{\sqrt{t_2(t_2+s)}} \lesssim \left(\frac{t_1}{s+t_1}\right)^{\frac{1}{2}}\left(1+ \log\Big(\frac{s+t_1}{s}\Big)\right)
$$
which can be easily checked by splitting into the two cases $t_1\leq s$ and $s\leq t_1$. The second part $A_2$ of the double integral where $t_2\leq r \leq t_1$ is controlled as follows
 \begin{align*} 
A_2(r,s) & = \int_r^1 \int_0^{r} \left(\frac{r}{r+t_1}\right)^b \left(\frac{t_1 t_2}{(t_1+s)(t_2+s)}\right)^{\frac{1}{2}} t_1^{-{d \over 2}} (t_1+s)^{-{d \over 2}} \, \frac{dt_2}{t_2} \frac{dt_1}{t_1} \\
& \lesssim \int_r^1  \left(\frac{r}{t_1}\right)^b \left(\frac{t_1}{t_1+s}\right)^{\frac{1}{2}} \left(\frac{r}{r+s}\right)^{\frac{1}{2}} \left(1+\log\Big(\frac{s+r}{s}\Big)\right)  (t_1+s)^{-{d \over 2}} t_1^{-{d \over 2}} \, \frac{dt_1}{t_1} \\
& \lesssim  \frac{r}{r+s} \big(r(r+s)\big)^{-{d \over 2}}  \left(1+\log\Big(\frac{s+r}{s}\Big) \right),
\end{align*}
where we used Lemma \ref{lem:lem-bis}. The third and last part $A_3$ of the double integral where $r\leq t_2 \leq t_1$ satisfies
 \begin{align*}
A_3(r,s) & = \int_r^1 \int_r^{t_1} \left(\frac{r}{t_1}\right)^b \left(\frac{r}{t_2}\right)^b \left(\frac{t_1t_2}{(t_1+s)(t_2+s)}\right)^{\frac{1}{2}} t_1^{-{d \over 2}} (t_1+s)^{-{d \over 2}} \, \frac{dt_2}{t_2} \frac{dt_1}{t_1} \\
 & \lesssim \int_r^1  \left(\frac{r}{t_1}\right)^b \left(\frac{t_1}{t_1+s}\right)^{\frac{1}{2}} \left(\frac{r}{r+s}\right)^{\frac{1}{2}}   (t_1+s)^{-{d \over 2}} t_1^{-{d \over 2}} \, \frac{dt_1}{t_1}.
 \end{align*}
We then use again Lemma \ref{lem:lem-bis}, to obtain
 \begin{align*}
A_3(r,s) & \lesssim  \frac{r}{r+s} (rs)^{-{d \over 2}}.
 \end{align*} 
Inequality \eqref{eq:A} comes by combining the above three estimates.
\end{Dem}

\medskip

\begin{lemma}  \label{lem:lem-bis}    {\sf 
For every $r,t\in(0,1)$ and any $0<\rho<\epsilon$
$$ 
\int_r^1 \left(\frac{t}{t+s}\right)^\rho t^{-\epsilon} \frac{dt}{t} \lesssim \left(\frac{r}{s+r}\right)^\rho r^{-\epsilon}. 
$$   }
\end{lemma}
 
\medskip 
 
\begin{Dem}
Indeed if $r\geq s$ then for every $t\in(r,1)$ we have $t \simeq t+s$ and so
\begin{align*}
\int_r^1 \left(\frac{t}{t+s}\right)^\rho t^{-\epsilon} \frac{dt}{t} & \lesssim \int_r^1  t^{-\epsilon} \frac{dt}{t} \simeq  r^{-\epsilon}.
\end{align*}
Now if $s\geq r$, we split the integral in two terms and we have
\begin{align*}
\int_r^1 \left(\frac{t}{t+s}\right)^\rho t^{-\epsilon} \frac{dt}{t} & \lesssim \int_r^s \left(\frac{t}{s}\right)^\rho t^{-\epsilon} \frac{dt}{t} + \int_s^1 t^{-\epsilon} \frac{dt}{t} \\
& \lesssim \left(\frac{r}{s}\right)^\rho r^{-\epsilon} + s^{-\epsilon} \lesssim \left(\frac{r}{s}\right)^\rho r^{-\epsilon},
\end{align*} 
where we used $\epsilon>A \rho$.
\end{Dem}

\bigskip

\section[\hspace{0.6cm} {\sf Extension of the theory}]{Extension of the theory} 
\label{part2}

Consider as above a doubling metric measure space $(M,d,\mu)$ equipped with a heat semigroup satisfying the upper gaussian estimates \eqref{UE}. We aim in this appendix at explaining how one can get the same conclusions as in the above main body of work
\begin{enumerate}
\item by weakening the Lipschitz regularity assumption on the heat kernel \eqref{Lipschitz}, assuming only some integrated estimates of the gradient of the heat kernel; \vspace{0.1cm}

\item by developing the theory of paracontrolled calculus in Sobolev spaces rather than in H\"older spaces. By Sobolev embedding, Sobolev spaces are included in some H\"older spaces, so it will be interesting to understand if starting from an initial data belonging to some Sobolev space, the solution of renormalized singular PDEs will lives in this same scale of Sobolev spaces. From a technical point of view, it is a bit more difficult since Sobolev spaces involve simultaneously all the frequencies, whereas for H\"older spaces we can work at a fixed frequency scale.
\end{enumerate}

\ssk

We give in section \ref{SubsectionAppRegularityAssumptions} the regularity assumptions on the heat kernel under which we shall work here, and reformulate and extend in section \ref{SubsectionAppFunctionalCalculus} the main continuity estimates on the operators $P^{(a)}_t, Q^{(a)}_t$ and $\Gamma$ needed to extend the paraproduct machinery to the present setting. The latter, together with some crucial commutator estimates in H\"older and Sobolev spaces, is investigated in section \ref{SubsectionAppParaproduct}. The last and short section \ref{SubsectionAppPAM} describes how these results can be used to extend the results of section \ref{sec:pam} to our optimal regularity setting. 

\medskip

This appendix was written jointly with Dorothee Frey.

\bigskip

\subsection[\hspace{-1cm} {\sf Regularity assumptions}]{Regularity assumptions}
\label{SubsectionAppRegularityAssumptions}

Rather than assuming the Lipschitz property (Lip) used above we shall assume here that the gradient /carr\'e du champ operator $\Gamma$ satisfies some $L^q$ estimates and the $L^q$-de Giorgi property recalled below in sections \ref{SubsubsectionLqEstimatesGradient} and \ref{SubsubsectionDeGiorgi}. We shall also assume that it satisfies a scale-invariant Poincar\'e inequality recalled in section \ref{SubsubsectionPoincare}.

\subsubsection{$L^q$-estimates of the gradient of the semigroup}
\label{SubsubsectionLqEstimatesGradient}

Given $q_0>2$, the uniform $L^{q_0}$-boundedness of the gradient (or ``carr\'e du champ") of the semigroup was introduced in \cite{ACDH}
\begin{equation} \label{Gp}
\sup_{t>0} \|\sqrt{t}\Gamma e^{-tL} \|_{q_0\to q_0} <+\infty. \tag{$\textrm{G}_{p_0}$}
\end{equation}

By definition of the carr\'e du champ operator,  $(\textrm{G}_2)$ holds trivially. It is known in that case that this global $L^2$-inequality can be improved into localized estimates, via $L^2$-Davies-Gaffney estimates. For every subset $E,F \subset M$ and every $t>0$, we have
$$ 
\big\|e^{-tL}\big\|_{L^2(E) \to L^2(F)} + \sqrt{t} \big\|\Gamma e^{-tL}\big\|_{L^2(E) \to L^2(F)} \lesssim e^{-c \frac{d^2(E,F)}{t}}.
$$
Assuming the volume doubling condition (VD) and the Gaussian upper bound (UE) for the heat kernel, one can interpolate the estimate $(\textrm{G}_{q_0})$ with the above $L^2$-Davies-Gaffney estimates and deduce that $(\textrm{G}_q)$ holds for every $q\in[2,q_0]$. More precisely, for every subset $E,F \subset M$ and every $t>0$, the inequality
$$ 
\big\|e^{-tL}\big\|_{L^q(E) \to L^q(F)} + \sqrt{t} \big\|\Gamma e^{-tL}\big\|_{L^q(E) \to L^q(F)} \lesssim e^{-c_q \frac{d^2(E,F)}{t}}
$$
holds for some positive constant $c_q$, only depending on $q\in [2,q_0)$. Following \cite[Proposition 1.10]{ACDH}, the latter estimate can be reformulated in terms of integral estimates of the gradient of the heat kernel. Denoting by $p_t$ the kernel of $e^{-tL}$, we have
\begin{equation*}
\sqrt{t} \big\| \Gamma_x\,p_t(.,y) \big\|_q \lesssim \left[ V\big(y,\sqrt{t}\big)\right]^{-(1-\frac{1}{q})},
\end{equation*} 
for $\mu$-almost all $y\in M$ and all positive times. By interpolation with the $L^2$-Gaffney estimates, there exists a positive constant $c$ such that
\begin{equation}\label{gradlp}
\sqrt{t} \Big\| e^{c\frac{d(\cdot,y)^2}{t}} \Gamma_x\,p_t(.,y) \Big\|_q \lesssim \left[ V\big(y,\sqrt{t}\big)\right]^{-(1-\frac{1}{q})}
\end{equation}
holds for $\mu$-almost all $y\in M$ and all positive times.
 
We refer the reader to \cite{ACDH} for more details about Property $(\textrm{G}_p)$ and the link with the boundedness of the Riesz transform; see also \cite{BF2} and references therein for more details.

\medskip

\subsubsection{$L^q$-de Giorgi Property}
\label{SubsubsectionDeGiorgi}

The so-called "de Giorgi property'',  or "Dirichlet property", on the growth of the Dirichlet integral for  harmonic functions was introduced by De Giorgi in \cite{DeG}, for second order divergence form differential operators on $\R^n$, with real coefficients. In de Giorgi's work, this property prescribes a(n at most) linear growth rate for the $L^2$-average of gradients of harmonic functions. This property was subsequently used in many works  and in various situations in order to prove H\"older regularity for  solutions of inhomogeneous elliptic equations and systems. An $L^q$-version was recently introduced in \cite{BCF1}, and we refer the reader to that work for more details about it.

\medskip

\begin{definition*}[$L^q$-de Giorgi property] {\sf 
Given $q\in[1,+\infty)$ and $\theta\in(0,1)$, we say that the operator $\Gamma$ satisfies the inequality \eqref{DG} if it satisfies the following estimate. For every positive $r\leq R$, every pair of concentric balls $B_r,B_R$ with radii $r$ and $R$, respectively, and for every function $f\in {\mathcal D}$, one has
\begin{equation}
\left(\aver{B_r} |\Gamma f|^q d\mu \right)^\frac{1}{q} \lesssim \left(\frac{R}{r}\right)^\theta \left\{\left(\aver{B_R} |\Gamma f|^q d\mu\right)^\frac{1}{q} + R \| L f \|_{L^\infty(B_R)} \right\}.
\tag{$D\textrm{G}_{q,\theta}$} \label{DG}
\end{equation}
We sometimes omit the parameter $\theta$, and write $(D\textrm{G}_q)$ if \eqref{DG} is satisfied for some $\theta \in (0,1)$.}
\end{definition*}

\medskip

As we always have 
$$ 
\left(\aver{B_r} |\Gamma f|^q d\mu \right)^\frac{1}{q} \lesssim \left(\frac{|B_R|}{|B_r|}\right)^\frac{1}{q} \left(\aver{B_R} |\Gamma f|^q d\mu\right)^\frac{1}{q}
$$
for every $f \in{\mathcal D}$ and $0<r<R$, if the space is doubling, with dimension $\nu$, the inequality $(D\textrm{G}_{q,\theta})$ holds for every $q>\nu$, with $\theta=\frac{\nu}{q}<1$.

\medskip

\subsubsection{Poincar\'e inequality}
\label{SubsubsectionPoincare}

Last, we shall assume that the carr\'e du champ operator $\Gamma$ satisfies the following scale-invariant Poincar\'e inequality
\begin{equation}\tag{$P_2$} 
\left( \aver{B} \left| f - \aver{B} f d\mu \right|^2 d\mu \right)^\frac{1}{2} \lesssim r \left(\int_B \Gamma(f)^2\, d\mu \right)^\frac{1}{2},  \label{P2}
\end{equation}
for every  $f\in {\mathcal D}_2(L)$ and every ball $B$ of radius $r$. We refer the reader to \cite{BCF1} for a precise study of the connection between Poincar\'e inequality, $L^p$-gradient estimates and de Giorgi property. Let us just point out that if $\Gamma$ satisfies the above Poincar\'e inequality and the gradient estimate $(\textrm{G}_{p_0})$, then there exists a parameter $\theta\in(0,1)$ such that the inequality $(D\textrm{G}_{p,\theta})$ holds for every $p\in [2,p_0)$. Note also that in the first and main part of this work, we assumed an upper Gaussian pointwise estimates for the gradient of the heat kernel equivalent to $(\textrm{G}_\infty)$. This assumption yields the Poincar\'e inequality $(P_2)$, the integrated gradient estimate $(\textrm{G}_{q_0})$ for every $q_0\in[2,\infty]$, and also de Giorgi property $(D\textrm{G}_{q,\theta})$ for every $\theta\in(0,1)$ and every $q\in[2,\infty)$. 

Our aim in this appendix will thus be to weaken the $(\textrm{G}_\infty)$ assumption made above into a combination of $(\textrm{G}_{q_0})$ and $(D\textrm{G}_{q,\theta})$, for some exponent $q_0,q$ and $\theta$. 

\bigskip

In the first and main part of this work, the paracontrolled calculus and its application to the 2-dimensional parabolic Anderson model equation was studied under the assumptions \eqref{UE} and $(\textrm{G}_\infty)$ that the heat kernel and its gradient satisfy pointwise Gaussian upper bounds. The aim of this section is to weaken the latter condition. Here are \textbf{examples} where the operator $\Gamma$ satisfies only the properties $(P_2), (\textrm{G}_q)$ and $(D\textrm{G}_q)$ for some $q>2$, and where $(\textrm{G}_\infty)$ does not hold.
 
 \medskip 
 
\begin{enumerate}
 \item \textbf{Conical manifolds.} Consider a compact Riemannian manifold $N$ of dimension $n-1\geq 1$, and define $M := (0,\infty) \times N$ as the conical manifold whose basis is $N$. It is known that $M$ is  a doubling manifold of dimension $n$ which satisfies \eqref{UE}. Moreover, as shown by Li in \cite{Li}, the operator $\Gamma$ satisfies $(\textrm{G}_q)$ if and only if
$$
q <q(N):=\left(\frac{1}{2} - \sqrt{\left(\frac{1}{2}-\frac{2}{n}\right)^2 + \frac{\lambda_1}{n}}\right)^{-1}
$$
where $\lambda_1$ is the first non-vanishing eigenvalue of the Laplace operator on $N$. As an example, if we consider $N= r{\mathbb S}^1$ the circle of radius $r>1$, then 
$$
q(N) = \frac{2r}{r-1}.
$$
So theorem \ref{thm:pam-bis} below allows us to solve the PAM equation on $M$ for $r$ sufficiently close to $1$. \vspace{0.2cm}
 
 \item \textbf{Elliptic perturbation of the Laplacien.} On the Euclidean space $\RR^d$, or any non-compact doubling Riemannian manifold satisfying Poincar\'e inequality $(P_2)$ and the Gaussian bound \eqref{UE}, we may consider a second order divergence form operator $\textrm{L} = -\textrm{div}(A\nabla)$ given by a map $A$ taking values in real symmetric matrices and satisfying the usual ellipticity condition. Then if $A$ is H\"older continuous, it is known that $-\textrm{L}$ generates a self-adjoint semigroup with \eqref{UE} and Gaussian pointwise bounds for the gradient of the semigroup $(\textrm{G}_\infty)$; see \cite{AMT}. In such a case we may apply the results proved in the first part. Following Auscher's work \cite{Aus}, we know that the combination of property \eqref{UE} with H\"older regularity of the heat kernel is stable under $L^\infty$ perturbation. So fix $A_0$ a H\"older continuous map with values in real symmetric matrices and satisfying usual ellipticity condition. Then for every $Q>2$, and any positive $\Theta$, there exists a positive constant $\epsilon$ such that for any map $A$ on the state space, with values in the space of real symmetric matrices, and such that $\|A-A_0\|_\infty \leq \epsilon$, the operator $L=-\textrm{div}(A\nabla)$ satisfies $(\textrm{G}_{q_0})$ for some $q_0>Q$, and has de Giorgi property $(D\textrm{G}_{q,\theta})$, for $\theta=d/q<\Theta$ and some $Q\leq q <q_0$. In such a situation, we may apply Theorem \ref{thm:pam-bis} and deduce that we can solve the PAM equation in such a $2$-dimensional context. \vspace{0.2cm}
 
 \item \textbf{Lipschitz domain with Neumann boundary conditions.} Similarly, consider an open and bounded subset $\Omega\subset \RR^2$ and consider for $\textrm{L}$ the self-adjoint Laplace operator associated with Neumann boundary conditions. Then by a change of variable, this situation is very similar to the previous one: if the boundary is sufficiently close (in a Lipschitz sense) to a smooth set (at least of regularity $C^{2}$), then we can solve the PAM equation.
\end{enumerate}

\bigskip

\subsection[\hspace{-1cm} {\sf Functional calculus and gradient estimates in H\"older and Sobolev spaces}]{Functional calculus and gradient estimates in H\"older and Sobolev spaces} 
\label{SubsectionAppFunctionalCalculus}

We start this section by quantifying the localization properties of the operators $P^{(a)}_t$ and $Q^{(a)}_t$, and their gradients, in $L^p$ spaces, before turning to the gradient estimates of the heat semigroup in the intrinsic H\"older and Sobolev spaces in section \ref{SubsubsectionAppGradientEstimates}.

\medskip

\subsubsection{Localization properties of the approximation operators $P^{(a)}_t$ and $Q^{(a)}_t$}
\label{SubsubsectionAppLocalizationProperties}

As we know, for every integer $a\geq 0$, the operators $P_t^{(a)}$ and $Q_t^{(a)}$ have a kernel satisfying Gaussian estimates \eqref{UE}. The above regularity assumptions $(\textrm{G}_{q}), (D\textrm{G}_{q}), (P_2)$ on the gradient operator actually imply much more.

\medskip

\begin{lemma} \label{lem:holderpt} 
{\sf Let $p_t$ stands for the kernel of $e^{-tL}$ or $P_t^{(a)}, Q_t^{(a)}$ for any integer $a\geq 1$. Under $(D\textrm{G}_{q,\theta})$ and $(\textrm{G}_{q_0})$ with Poincar\'e inequality $(P_2)$ for some $2\leq q <q_0$, we have the following H\"older regularity estimate for the heat kernel. For every $\eta\in(0,1-\theta]$, $t>0$ and almost every $x,y,z\in M$
$$ 
\left| p_t(x,z) - p_t(y,z) \right| \lesssim \left( \frac{d(x,y)}{\sqrt{t}}\right)^\eta V(z, \sqrt{t})^{-1} e^{-c \frac{d(x,z)^2}{t}}.
$$
} 
\end{lemma}

\medskip

We only sketch the proof and refer the reader to \cite{BCF1} for details.

\ssk

\begin{Dem}
We follow the argument of Morrey's inequality, which relies oscillation estimates to some gradient bounds. Let $x,y\in M$ be  Lebesgue points for $f=p_t(\cdot,z)$ with $d(x,y)\leq \sqrt{t}$, otherwise there is nothing to be done. Let $B_i(x) = B\big(x,2^{-i}d(x,y)\big)$, for $i\in\N$. Note that for all $i\in \N$,  $B_i(x)\subset B_0(x)$. By Poincar\'e's inequality, this yields
$$ \left|f(x) - \aver{B_0(x)} f d\mu \right| \lesssim \sum_{i\geq 0} 2^{-i}d(x,y)  \left(\aver{B_i(x)} |\Gamma f|^q \, d\mu \right)^{\frac{1}{q}}. $$
By considering $B_{\sqrt{t}}$ a ball of radius $\sqrt{t}$ containing both $x,y$,  $(D\textrm{G}_{q,\theta})$ yields
\begin{align*}
\left(\aver{B_i(x)} |\Gamma f|^q \, d\mu \right)^{\frac{1}{q}} \lesssim  \left(\frac{\sqrt{t}}{2^{-i} d(x,y)}\right)^ \theta \left[\left(\aver{B_{\sqrt{t}}} |\Gamma f|^q \, d\mu \right)^{\frac{1}{q}} + \sqrt{t} \|Lf\|_{L^\infty(B_{\sqrt{t}})} \right].
\end{align*}
Since $f=p_t(\cdot,z)$, by $(\textrm{G}_{q_0})$ and \eqref{UE} we know that 
$$
\left(\aver{B_{\sqrt{t}}} \big|\sqrt{t} \Gamma f\big|^q \, d\mu \right)^{\frac{1}{q}} +  \big\|tLf\big\|_{L^\infty(B_{\sqrt{t}})} \lesssim V\big(z, \sqrt{t}\big)^{-1} e^{-c \frac{d(x,z)^2}{t}}
$$
so we can conclude the proof by summing over $i$, since $\theta\in(0,1)$.
\end{Dem}

\medskip

Under the sole assumption (UE) that the kernels of the operators $P_t^{(a)}$ and $Q_t^{(a)}$ have Gaussian upper bounds, these operators are bounded in every $L^p$ space for $p\in[1,\infty]$, uniformly with respect to $t\in (0,1]$. Moreover, for every $p_1,p_2\in[1,\infty]$ and $t>0$, they satisfy the following $L^{p_1}$-$L^{p_2}$ off-diagonal estimates at scale $\sqrt{t}$, which quantify the localization properties of these operators. For every ball $B_1,B_2$ of radius $\sqrt{t}$, and for every function $f\in L^{p_1} (B_1)$, we have 
$$  
\left( \aver{B_2} |P_t^{(a)}f|^{p_2} \, d\mu \right)^\frac{1}{p_2} + \left( \aver{B_2} |Q_t^{(a)}f|^{p_2} \, d\mu \right)^\frac{1}{p_2} \lesssim e^{-c \frac{d(B_1,B_2)^2}{t}} \left(\aver{B_1} |f|^{p_1} \, d\mu \right)^\frac{1}{p_1}.
$$
One can refine this estimate by using off-diagonal estimates, such as done in \cite[Lemma 2.5, Lemma 2.6]{BCF2}.

\medskip

\begin{proposition}  \label{prop:gradientsemigroup-bis}   {\sf 
Assume $(\textrm{G}_{q_0})$ for some $q_0>2$.
\begin{itemize}
   \item[{\bf (i)}] For every non-negative \textit{integer} $a$ and every $p\in [2,q_0)$, the operators $\Gamma P_t^{(a)}$ and $\Gamma Q_t^{(a)}$ satisfy $L^2$-$L^p$ the following off-diagonal estimates at the scale $\sqrt{t}$. For every ball $B_1,B_2$ of radius $\sqrt{t}$ and every function $f\in L^2 (B_1)$, we have
$$ 
\left( \aver{B_2} |\sqrt{t} \Gamma Q_t^{(a)}f|^p \, d\mu \right)^\frac{1}{p} + \left( \aver{B_2} |\sqrt{t} \Gamma P_t^{(a)}f|^p \, d\mu \right)^\frac{1}{p} \lesssim e^{-c\frac{d(B_1,B_2)^2}{t}} \left(\aver{B_1} |f|^2 \, d\mu \right)^\frac{1}{2}.
$$
It follows in that we have
$$ 
\sup_{t>0} \, \left\{\Big\| \big(\sqrt{t} \Gamma\big)\big(P_t^{(a)}\cdot\big) \Big\|_{p \to p} + \Big\| \big(\sqrt{t} \Gamma\big)\big(Q_t^{(a)}\cdot\big) \Big\|_{p \to p} \right\} < \infty
$$ 
for every $p\in[2,q_0)$. \vspace{0.1cm}

   \item[{\bf (ii)}] For every positive real number $a$ and every positive $t$, the operator $Q_t^{(a)}$ is an integral operator with kernel $k_t^{(a)}$ that satisfies the inequality 
\begin{equation} \label{kernel-est}
	\abs{k_t^{(a)}(x,y)}\lesssim \frac{1}{V(x,\sqrt{t})^\lambda V(y,\sqrt{t})^{1-\lambda} } \left(1+\frac{d^2(x,y)}{t}\right)^{-a}
\end{equation}
for all $\lambda\in[0,1]$ and $\mu$-almost all $x,y \in M$. As a consequence the operator $Q_t^{(a)}$ satisfies the following $L^{p_1}$-$L^{p_2}$ off-diagonal bounds of order $a$ at scale $\sqrt{t}$, for every $p_1,p_2\in[1,{+\infty}]$. Given any balls $B_1,B_2$ of radius $\sqrt{t}$, and any function $f\in L^{p_1} (B_1)$, we have
$$ 
\left( \aver{B_2} \big|Q_t^{(a)}f\big|^{p_2} \, d\mu \right)^\frac{1}{p_2} \lesssim \left(1+\frac{d(B_1,B_2)^2}{t}\right)^{-a} \left(\aver{B_1} |f|^{p_1} \, d\mu \right)^\frac{1}{p_1}.
$$ 
\end{itemize}  }
\end{proposition}

\medskip 

Besides these localization property in the physical space, the approximation operators $Q_t^{(a)}$ satisfy some orthogonality properties, which will be of crucial use in proving the continuity properties of the paraproduct and resonant operators below, and which can be viewed as an analog of the Littlewood-Paley theory, as made clear in \cite[Proposition 2.13 and Lemma 2.15]{BCF2}.

\medskip

\begin{lemma}\label{lem:ortho} {\sf 
Let $a$ be a positive real number. Set
$$
\widetilde Q_t:= (t\mathcal{L})^\frac{a}{2} e^{-\frac{t}{2}\mathcal{L}}=2^\frac{a}{2} Q_\frac{t}{2}^\frac{a}{2},
$$ 
so that $Q_t^{(a)}=\widetilde Q_t^2$. Assume the Gaussian upper bound \eqref{UE} holds. Let also $F: (0,{+\infty})\times M\to\R$ be a measurable function and write $F_t(x)$ for $F(t,x)$. Then for every $p\in(1,{+\infty})$, one has
$$ 
\left\| \int_0^{+\infty} Q_t^{(a)}F_t \, \frac{dt}{t} \right\|_{p} \lesssim \left\| \left( \int_0^{+\infty} \big|\widetilde Q_t F_t\big|^2 \, \frac{dt}{t}\right)^\frac{1}{2} \right\|_{p},
$$
whenever the right hand side has a meaning and is finite. If $F=f$ does not depend on $t$, we have the following $L^p$-boundedness of the vertical square function
$$ 
\left\| \int_0^{+\infty} \big|Q_t^{(a)}f\big| \, \frac{dt}{t} \right\|_{p} \simeq \|f\|_{p}.
$$
}
\end{lemma}

\bigskip

\subsubsection{Gradient estimates in H\"older and Sobolev spaces}
\label{SubsubsectionAppGradientEstimates}

As said above, we shall now work in the following setting, strictly weaker than the geometrical setting used in the first five sections of this work.

\medskip

\noindent {\sf \textbf{Regularity assumptions} 
\begin{itemize}
   \item[\textbf{(i)}] The metric measure space $(M,d,\mu)$ is doubling and the semigroup satisfies the Gaussian bound \eqref{UE}. \vspace{0.1cm}

    \item[\textbf{(ii)}] The gradient operator $\Gamma$ satisfies $(\textrm{G}_{q_0})$ and $(D\textrm{G}_{q,\theta})$ for some $2\leq q <q_0 \leq \infty$, and the scale-invariant Poincar\'e inequality.   
\end{itemize}} 

\ssk

If $q_0=q=2$, we require that the $L^2$ Davies-Gaffney estimates hold instead of $(D\textrm{G}_{q,\theta})$. As we shall see below, one can extend the machinery of paracontrolled calculus to that setting in H\"older and Sobolev spaces. Recall the definition of the spaces $\Lambda^\sigma$ and $\calC^\sigma$ given in section \ref{SubsectionHolder}. The parameter $\theta$ is involved in the property $(D\textrm{G}_{q,\theta})$. The following embedding is proved as Proposition \ref{prop:caraholder} by using the fact proved in Lemma \ref{lem:holderpt} that the heat kernel is H\"older continuous, with exponent $1-\theta$, instead of its Lipschitz character. 
 
\medskip

\begin{proposition}  \label{prop:caraholder-bis}   {\sf 
For $\sigma\in(0,1)$, the space $\Lambda^\sigma$ is continuously embedded into $\calC^\sigma$. If $\sigma \in(0,1-\theta)$, the two spaces are the same with equivalent norms. }
\end{proposition}

\medskip

Sobolev spaces are naturally defined in terms of $L$ as follows.

\medskip

\begin{definition*}   {\sf
Fix an exponent $p\in (1,\infty)$, and $s\in{\mathbb R}$. A distribution $f\in \SSS$, is said to belong to the inhomogeneous \emph{\textbf{Sobolev space}} $W^{s,p}$ if
$$
\|f\|_{W^{s,p}} := \Big\|\big(1+L\big)^\frac{s}{2} f\Big\|_p \simeq \big\|e^{-L}f\big\|_p + \Big\|\big(1+L\big)^\frac{s}{2} f\Big\|_p <\infty.
$$   }
\end{definition*}

\medskip

\begin{proposition} \label{prop:gradient-Holder}  {\sf For $\sigma\in(-\infty,1-\theta)$ and $f \in \calC^\sigma$, we have
$$ 
\sup_{x\in M} \left( \aver{B(x,\sqrt{t})} \big|\sqrt{t} \Gamma e^{-tL} f\big|^q \, d\mu \right)^{\frac{1}{q}} \lesssim t^\frac{\sigma}{2} \|f\|_{\calC^\sigma}. 
$$
The same conclusion holds with any of the operators $P_t^{(a)}$, with an integer $a\geq 1$, in the role of $e^{-tL}$.   }
\end{proposition}

\medskip

\begin{Dem}
Consider $b\geq 1$, and write 
$$ 
\sqrt{t} \Gamma e^{-tL} f = \int_0^1 \sqrt{t} \Gamma Q_s^{(b)} e^{-tL} f \, \frac{ds}{s} + \sqrt{t}\Gamma P_1^{(b)} e^{-tL} f.
$$
For $s\leq t$, we have
\begin{align*}
\left( \aver{B(x,\sqrt{t})}  \Big| \sqrt{t} \Gamma Q_s^{(b)} e^{-tL} f\Big|^q \, d\mu \right)^{\frac{1}{q}} & = \left(\frac{s}{s+t}\right)^b \left( \aver{B(x,\sqrt{t})}  \Big| \sqrt{t} \Gamma Q_{s+t}^{(b)} f\Big|^q \, d\mu \right)^{\frac{1}{q}} \\
& \lesssim \left(\frac{s}{t}\right)^b  \sum_{\ell\geq 0} \gamma_{\ell}  \left( \aver{2^\ell B(x,\sqrt{t})}  \Big| Q_{\frac{s+t}{2}}^{(b)} f\Big|^q \, d\mu \right)^{\frac{1}{q}},
\end{align*}
where $\gamma_\ell$ are exponentially decreasing coefficients and where we used $L^q$-$L^q$ off-diagonal estimates of $\Gamma e^{-\frac{s+t}{2}L}$ (at the scale $\sqrt{s+t}\simeq \sqrt{t}$) with the relation
$$ 
Q_{s+t}^{(b)} = 2^{b} e^{-\frac{s+t}{2}L} Q_{\frac{s+t}{2}}^{(b)}. 
$$
So we have
\begin{align*}  
\left( \aver{B(x,\sqrt{t})}  \Big| \sqrt{t} L^{\frac{1}{2}} Q_s^{(b)} e^{-tL} f\Big|^q \, d\mu \right)^{\frac{1}{q}} & \lesssim \left(\frac{s}{t}\right)^b  \sum_{\ell\geq 0} \gamma_{\ell}  \big\| Q_{\frac{s+t}{2}}^{(b)} f \big\|_\infty \\
& \lesssim  \left(\frac{s}{t}\right)^b t^{\frac{\sigma}{2}} \|f\|_{\calC^\sigma},
\end{align*}
and we can integrate this inequality on the interval $s\in(0,t)$. For $s\geq t$, we use Property $(D\textrm{G}_{q,\theta})$ to have
\begin{align*}
\left( \aver{B(x,\sqrt{t})}  \Big| \sqrt{t} \Gamma Q_s^{(b)} e^{-tL} f\Big|^q \, d\mu \right)^{\frac{1}{q}} & \lesssim \left(\frac{t}{s}\right)^{\frac{1-\theta}{2}} \left( \aver{B(x,\sqrt{s})}  \Big| \sqrt{s} \Gamma Q_s^{(b)} e^{-tL} f\Big|^q \, d\mu \right)^{\frac{1}{q}} \\
& \qquad + \left(\frac{t}{s}\right)^{\frac{1-\theta}{2}} \left\| Q_s^{(b+1)} e^{-tL} f \right\|_{L^\infty(B(x,\sqrt{s}))} \\
& \lesssim \left(\frac{t}{s}\right)^{\frac{1-\theta}{2}} s^{\frac{\sigma}{2}} \|f\|_{\calC^\sigma},
\end{align*}
where we have used $ Q_s^{(b)}  = 2^b Q_{\frac{s}{2}}^{(b/2)} Q_{\frac{s}{2}}^{(b/2)}$ with $L^q$-$L^q$ (resp. $L^q$-$L^\infty$) off-diagonal estimates for $\Gamma Q_{s/2}^{(b/2)}$ (resp. $(sL)Q_{s/2}^{(b/2)}$), provided $b$ is large enough.
This inequality can be then integrated along $s\in(t,1)$ as soon as $\theta+\sigma<1$.

\ssk

We perform the same analysis for the term $\sqrt{t}\Gamma P_1^{(b)} e^{-tL} f$, which gives
\begin{align*}
\left( \aver{B(x,\sqrt{t})}  \big| \sqrt{t} \Gamma P_1^{(b)} e^{-tL} f\big|^q \, d\mu \right)^{\frac{1}{q}} & \lesssim t^\frac{1-\theta}{2} \left( \aver{B(x,1)}  \big| \Gamma P_1^{(b)}e^{-tL} f\big|^q \, d\mu \right)^{\frac{1}{q}} \\
& \qquad + t^\frac{1-\theta}{2} \left\| LP_1^{(b)} e^{-tL} f \right\|_{L^\infty(B(x,1))} \\
& \lesssim t^\frac{1-\theta}{2}  \|f\|_{\calC^\sigma}.
\end{align*}
The conclusion follows from this inequality since $t\in(0,1)$ and $\sigma<1-\theta$. 
\end{Dem}

\medskip

\begin{proposition} \label{prop} {\sf For $\alpha\in(0,1-\theta)$ and $0<2\delta<1-\theta-\alpha$, we have uniformly in $x\in M$ and $t>0$
$$
\left(\aver{B(x,\sqrt{t})} \big|\sqrt{t}\Gamma e^{-tL} f\big|^q \, d\mu \right)^{\frac{1}{q}} \lesssim \calM \Big[(tL)^{\frac{\alpha}{2}} Q_{\frac{t}{2}}^{(\delta)}f \Big],
$$
where $\calM$ is the Hardy-Littlewood maximal function. The same conclusion holds with any of the operators $P_t^{(a)}$, with an integer $a\geq 1$, in the role of $e^{-tL}$ and also by replacing $B(x,\sqrt{t})$ by any bigger ball $B\big(x,K\sqrt{t}\big)$ and the estimates are uniform with respect to $K\geq 1$.}
\end{proposition}

\medskip

\begin{Dem}
We write (for $a$ chosen large enough)
$$
\sqrt{t} \Gamma e^{-tL} =  \sqrt{t} \Gamma L^{-\frac{\alpha}{2}} e^{-tL} L^{\frac{\alpha}{2}} f = \sqrt{t} \int_0^\infty \Gamma Q_s^{(a)} e^{-tL} L^{\frac{\alpha}{2}} f \frac{ds}{s^{1-\frac{\alpha}{2}}}.
$$
For $s<t$, we then write
$$
Q_s^{(a)} e^{-tL} = \left(\frac{s}{t}\right)^{a} Q_{t}^{(a)} e^{-sL} =  \left(\frac{2s}{t}\right)^{a} Q_{\frac{t}{2}}^{(a-\delta)}e^{-sL}  Q_{\frac{t}{2}}^{(\delta)}  
$$
and using $L^1$-$L^q$ off-diagonal estimates of the carr\'e du champ of the semigroup, this yields
\begin{align*}
 \left(\aver{B(x,\sqrt{t})} \big|\sqrt{t}\Gamma Q_s^{(a)} e^{-tL} L^{\frac{\alpha}{2}} f \big|^q \, d\mu \right)^{\frac{1}{q}} & \lesssim \left(\frac{s}{t}\right)^{a} \calM\Big[L^{\frac{\alpha}{2}} Q_{\frac{t}{2}}^{(\delta)} f\Big](x). 
\end{align*}
For $t\leq s$ we have by $(D\textrm{G}_{q,\theta})$
\begin{align*}
 \left(\aver{B(x,\sqrt{t})} \big|\sqrt{t}\Gamma Q_s^{(a)} e^{-tL} L^{\frac{\alpha}{2}} f \big|^q \, d\mu \right)^{\frac{1}{q}} & \lesssim \left(\frac{s}{t}\right)^{\frac{\theta-1}{2}} \left(\aver{B(x,\sqrt{s})} \big|\sqrt{s}\Gamma Q_s^{(a)} e^{-tL} L^{\frac{\alpha}{2}} f \big|^q \, d\mu \right)^{\frac{1}{q}} \\
 & \qquad + \left(\frac{s}{t}\right)^{\frac{\theta-1}{2}} \left\| Q_s^{(a+1)} e^{-tL} L^{\frac{\alpha}{2}} f \right\|_{L^\infty(B(x,\sqrt{s}))} \\
 & \lesssim \left(\frac{s}{t}\right)^{\frac{\theta-1}{2}+\delta} \calM\Big[L^{\frac{\alpha}{2}} Q_t^{(\delta)}f\Big](x),
\end{align*}
where we used that $Q_s^{(a)} e^{-tL} = \big(\frac{s}{t}\big)^\delta Q_s^{(a-\delta)} Q_t^{(\delta)}$ with $L^1$-$L^q$ (resp. $L^1$-$L^\infty$) off-diagonal estimates for $\Gamma Q_s^{(a-\delta)}$ (resp. $Q_s^{(a+1-\delta)}$), provided $a$ is large enough. Hence, 
{\small   \begin{align*}
\left(\aver{B(x,\sqrt{t})} \Big|\sqrt{t}\Gamma e^{-tL} f\Big|^q \, d\mu \right)^{\frac{1}{q}} & \lesssim \left[ \int_0^t \left(\frac{s}{t}\right)^{a}  \frac{ds}{s^{1-\frac{\alpha}{2}}}  + \int_t^1 \left(\frac{s}{t}\right)^{\frac{\theta-1}{2}+\delta} \frac{ds}{s^{1-\frac{\alpha}{2}}} \right] \calM\Big[L^{\frac{\alpha}{2}} Q_t^{(\delta)} f\Big](x)  \\
& \lesssim \calM \Big[(tL)^{\frac{\alpha}{2}} Q_t^{(\delta)} f \Big](x),
\end{align*}   } 
due to $a$ large enough and $2\delta<1-\theta-\alpha$. We let the reader check the straightforward modifications that are required to deal with a bigger ball $B(x,K\sqrt{t})$, and that the estimates are uniform with respect to $K\geq 1$.
\end{Dem}

\medskip

Replacing the $L^1$-$L^q$ off-diagonal estimates by $L^p$-$L^q$ estimates, the same proof as above leads to the following result.

\medskip

\begin{proposition} \label{prop:gradient-Sobolev} {\sf 
Assume the local Ahlfors regularity.
Suppose $p\in[1,\infty)$ and $\alpha\in(-\infty,1-\theta+\frac{\nu}{p})$, and $f\in W^{\alpha,p}$. Then, uniformly in $x\in M$ and $t\in(0,1]$,
$$ 
\left( \aver{B(x,\sqrt{t})} \Big|\sqrt{t} \Gamma e^{-tL} f\Big|^q \, d\mu \right)^{\frac{1}{q}} \lesssim t^{-\frac{\nu}{2p}+ \frac{\alpha}{2}} \, \|f\|_{W^{\alpha,p}}.
$$
The same conclusion holds with any of the operators $P_t^{(a)}$, with an integer $a\geq 1$, in the role of $e^{-tL}$.   }
\end{proposition}

\medskip

\begin{Dem}
As previously, we write 
$$
\sqrt{t} \Gamma e^{-tL} = \sqrt{t} \int_0^1 \Gamma Q_s^{(a)} e^{-tL} L^{\frac{\alpha}{2}} f \frac{ds}{s^{1-\frac{\alpha}{2}}} + \sqrt{t}\Gamma P_1^{(a)} e^{-tL} f,
$$
for $a$ a large enough integer. For $s<t$, we then have
\begin{align*}
\left(\aver{B(x,\sqrt{t})} \Big|\sqrt{t}\Gamma Q_s^{(a)} e^{-tL} L^{\frac{\alpha}{2}} f \Big|^q \,d\mu \right)^{\frac{1}{q}} & \lesssim \left(\frac{s}{t}\right)^{a} t^{-\frac{\nu}{2p}}  \|f\|_{W^{\alpha,p}}
\end{align*}
and for $t\leq s$ we have by $(D\textrm{G}_{q,\theta})$
\begin{align*}
\left(\aver{B(x,\sqrt{t})} \Big|\sqrt{t}\Gamma Q_s^{(a)} e^{-tL} L^{\frac{\alpha}{2}} f \Big|^q \, d\mu \right)^{\frac{1}{q}} & \lesssim \left(\frac{s}{t}\right)^{\frac{\theta-1}{2}} \left(\aver{B(x,\sqrt{s})} |\sqrt{s}\Gamma Q_s^{(a)} e^{-tL} L^{\frac{\alpha}{2}} f |^q \, d\mu \right)^{\frac{1}{q}} \\
& \qquad + \left(\frac{s}{t}\right)^{\frac{\theta-1}{2}} \left\| Q_s^{(a+1)} e^{-tL} L^{\frac{\alpha}{2}} f \right\|_{L^\infty(B(x,\sqrt{s}))} \\
& \lesssim \left(\frac{s}{t}\right)^{\frac{\theta-1}{2}}  s^{-\nu/(2p)}  \|f\|_{W^{\alpha,p}}.\end{align*}
For the low frequency part, we have
\begin{align*}
\left(\aver{B(x,\sqrt{t})} \Big|\sqrt{t}\Gamma P_1^{(a)} e^{-tL} f \Big|^q \, d\mu \right)^{\frac{1}{q}} & \lesssim t^{\frac{1-\theta}{2}}  \|f\|_{W^{\alpha,p}}. 
\end{align*}
Hence, 
\begin{align*}
&  \left(\aver{B(x,\sqrt{t})} \Big|\sqrt{t}\Gamma e^{-tL} f\Big|^p \, d\mu \right)^\frac{1}{p}   \\
& \qquad \lesssim \left[ \int_0^t \left(\frac{s}{t}\right)^{a}  \frac{ds}{s^{1-\frac{\alpha}{2}}}  + \int_t^1 \left(\frac{s}{t}\right)^{\frac{\theta-1}{2} - \frac{\nu}{2p}} \frac{ds}{s^{1-\frac{\alpha}{2}}} + t^{\frac{1-\theta}{2}} \right] t^{-\frac{\nu}{2p}}  \|f\|_{W^{\alpha,p}}  \\
& \qquad \lesssim  t^{-\frac{\nu}{2p}+ \frac{\alpha}{2}}  \|f\|_{W^{\alpha,p}} ,
\end{align*}
due to $a$ large enough and $\alpha<1-\theta+\frac{\nu}{p}$. 
\end{Dem}

\medskip

\subsection[\hspace{-1cm} {\sf Paraproduct and commutator estimates in H\"older-Sobolev spaces}]{Paraproduct and commutator estimates in H\"older-Sobolev spaces} 
\label{SubsectionAppParaproduct}

This subsection is devoted to the statement and proofs of the main estimates about Paraproducts and commutators, in the current more general framework.

\medskip

\subsubsection{Paraproduct estimates}
\label{SubsubsectionAppParaproductEstimates}

We state in this paragraph the basic continuity estimates satisfied by the maps defined by the low frequency part, the paraproduct and the resonant terms -- see Subsection \ref{SubsectionParaproducts} for the precise definition of these quantities and for detailed proofs. The low-frequency part is easily bounded.

\medskip

\begin{proposition} \label{prop:low2-bis}  {\sf 
Fix an integer $b\geq 2$. For any $\alpha,\beta\in\R$ and every $\gamma>0$ we have for every $f \in \calC^\alpha$ and $g\in \calC^\beta$
\begin{equation} 
\big\| \Delta_{-1}(f,g)\big\|_{\calC^\gamma} \lesssim  \|f\|_{\calC^\alpha} \|g\|_{\calC^\beta}.  
\end{equation}   

If the space $(M,d,\mu)$ is locally Alhfors regular, then for every $\alpha,\beta,\gamma \in\R$ and $p\in[1,\infty)$, we have for every $f \in W^{\alpha,p}$ and $g\in W^{\beta,p}$
\begin{equation} 
\big\| \Delta_{-1}(f,g)\big\|_{W^{\gamma,p}} \lesssim  \|f\|_{W^{\alpha,p}} \|g\|_{W^{\beta,p}}.  
\end{equation}   }
\end{proposition}   

\medskip

The continuity properties of the paraproduct are given by the following statement.

\medskip

\begin{proposition}  \label{prop:para-Cs}   {\sf 
Fix an integer $b\geq 2$. For any $\alpha\in (-2,1)$ and $f\in \calC^\alpha$, we have
\begin{itemize}
 \item for every $g\in L^\infty$
\begin{equation} 
\Big\| \Pi^{(b)}_g(f)\Big\|_{\calC^\alpha} \lesssim \|g\|_\infty \|f\|_{\calC^\alpha} 
\label{eq:para1-Cs} \end{equation}
 \item for every $g\in\calC^\beta$ with $\beta<0$ and $\alpha+\beta\in (-2,1)$  
\begin{equation} 
\Big\| \Pi^{(b)}_g(f)\Big\|_{\calC^{\alpha+\beta}} \lesssim \|g\|_{\calC^\beta} \|f\|_{\calC^\alpha}. 
\label{eq:para2-Cs} \end{equation}
\end{itemize}  }
\end{proposition}  

\medskip

The proof is already given for Proposition \ref{prop:para} -- and only relies on \eqref{UE} (which is also assumed here). We then state the analog in Sobolev spaces.

\medskip

\begin{proposition}  \label{prop:para-Ws}   {\sf  Assume local Alhfors regularity.
Fix an integer $b\geq 2$ and $p\in[1,\infty)$. For any $\alpha\in (-2,1)$ and $f\in W^{\alpha,p}$, we have
\begin{itemize}
 \item for every $g\in W^{\beta,p}$ with $\frac{\nu}{p}<\beta<1$
\begin{equation} 
\Big\| \Pi^{(b)}_g(f)\Big\|_{W^{\alpha,p}} \lesssim \|g\|_{W^{\beta,p}} \|f\|_{W^{\alpha,p}} 
\label{eq:para1-W} \end{equation}
 \item for every $g\in W^{\beta,p}$ with $\beta<\frac{\nu}{p}$ and $\alpha+\beta-\frac{\nu}{p}\in (-2,1)$ 
\begin{equation} 
\Big\| \Pi^{(b)}_g(f)\Big\|_{W^{\alpha+\beta-\frac{\nu}{p},p}} \lesssim \|g\|_{W^{\beta,p}} \|f\|_{W^{\alpha,p}}. 
\label{eq:para2-W} \end{equation}
\end{itemize}  }
\end{proposition}  

\medskip

Even if the proof is not very difficult, we give the details here in order to explain how to use the $L^p$-orthogonality property put forward in Lemma \ref{lem:ortho}.

\medskip

\begin{Dem}
Recall that
$$  
\Pi^{(b)}_g(f) = \frac{1}{\gamma_b} \int_0^1  (tL)P_t^{(b)} \left(Q_t^{(b-1)} f\cdot P_t^{(b)} g\right) + Q_t^{(b-1)} \left((tL)P_t^{(b)} f\cdot P_t^{(b)} g\right) \, \frac{dt}{t}.
$$
With $s=\alpha+(\beta-\frac{\nu}{p})_{-}>-2$, Lemma \ref{lem:ortho} yields
\begin{align*}
 \|\Pi^{(b)}_g(f) \|_{W^{s,p}} & \lesssim \left \| \left( \int_0^1 t^{-s} \left|Q_t^{(b-1)} f\cdot P_t^{(b)} g\right|^2 + t^{-s} \left|(tL)P_t^{(b)} f\cdot P_t^{(b)} g \right|^2 \frac{dt}{t}\right)^\frac{1}{2} \right\|_p.
\end{align*}
If $\beta>\nu/p$ (and so $s=\alpha$) then uniformly with respect to $t>0$ we have due to the local Ahlfors regularity (which allows us to use a Sobolev embedding, see \cite[Lemma 10.5]{BCF2})
$$ 
\big\|P_t^{(b)} g\big\|_\infty \lesssim \|g\|_{\infty} \lesssim \|g\|_{W^{\beta,p}}
$$
and so
\begin{align*}
\|\Pi^{(b)}_g(f) \|_{W^{s,p}} & \lesssim \left \| \left( \int_0^1 t^{-s} \left|Q_t^{(b-1)} f\right|^2 + t^{-s} \left|(tL)P_t^{(b)} f\right|^2 \frac{dt}{t}\right)^\frac{1}{2} \right\|_p \|g\|_{W^{\beta,p}} \\
& \lesssim \|f\|_{W^{s,p}}\|g\|_{W^{\beta,p}},
\end{align*}
where we used again Lemma \ref{lem:ortho}.
If $\beta < \frac{\nu}{p}$ (and so $s=\alpha+\beta-\frac{\nu}{p}$), then 
\begin{align*}
 \big\|P_t^{(b)} g\big\|_\infty & \lesssim \int_t^1 \big\|Q_s^{(b)} g\big\|_\infty \frac{ds}{s} + \big\|P_1^{(b)}g\big\|_\infty \\
 & \lesssim \int_t^1 s^{\frac{\beta}{2}} \big\|Q_s^{(b-\frac{\beta}{2})} L^{\frac{\beta}{2}} g \big\|_\infty \frac{ds}{s} +  \|g\|_{W^{\beta,p}} \\
 & \lesssim \left(1+\int_t^1 s^{\frac{1}{2}(\beta-\frac{\nu}{p})} \frac{ds}{s} \right) \|g\|_{W^{\beta,p}} \\
 & \lesssim t^{\frac{1}{2}(\beta-\frac{\nu}{p})}\|g\|_{W^{\beta,p}}.
\end{align*}
Hence, we conclude with Lemma \ref{lem:ortho} since
\begin{align*}
\big\|\Pi^{(b)}_g(f) \big\|_{W^{s,p}} & \lesssim \left \| \left( \int_0^1 t^{-\alpha} \left|Q_t^{(b-1)} f\right|^2 + t^{-\alpha} \left|(tL)P_t^{(b)} f\right|^2 \frac{dt}{t}\right)^\frac{1}{2} \right\|_p \|g\|_{W^{\beta,p}} \\
 & \lesssim \|f\|_{W^{\alpha,p}}\|g\|_{W^{\beta,p}}.
\end{align*}
\end{Dem}

\medskip

\begin{proposition} \label{prop:diag-Cs}   {\sf 
Fix an integer $b>2$. For any $\alpha,\beta\in (-\infty,1-\theta)$ with $\alpha+\beta>0$, for every $f\in\calC^\alpha$ and $g\in \calC^\beta$, we have the continuity estimate
$$ 
\Big\| \Pi^{(b)}(f,g) \Big\|_{\calC^{\alpha+\beta}} \lesssim \|f\|_{\calC^\alpha} \|g\|_{\calC^\beta}.
$$   }
\end{proposition}   

\medskip

\begin{Dem}
We only study the most difficult term in the resonant term $\Pi^{(b)}(f,g)$, which takes the form
\begin{equation} \label{def:A}
A(f,g):=\int_0^1 P_t^{(b-1)} \Gamma\left(\sqrt{t} P_t^{(b-1)} f, \sqrt{t}  P_t^{(b-1)} g\right) \, \frac{dt}{t}.
\end{equation}
$P_t^{(b-1)}$ satisfies $L^1$-$L^\infty$ off-diagonal estimates at order $N$ ($N$ can be chosen arbitrarily large, since $b$ is an integer) and so
$$ 
\left|P_t^{(b-1)} (h) (x)\right| \lesssim \sum_{\ell \geq 0} 2^{-\ell N} \left( \aver{2^\ell B(x,\sqrt{t})} |h| \, d\mu \right).
$$
With $h=\sqrt{t} \Gamma P_t^{(b-1)}f \cdot \sqrt{t} \Gamma P_t^{(b-1)} g$ and H\"older's inequality, we deduce that
{\small
$$ 
A(f,g)(x)  \lesssim \sum_{\ell \geq 0} 2^{-\ell N} \int_0^1 \left( \aver{B(x,2^\ell \sqrt{t})} |\sqrt{t} \Gamma P_t^{(b-1)}f|^2 \, d\mu \right)^\frac{1}{2} \left( \aver{B(x,2^\ell \sqrt{t})} |\sqrt{t} \Gamma P_t^{(b-1)}g|^2 \, d\mu \right)^\frac{1}{2} \frac{dt}{t}.
$$   }
We then conclude as previously, with Proposition \ref{prop:gradient-Holder}.
\end{Dem}

\medskip

We then give the analog estimate in Sobolev spaces.

\medskip

\begin{proposition} \label{prop:diag-Ws}   {\sf 
Assume the local Ahlfors regularity. Fix an integer $b>2$ and $p\in (1,\infty)$. For any $\alpha,\beta\in (-\infty,1-\theta)$ with $\alpha+\beta>\frac{\nu}{p}$, for every $f\in W^{\alpha,p}$ and $g\in W^{\beta,p}$, we have the continuity estimate
$$ 
\Big\| \Pi^{(b)}(f,g) \Big\|_{W^{\alpha+\beta-\frac{\nu}{p},p}} \lesssim \|f\|_{W^{\alpha,p}} \|g\|_{W^{\beta,p}}.
$$   }
\end{proposition}   

\medskip

\begin{Dem}
Again, we only study the most difficult term $A(f,g)$ defined in \eqref{def:A}. With $s:=\alpha+\beta-\frac{\nu}{p}>0$, we have by Lemma \ref{lem:ortho}
$$ 
\big\|L^\frac{s}{2} A(f,g)\big\|_{p}  \lesssim \left\| \left(\int_0^1  t^{-s} \left|(tL)^\frac{s}{2}P_t^{(b-1)} \Gamma\left(\sqrt{t} P_t^{(b-1)} f, \sqrt{t}  P_t^{(b-1)} g\right) \right|^2 \frac{dt}{t}\right)^\frac{1}{2} \right\|_p.
$$
Since $s>0$, $(tL)^\frac{s}{2} P_t^{(b-1)}$ satisfies $L^1$-$L^\infty$ off-diagonal estimates at order $\frac{s}{2}$ (see Proposition \ref{prop:gradientsemigroup-bis}) and so
$$ 
\left|(tL)^\frac{s}{2} P_t^{(b-1)} (h) (x)\right| \lesssim \sum_{\ell \geq 0} 2^{-\ell \frac{s}{2}} \left( \aver{2^\ell B(x,\sqrt{t})} |h| \, d\mu \right).
$$
With 
$$
h=\sqrt{t} \Gamma P_t^{(b-1)}f \cdot \sqrt{t} P_t^{(b-1)} g
$$ 
and H\"older's inequality, we deduce that $\left\|L^\frac{s}{2}A(f,g) \right\|_p$ is bounded by
{\small   
\begin{equation} 
\label{eq:hold}
\sum_{\ell \geq 0} 2^{-\ell s} \left\| \left( \int_0^1 t^{-s} \left( \aver{B(x,2^\ell \sqrt{t})} |\sqrt{t} \Gamma P_t^{(b-1)}f|^2 \, d\mu \right) \left( \aver{B(x,2^\ell \sqrt{t})} |\sqrt{t} \Gamma P_t^{(b-1)}g|^2 \, d\mu \right) \frac{dt}{t} \right)^\frac{1}{2} \right\|_p .  
\end{equation}   } 
Then using Proposition \ref{prop:gradient-Sobolev} with the Ahlfors regularity, we have
$$ 
\left( \aver{B(x,2^{\ell} \sqrt{t})} \big|\sqrt{t} \Gamma P_t^{(b-1)} g\big|^2 \, d\mu \right)^\frac{1}{2} \lesssim t^{-\frac{\nu}{2p}} t^{\frac{\beta}{2}} \|g\|_{W^{\beta,p}}.
$$
By Proposition \ref{prop}, we get
$$  
\left(\aver{B(x,2^\ell \sqrt{t})} \big|\sqrt{t}\Gamma P_t^{(b-1)} f\big|^2 \, \right)^\frac{1}{2} \lesssim \calM \Big[(tL)^{\frac{\alpha}{2}} Q_t^{(\delta)}f \Big](x),
$$
for some $\delta>0$. Hence
\begin{align*}
 \left\| L^\frac{s}{2} A(f,g)\right\|_p  & \lesssim \|g\|_{W^{\beta,p}} \sum_{\ell \geq 0} 2^{-\ell s} \left\|\left( \int_0^1 t^{-s} \left|\calM\Big[(tL)^{\frac{\alpha}{2}} Q_t^{(\delta)}f \Big]\right|^2 t^{\beta - \frac{\nu}{p}} \, \frac{dt}{t} \right)^\frac{1}{2} \right\|_p \\
 & \lesssim \|g\|_{W^{\beta,p}} \left\| \left( \int_0^1  \left|\calM\Big[ Q_t^{(\delta)}L^{\frac{\alpha}{2}}f\Big]\right|^2 \frac{dt}{t} \right)^\frac{1}{2} \right\|_p.
 \end{align*}
Using the Fefferman-Stein inequality on the maximal function together with the $L^p$-boundedness of the square function -- see Lemma \ref{lem:ortho}, we deduce that
$$ 
\big\| L^\frac{s}{2} A(f,g)\big\|_{p} \lesssim \big\| L^{\frac{\alpha}{2}}f \big\|_p \|g\|_{W^{\beta,p}} \lesssim \|f\|_{W^{\alpha,p}}\|g\|_{W^{\beta,p}}.
$$
By a similar reasoning, we have
\begin{align*}
 \big\|A(f,g)\big\|_{p}  & \lesssim \int_0^1 \left\| P_t^{(b-1)} \Gamma\left(\sqrt{t} P_t^{(b-1)} f, \sqrt{t}  P_t^{(b-1)} g\right) \right\|_p \, \frac{dt}{t} \\
  & \lesssim \|g\|_{W^{\beta,p}} \int_0^1 \left\|\calM\Big[ Q_t^{(\delta)}L^{\frac{\alpha}{2}}f\Big] \right\|_p t^{-\frac{\nu}{2p}} t^{\frac{\beta}{2}} \, \frac{dt}{t} \\
  & \lesssim \|f\|_{W^{\alpha,p}}\|g\|_{W^{\beta,p}} \left(\int_0^1 t^\frac{s}{2} \frac{dt}{t}\right) \\
  & \lesssim \|f\|_{W^{\alpha,p}}\|g\|_{W^{\beta,p}},
\end{align*}
where we used that $s>0$ and the $L^p$-boundedness of the approximation operators. That concludes the proof of the estimate
$$  
\big\|A(f,g)\big\|_{W^{s,p}} \lesssim \|f\|_{W^{\alpha,p}}\|g\|_{W^{\beta,p}}.
$$
Since the resonant part $\Pi^{(b)}$ can be split into a finite number of terms similar to $A(f,g)$, we then deduce the Sobolev boundedness of the resonent part.
\end{Dem}

\bigskip

\subsubsection{Commutator estimates}

We now focus on the adaptation of the commutator estimates given above in Proposition \ref{prop:commutator}.

\medskip

\begin{proposition}  \label{prop:commutator-Cs}   {\sf 
Consider the a priori unbounded trilinear operator
$$ 
C(f,g,h):=\Pi^{(b)}\Big(\Pi^{(b)}_g(f),h\Big) - g \Pi^{(b)}(f,h), 
$$
on $\SSS$. Let $\alpha,\beta,\gamma$ be H\"older regularity exponents with $\alpha\in(-1,1-\theta), \beta\in(0,1-\theta)$ and $\gamma\in (-\infty,1]$. If
$$
0 < \alpha + \beta + \gamma  \qquad \textrm{ and } \qquad \alpha + \gamma < 0
$$
then, setting $\delta := (\alpha+\beta)\wedge (1-\theta) +\gamma$, we have 
\begin{equation}
\label{EqCommutator-bis}
\big\| C(f,g,h) \big\|_{\calC^{\delta}} \lesssim \|f\|_{\calC^\alpha} \, \|g\|_{\calC^\beta} \, \|h\|_{\calC^\gamma},
\end{equation}
for every $f\in \calC^\alpha$ ,$g\in\calC^\beta$ and $h\in \calC^\gamma$; so the commutator defines a trilinear map from $\calC^\alpha\times\calC^\beta\times\calC^\gamma$ to $\calC^\delta$.  }
\end{proposition}

\medskip

\begin{Dem}
We refer to the proof of Proposition \ref{prop:commutator} for details and we keep the same notations. So it suffices to focus on a generic term of the form
$$ 
D(f,g,h) := \mcR\big(\mcA(f,g),h\big) - g\,\mcR(f,h)
$$
and prove the continuity estimate \eqref{EqCommutator-bis} for it. As previously done, we split the proof of the commutator estimate \eqref{EqCommutator-bis} for $D$ in two steps, and introduce an intermediate quantity 
$$
\mcS(f,g,h) := \int_0^1 \mcP^1_t\left( \Gamma\big(\sqrt{t} \mcP^2_t f , \sqrt{t} \mcP_t^3 h\big) \cdot P_t g\right) \, \frac{dt}{t},
$$
for which we shall prove that we have both 
\begin{equation}
\label{EqEstimateD1-bis}
\big\| g\,\mcR(f,h) - \mcS(f,g,h)\big\|_{\calC^\delta} \lesssim \|f\|_{\calC^\alpha} \, \|g\|_{\calC^\beta} \, \|h\|_{\calC^\gamma}
\end{equation}
and 
\begin{equation}
\label{EqEstimateD2-bis}
\big\| D(f,g,h) - \mcS(f,g,h)\big\|_{\calC^\delta} \lesssim \|f\|_{\calC^\alpha} \, \|g\|_{\calC^\beta} \, \|h\|_{\calC^\gamma}.
\end{equation}

\medskip

\noindent {\bf Step 1 -- proof of \eqref{EqEstimateD1-bis}.} 
This part is very similar to Step 1 of Proposition \ref{prop:commutator}, so we only point out the modifications.
Using Gaussian pointwise estimates for the kernel of $\mcP^1_t$, we have for almost every $x\in M$
\begin{align*}
 & \Big| \mcP^1_t \Big(  \Gamma\big(\sqrt{t} \mcP^2_t f , \sqrt{t} \mcP_t^3 h\big) \cdot \big(g(x) - P_t g\big)\Big)(x) \Big| \lesssim \sum_{\ell \geq 0} e^{-c 4^\ell} \left( \sup_{d(x,y) \leq 2^\ell \sqrt{t}} \ \big|g(x) - P_t g(y)\big|\right)... \\
 & \qquad ...\left( \aver{2^\ell B(x,\sqrt{t})} \big|\sqrt{t} \Gamma\big(\mcP^2_t f\big)\big|^2 \, d\mu \right)^\frac{1}{2} \left( \aver{2^\ell B(x,\sqrt{t})} \big|\sqrt{t} \Gamma\big(\mcP^3_t h\big)\big|^2 \, d\mu \right)^\frac{1}{2}. 
\end{align*}
By using the $\calC^\beta$-regularity of $g$ as well as Proposition \ref{prop:gradient-Holder} to estimate the $L^2$ averages, we get
\begin{align*}
 \Big| \mcP^1_t \Big(  \Gamma\big(\sqrt{t} \mcP^2_t f , \sqrt{t} \mcP_t^3 h\big) \cdot \big(g(x) - P_t g\big)\Big)(x) \Big|  & \lesssim \left(\sum_{\ell \geq 0} e^{-c 4^\ell} (4^\ell t)^{\frac{\beta}{2}} t^{\frac{\alpha}{2}} t^{\frac{\gamma}{2}}\right) \|f\|_{\calC^\alpha}\|g\|_{\calC^\beta}\|h\|_{\calC^\gamma} \\
&   \lesssim t^\frac{\alpha+\beta+\gamma}{2} \|f\|_{\calC^\alpha}\|g\|_{\calC^\beta}\|h\|_{\calC^\gamma}.
\end{align*}
Consequently, the continuity estimate \eqref{EqEstimateD1-bis} in $L^\infty$ comes from integrating with respect to time, taking into account the fact that $\alpha+\beta+\gamma > 0$.

\ssk

Then to estimate the regularity of  $g \mcR(f,h) - \mcS(f,g,h)$, one can exactly reproduce the same reasoning as for Proposition \ref{prop:commutator} by using the H\"older regularity of the heat kernel (Lemma \ref{lem:holderpt}), which involves the condition $ \alpha+\beta+\gamma<1-\theta$ (since $\beta<1-\theta$ and $\alpha+\gamma<0$).

\medskip

\noindent {\bf Step 2 -- proof of \eqref{EqEstimateD2-bis}}. Given the collection $\big(Q_r := Q_r^{(1)}\big)_{r\in(0,1]}$ of operators, we need to prove that we have
\begin{equation}
\label{EqSecondBitCommutator-bis}
\Big\| Q_r\Big(\mcR\big(\mcA(f,g),h\big) - \mcS(f,g,h)\Big)\Big\|_\infty \lesssim r^\frac{\delta}{2}.
\end{equation}
for every $r\in(0,1]$, and where
\begin{equation}
\label{eq:eq4-bis}
\begin{split}
&\mcR\big(\mcA(f,g),h\big) - \mcS(f,g,h) \\ 
&= \int_0^1 \mcP^1_t \Gamma \left( \sqrt{t}\left\{\int_0^1 \mcP^2_t  \mcQ^1_s\left(\mcQ^2_s f \cdot \mcP_s^3 g\right) \, \frac{ds}{s} - P_t g \cdot\mcP^2_t f\right\} \,,\,\sqrt{t} \mcP_t^3 h\right) \, \frac{dt}{t}. 
\end{split}
\end{equation}
We are going to follow the same argument as for Proposition \ref{prop:commutator} and we only detail the modifications. So we set
$$ 
A_{t}(f,g) := \sqrt{t} \Gamma\left(\int_0^1 \mcP^2_t  \mcQ^1_s\left(\mcQ^2_s f \cdot \mcP_s^3 g \right) \, \frac{ds}{s} - P_t g \mcP^2_t f\right).
$$
and using the $L^1$-$L^\infty$ off-diagonal estimates of $\mcP_t^1$, we then deduce that for almost every $x_0\in M$
\begin{align} 
\label{eq:spl1-bis}
& \Big| \mcR\big(\mcA(f,g),h\big) - \mcS(f,g,h)\Big)(x_0)\Big| \nonumber \\
& \lesssim \sum_{\ell \geq 0} \int_0^1 e^{-c4^\ell} \left( \aver{B(x_0, 2^\ell \sqrt{t})} |A_t(f,g)|^2 \, d\mu \right)^\frac{1}{2} \left(\aver{B(x_0, 2^\ell \sqrt{t})} |\sqrt{t} \Gamma \mcP_t^4 h|^2 \, d\mu \right)^\frac{1}{2}.  
\end{align}
Using a suitable normalization of the operators $\int_0^1 \mcQ^1_s \mcQ^2_sf  \, \frac{ds}{s}=f-\mcP_1f$ for some operator $\mcP_1$, it yields for every $x\in M$
\begin{align*}
A_t (f,g)(x) \leq \int_0^1  \sqrt{t} \Gamma \mcP^2_t  \mcQ^1_s\Big(\mcQ^2_s f  \big(\mcP_s^3 g - P_t g(x)\big)\Big)(x) \, \frac{ds}{s} + \big|P_t g(x)\big| \sqrt{t}\Gamma\big(\mcP^2_t\mcP_1f\big)(x).
\end{align*}
This quantity will then be integrated on $B_\ell:=B(x_0,2^\ell \sqrt{t})$, so we first aim to replace $P_tg(x)$ by $\aver{B_\ell} P_t g\, d\mu$.
Observe that
\begin{align} 
\label{eq:spl-bis} 
A_t(f,g)(x)  & \leq \int_0^1  \sqrt{t} \Gamma \mcP^2_t  \mcQ^1_s\left(\mcQ^2_s f  \big(\mcP_s^3 g - \aver{B_\ell} P_t g\, d\mu \big)\right)(x) \, \frac{ds}{s} \nonumber \\
 & + \left| P_tg(x) - \aver{B_\ell} P_tg \, d\mu \right|\sqrt{t} \Gamma \mcP^2_t f(x) + |P_t g(x)| \sqrt{t}\Gamma[\mcP^2_t\mcP_1f](x). 
\end{align}
As before, we use $\beta>0$ and the $\calC^\beta$ regularity of $g$ to have
\begin{align*} \left| P_tg(x) - \aver{B_\ell} P_tg \, d\mu\right| \lesssim \left(2^\ell \sqrt{t}\right)^\beta  \|g\|_{\calC^\beta}, 
\end{align*}
and uniformly in $y\in M$
\begin{align*}
\Big| \mcP_s^3 g(y) - \aver{B(x_0,\sqrt{t})} P_tg \, d\mu \Big|  & \lesssim \left( \max(s,t)^{\frac{\beta}{2}} + d(x_0,y)^\beta\right)  \|g\|_{\calC^\beta}.
\end{align*}
Moreover, it follows from the composition of $L^2$ off-diagonal estimates (corresponding to a $L^2$ analog of Lemma \ref{lem:gauss} -- Part1, see also \cite[Lemma 2.5]{BCF2}), that the operator $\sqrt{t}\Gamma\Big(\mcP^2_t  Q^1_s \Big)$ satisfies $L^2$ off-diagonal estimates at the scale $\max(s,t)$ with an extra factor $\left(\frac{\min(s,t)}{\max(s,t)}\right)$; so if one sets $\tau:=\max(s,t)$, we have with \eqref{eq:spl-bis}
{\small \begin{align*}
& \left( \aver{B(x_0,2^\ell \sqrt{t})} |A_t(f,g)|^2 \, d\mu \right)^\frac{1}{2} \\
& \qquad\lesssim \int_0^1 \sum_{k\geq 0} \left(\frac{\min(s,t)}{\max(s,t)}\right) 
e^{-c 4^k} (4^{k+\ell} \tau)^{\frac{\beta}{2}} \left( \aver{B(x_0, 2^{k+\ell} \sqrt{\tau})} |\mcQ^2_s f|^2 \, d\mu \right)^\frac{1}{2}  \|g\|_{\calC^\beta} \, \frac{ds}{s} \\
& \qquad \qquad + \left(\aver{B(x_0,2^\ell \sqrt{t})} |\sqrt{t} \Gamma \mcP^2_t f|^2 \, d\mu \right)^\frac{1}{2} \left(2^\ell \sqrt{t}\right)^\beta  \|g\|_{\calC^\beta} + t^{\frac{1-\theta}{2}}\|g\|_\infty\|f\|_{\calC^\alpha} \\
& \qquad \lesssim 2^{\ell\beta} \left[ \int_0^1 \left(\frac{\min(s,t)}{\max(s,t)}\right) s^{\frac{\alpha}{2}} \tau^{\frac{\beta}{2}} \, \frac{ds}{s} + t^{\frac{\sigma}{2}} \right] \|f \|_{\calC^\alpha} \|g\|_{\calC^\beta}  \\
& \qquad \lesssim 2^{\ell \beta}  t^{\frac{\alpha+\beta}{2}} \|f \|_{\calC^\alpha} \|g\|_{\calC^\beta},
\end{align*}}
where we used Proposition \ref{prop:gradient-Holder} and the fact that $\alpha>-1$ and $\alpha+\beta\leq 1-\theta$ to estimate the integral over $s$.
Observe that in the case where $\alpha+\beta\geq 1-\theta$, we get
\begin{align*}
\left( \aver{B(x_0,2^\ell \sqrt{t})} |A_t(f,g)|^2 \, d\mu \right)^\frac{1}{2} \lesssim 2^{\ell \beta} t^{\frac{1-\theta}{2}}\|f \|_{\calC^\alpha} \|g\|_{\calC^\beta}.
\end{align*}
Coming back to the identity \eqref{eq:spl1-bis}, with Proposition \ref{prop:gradient-Holder} we have
\begin{align*} 
\Big| \mcR\big(\mcA(f,g),h\big) - \mcS(f,g,h)\Big)(x_0)\Big| & \lesssim \left(\int_0^1 t^\frac{\alpha+\beta+\gamma}{2} \, \frac{dt}{t} \right) \|f \|_{\calC^\alpha} \|g\|_{\calC^\beta} \|h\|_{\calC^\gamma} \\
& \lesssim \|f \|_{\calC^\alpha} \|g\|_{\calC^\beta} \|h\|_{\calC^\gamma},
\end{align*}
since $\alpha+\beta+\gamma>0$, uniformly for every $x_0\in M$. We then conclude to 
$$ 
\Big\| \mcR\big(\mcA(f,g),h\big) - \mcS(f,g,h)\Big)\Big\|_\infty \lesssim \|f \|_{\calC^\alpha} \|g\|_{\calC^\beta} \|h\|_{\calC^\gamma}.
$$
Moreover, taking into account that we have $Q_r^{(1)} \mcP_t^1 = \frac{r}{t} \mcP_r Q^{(1)}_t$  for $t\geq r$, we see that the estimate \eqref{EqSecondBitCommutator-bis} holds true (see thee proof of Proposition \ref{prop:commutator}). 
\end{Dem}

\medskip

We then aim to have a similar commutator estimate in Sobolev spaces.

\medskip

\begin{proposition}  \label{prop:commutator-Ws}   {\sf Assume the local Ahlfors regularity. Let $\alpha,\beta,\gamma$ be regularity exponents and $p\in(1,\infty)$ with $\alpha\in(-1,1-\theta), \beta\in\big(\frac{\nu}{p} ,1-\theta\big)$ and $\gamma\in (-\infty,1]$. If
$$
\frac{2\nu}{p} < \alpha + \beta + \gamma  \qquad \textrm{ and } \qquad \alpha + \gamma < \frac{\nu}{p}
$$
then, setting $\delta := (\alpha+\beta-\frac{\nu}{p})\wedge 1 +\gamma - \frac{\nu}{p} > 0$ and assume that $2\delta > \beta - \frac{\nu}{p}$. We have 
\begin{equation}
\label{EqCommutator-Ws}
\big\| C(f,g,h) \big\|_{W^{\delta,p}} \lesssim \|f\|_{W^{\alpha,p}} \, \|g\|_{W^{\beta,p}} \, \|h\|_{W^{\gamma,p}},
\end{equation}
for every $f\in W^{\alpha,p}$ ,$g\in W^{\beta,p}$ and $h\in W^{\gamma,p}$; so the commutator defines a trilinear map from $W^{\alpha,p}\times W^{\beta,p} \times W^{\gamma,p}$ to $W^{\delta,p}$.  }
\end{proposition}

\medskip

We follow the exact same proof as previously, so we keep the same notations and only focus on the modifications.

\medskip

\begin{Dem}
Consider a generic term of the form
$$ 
D(f,g,h) := \mcR\big(\mcA(f,g),h\big) - g\,\mcR(f,h)
$$
and prove the continuity estimate \eqref{EqCommutator-Ws} for it. Aiming that, we split into two terms by introducing the quantity 
$$
\mcS(f,g,h) := \int_0^1 \mcP^1_t\left( \Gamma\big(\sqrt{t} \mcP^2_t f , \sqrt{t} \mcP_t^3 h\big) \cdot P_t g\right) \, \frac{dt}{t}.
$$
for which we shall prove that we have both 
\begin{equation}
\label{EqEstimateD1-Ws}
\big\| g\,\mcR(f,h) - \mcS(f,g,h)\big\|_{W^{\delta,p}} \lesssim \|f\|_{W^{\alpha,p}} \, \|g\|_{W^{\beta,p}} \, \|h\|_{W^{\gamma,p}}
\end{equation}
and 
\begin{equation}
\label{EqEstimateD2-Ws}
\big\| D(f,g,h) - \mcS(f,g,h)\big\|_{W^{\delta,p}} \lesssim \|f\|_{W^{\alpha,p}} \, \|g\|_{W^{\beta,p}} \, \|h\|_{W^{\gamma,p}}.
\end{equation}

\medskip

\noindent {\bf Step 1 -- proof of \eqref{EqEstimateD1-Ws}.} We first prove a weaker version of the continuity estimate \eqref{EqEstimateD1-Ws}, under the form of the inequality
\begin{equation} 
\label{eq:eq1-Ws}
\big\| g\,\mcR(f,h) - \mcS(f,g,h)\big\|_p \lesssim \|f\|_{W^{\alpha,p}} \, \|g\|_{W^{\beta,p}} \, \|h\|_{W^{\gamma,p}}.
\end{equation}
As previously, we have 
\begin{equation} 
\label{eq:eq2-Ws} 
\big(g \mcR(f,h) - \mcS(f,g,h)\big)(x) = \int_0^1 \mcP^1_t\Big( \Gamma\big(\sqrt{t} \mcP^2_t f , \sqrt{t} \mcP_t^3 h\big) \cdot \big(g(x) - P_t g\big)\Big)(x) \, \frac{dt}{t},
\end{equation}
for $\mu$-almost every $x\in M$. Since $g\in \calC^\beta$, with $\beta > \frac{\nu}{p}$ then $g\in W^{\beta,p} \subset \calC^{\beta-\frac{\nu}{p}}$ and so
\begin{align*}
\| P_t g -g \|_\infty & \lesssim t^{\frac{1}{2}\big(\beta-\frac{\nu}{p}\big)}\,\|g\|_{W^{\beta,p}}.
\end{align*} 
Hence
\begin{align}
\label{eq:ppp-Ws}
\big|P_t g(y) - g(x)\big| & \lesssim \left(\sqrt{t}+ d(x,y)\right)^{\beta-\frac{\nu}{p}}  \|g\|_{W^{\beta,p}}, 
\end{align}
for every $x,y\in M$. Coming back to equation \eqref{eq:eq2-Ws} and using Gaussian pointwise estimates for the kernel of $\mcP^1_t$, we have for almost every $x\in M$
\begin{align*}
\Big| \mcP^1_t \Big(  \Gamma\big(\sqrt{t} \mcP^2_t f , \sqrt{t} \mcP_t^3 h\big) \cdot \big(g(x) - P_t g\big)\Big)(x) \Big|
\end{align*}
is bounded above by
\begin{align*}
t^{\frac{1}{2}\big(\beta-\frac{\nu}{p}\big)}\|g\|_{W^{\beta,p}} \sum_{\ell \geq 0} e^{-c 4^\ell} \left( \aver{2^\ell B(x,\sqrt{t})} \big|\sqrt{t} \Gamma\big(\mcP^2_t f\big)\big|^2 \, d\mu \right)^\frac{1}{2} \left( \aver{2^\ell B(x,\sqrt{t})} \big|\sqrt{t} \Gamma\big(\mcP^3_t h\big)\big|^2 \, d\mu \right)^\frac{1}{2}.  
\end{align*}
So using Propositions \ref{prop} and \ref{prop:gradient-Holder}, we deduce that
\begin{align*}
\Big| \mcP^1_t \Big(  \Gamma\big(\sqrt{t} \mcP^2_t f , \sqrt{t} \mcP_t^3 h\big) \cdot \big(g(x) - P_t g\big)\Big)(x) \Big|  \lesssim  \calM\big(L^{\frac{\alpha}{2}} f\big)(x) \, \|g\|_{W^{\beta,p}} \|h\|_{W^{\beta,p}} t^{\frac{\alpha + \beta+ \gamma)}{2}- \frac{\nu}{p}}.
\end{align*}
Then the continuity estimate \eqref{eq:eq1-Ws} comes from integrating with respect to time, taking into account the fact that $\alpha+\beta+\gamma > \frac{2\nu}{p}$.

\medskip

Let us then estimate the regularity of  $F:=g \mcR(f,h) - \mcS(f,g,h)$. It is known, see \cite[Section 2.1.1]{CRT},\cite[Section 5.2]{BBR} or \cite[Proposition 9.7]{BCF2}, that
$$ 
\|F\|_{W^{\delta,p}} \lesssim \|F\|_p + \| S_\delta(F) \|_p,
$$
where $S_\delta$ is the Strichartz functional of index $\delta\in(0,1)$:
$$ 
S_\delta(F)(x) := \left(\int_0^1 r^{-2\delta} \left(\aver{B(x,r)} |F(x) - F(y)| \, d\mu(y)\right)^2 \, \frac{dr}{r}\right)^\frac{1}{2}.
$$
Fix $r>0$ and two points $x,y\in M$ with $d(x,y)\leq r$. Then as previously, we write
\begin{align*} 
F(x) - F(y) & = \big(g(x) \mcR(f,h) - \mcS(f,g,h)\big)(x) -  \big(g(y) \mcR(f,h) - \mcS(f,g,h)(y)\big)(y)\\
&  =: U+ V,
\end{align*}
with $U$ defined by the formula
{\small $$
\int_0^{r^2} \left\{\mcP^1_t\left(  \Gamma(\sqrt{t} \mcP^2_t f , \sqrt{t} \mcP_t^3 h) \cdot \big(g(x) - P_t g\big)\right)(x) - \mcP^1_t\Big(  \Gamma(\sqrt{t} \mcP^2_t f , \sqrt{t} \mcP_t^3 h) \cdot \big(g(y) - P_t g\big)\Big)(y)\right\}\,\frac{dt}{t},
$$}
and $V$ by
{\small $$
\int_{r^2}^1 \left\{\mcP^1_t\left(  \Gamma(\sqrt{t} \mcP^2_t f , \sqrt{t} \mcP_t^3 h) \cdot \big(g(x) - P_t g\big)\right)(x) - \mcP^1_t\Big(  \Gamma(\sqrt{t} \mcP^2_t f , \sqrt{t} \mcP_t^3 h) \cdot \big(g(y) - P_t g\big)\Big)(y)\right\}\,\frac{dt}{t}.
$$}
By repeating previous arguments, we easily bound $U$ as follows
\begin{align*}
U & \lesssim \left(\int_0^{r^2} t^{\frac{\alpha+\beta+\gamma}{2}-\frac{\nu}{p}} \left[\calM\big(Q_t^{(\epsilon)}L^{\frac{\alpha}{2}}f\big)(x)+ \calM\big(Q_t^{(\epsilon)} L^{\frac{\alpha}{2}}f\big)(y)\right] \, \frac{dt}{t}\right) \|g\|_{W^{\beta,p}}\|h\|_{W^{\gamma,p}},
\end{align*}
for some $\epsilon>0$ satisfying $\alpha+\epsilon<1-\theta$. For the second part, we use the estimate
$$ 
|V| \leq A+B,
$$
with $A$ equal to
{\small 
$$
\left| \int_{r^2}^1 \Big\{ \mcP^1_t\Big( \Gamma\big(\sqrt{t} \mcP^2_t f , \sqrt{t} \mcP_t^3 h\big) \cdot \big(g(x) - P_t g\big)\Big)(x) - \mcP^1_t\Big(\Gamma\big(\sqrt{t} \mcP^2_t f , \sqrt{t} \mcP_t^3 h) \cdot \big(g(x) - P_t g\big)\Big)(y)\Big\}\,\frac{dt}{t}\right|
$$}
and 
$$
B:= \int_{r^2}^1 \big|g(x)-g(y)\big|\cdot \Big| \mcP^1_t\Big(\Gamma\big(\sqrt{t} \mcP^2_t f , \sqrt{t} \mcP_t^3 h\big)\Big)(y)\Big| \, \frac{dt}{t}.
$$
The last quantity is bounded (following the same estimates as previously since $g\in \calC^{\beta-\frac{\nu}{p}}$) by
\begin{align*}
B  \lesssim r^{\beta-\nu/p} \left(\int_{r^2}^1 t^{\frac{1}{2}\big(\alpha+\gamma-\frac{\nu}{p}\big)} \calM[Q_t^{(\epsilon)}L^{\frac{\alpha}{2}}f](y) \, \frac{dt}{t}\right)
\|g\|_{W^{\beta,p}}\|h\|_{W^{\gamma,p}}.
\end{align*}
For the quantity $A$, we combine the previous argument with the H\"older regularity of the heat kernel, Lemma \ref{lem:holderpt}, to get the upper bound
{\small \begin{align*}
A & \lesssim  \left(\int_{r^2}^1 \left(\frac{r}{\sqrt{t}}\right)^{1-\theta}  t^{\frac{1}{2}\big(\alpha+\beta+\gamma\big)-\frac{\nu}{p}}  \calM\big(Q_t^{(\epsilon)}L^{\frac{\alpha}{2}}f\big)(y) \, \frac{dt}{t} \right) \|g\|_{W^{\beta,p}} \|h\|_{W^{\gamma,p}}.
\end{align*}}
The combination of all the previous estimates yields
\begin{align*}
\left|F(x)-F(y)\right| & \leq |U|+A+B \\
& \lesssim  \left[\left(\int_0^{r^2} t^{\delta/2} \left[\calM\big(Q_t^{(\epsilon)}L^{\frac{\alpha}{2}}f\big)(x)+ \calM[Q_t^{(\epsilon)}L^{\frac{\alpha}{2}}f](y)\right] \, \frac{dt}{t}\right)  \right. \\
& \qquad \left. +  \left(\int_{r^2}^1 \left(\frac{r}{\sqrt{t}}\right)^{1-\theta}  t^\frac{\delta}{2}  \calM\big(Q_t^{(\epsilon)}L^{\frac{\alpha}{2}}f\big)(y) \, \frac{dt}{t} \right) \right. \\
& \qquad \left. + \left(\int_{r^2}^1  r^{\beta-\frac{\nu}{p}}t^{\frac{1}{2}\big(\alpha+\gamma-\frac{\nu}{p}\big)} \calM[Q_t^{(\epsilon)}L^{\frac{\alpha}{2}}f](y)  \, \frac{dt}{t} \right) \right] \|g\|_{W^{\beta,p}}\|h\|_{W^{\gamma,p}}.
\end{align*}
This estimate holds uniformly for every $y\in B(x,r)$ and so can be averaged on this ball. We then conclude by Hardy's inequality (with $\delta>0$, $\delta<1-\theta$ and $\beta>\frac{\nu}{p}$) that
$$ 
S_\delta(F) \lesssim \left( \int_0^1 \left|\calM \calM\big(Q_t^{(\epsilon)}L^{\frac{\alpha}{2}}f\big)\right|^2 \, \frac{dt}{t} \right)^\frac{1}{2}\|g\|_{W^{\beta,p}}\|h\|_{W^{\gamma,p}}.
$$
Using Fefferman-Stein's inequality and the $L^p$-boundedness of the vertical square function (see Lemma \ref{lem:ortho}), we then deduce that
$$ 
\|F\|_{W^{\delta,p}} \lesssim \|F\|_p + \| S_\delta(F)\|_p \lesssim \|f\|_{W^{\alpha,p}} \|g\|_{W^{\beta,p}}\|h\|_{W^{\gamma,p}},
$$
which concludes the proof of the continuity estimate \eqref{EqEstimateD1-Ws}.

\medskip

\noindent {\bf Step 2 -- proof of \eqref{EqEstimateD2-Ws}}. We need to prove that we have
\begin{equation}
\label{EqSecondBitCommutator-Ws}
(*):=\Big\| L^\frac{\delta}{2} \Big(\mcR\big(\mcA(f,g),h\big) - \mcS(f,g,h)\Big)\Big\|_p \lesssim \|f\|_{W^{\alpha,p}} \|g\|_{W^{\beta,p}}\|h\|_{W^{\gamma,p}},
\end{equation}
where $\mcR\big(\mcA(f,g),h\big) - \mcS(f,g,h)$ is equal to 
\begin{equation}
\label{eq:eq4-Ws}
\int_0^1 \mcP^1_t \Gamma \left( \sqrt{t}\left\{\int_0^1 \mcP^2_t  \mcQ^1_s\left(\mcQ^2_s f \cdot \mcP_s^3 g\right) \, \frac{ds}{s} - P_t g \cdot\mcP^2_t f\right\} \,,\,\sqrt{t} \mcP_t^4 h\right) \, \frac{dt}{t}. 
\end{equation}
Using Lemma \ref{lem:ortho} with the $L^1$-$L^\infty$ off-diagonal estimates of $\mcP_t^1$, we deduce that quantity $(*)$ is bounded above by a constant multiple of 
$$ 
\sum_{\ell \geq 0} 2^{-\ell \delta} \left\| \left(\int_0^1 t^{-\delta} \left( \aver{B(x_0, 2^\ell \sqrt{t})} |A_t(f,g)|^2 \, d\mu \right) \left(\aver{B(x_0, 2^\ell \sqrt{t})} |\sqrt{t} \Gamma \mcP_t^3 h|^2 \, d\mu \right) \, \frac{dt}{t} \right)^\frac{1}{2} \right\|_{p},
$$
with
$$ 
A_{t}(f,g) := \sqrt{t} \Gamma\left(\int_0^1 \mcP^2_t  \mcQ^1_s\left(\mcQ^2_s f \cdot \mcP_s^3 g \right) \, \frac{ds}{s} - P_t g \mcP^2_t f\right).
$$
Following the reasoning in the previous proof of Proposition \ref{prop:commutator-Cs}, by combining with Proposition \ref{prop} we can obtain that for almost every $x_0$
\begin{align*}
 \left( \aver{B(x_0,2^\ell \sqrt{t})} |A_t(f,g)|^2 \, d\mu \right)^\frac{1}{2} & \lesssim 2^{\ell \big(\beta-\frac{\nu}{p}\big)} t^{\frac{1}{2}\big(\beta-\frac{\nu}{p}\big)}  \calM[Q_t^{(\epsilon)} (tL)^{\frac{\alpha}{2}} f](x_0)  \|g\|_{W^{\beta,p}} \\
  & + t^\frac{1-\theta}{2}\calM\big(Q_1^{(\epsilon)} L^{\frac{\alpha}{2}} f]\big)(x_0)  \|g\|_{W^{\beta,p}}.
\end{align*}
Hence, since $2\delta>\beta-\frac{\nu}{p}$ we obtain
\begin{align*}
(*) \lesssim \|g\|_{W^{\beta,p}} \|h\|_{W^{\gamma,p}} \left\{\left\| \left(\int_0^1   \calM[Q_t^{(\epsilon)} L^{\frac{\alpha}{2}}f]^2  \, \frac{dt}{t} \right)^\frac{1}{2} \right\|_{p} + \|f\|_p \right\}.
 \end{align*}
We then conclude to \eqref{EqSecondBitCommutator-Ws} by the Fefferman-Stein inequality and Lemma \ref{lem:ortho}.
\end{Dem}

\bigskip

\subsubsection{Composition estimates}

The above continuity estimates are the main estimates used in the main part of this work to prove the paralinearisation and composition estimates for paraproduct. We state these results here in H\"older and Sobolev spaces under our relaxed assumptions and leave the reader the task of checking that the proofs of section \ref{SubsectionParalinearization} are easily adapted.

\medskip

\begin{theorem}  \label{thm:paralinearization-Cs}    {\sf
Fix an integer $b\geq2$ and a nonlinearity $\textrm{F}\in C_b^{3}$.

\begin{enumerate}
   \item Let $\alpha\in(0,1-\theta)$ be given. For every $f\in \calC^\alpha$, we have $\textrm{F}(f)\in\calC^\alpha$ and
$$ 
R_F(f) := \textrm{F}(f) - \Pi^{(b)}_{\textrm{F}'(f)}(f) \in \calC^{2\alpha}.
$$
More precisely,
$$ 
\Big\|\textrm{F}(f) - \Pi^{(b)}_{\textrm{F}'(f)}(f)\Big\|_{\calC^{2\alpha}} \lesssim \|\textrm{F}\|_{C^3_b} \Big(1+\|f\|_{\calC^\alpha}^2\Big).
$$
If $F\in C^4_b$, then the remainder term $R_F(f)$ is Lipschitz with respect to $f$, in so far as we have
$$ 
\big\|R_\textrm{F}(f)-R_\textrm{F}(g)\big\|_{\calC^{2\alpha}} \lesssim \|\textrm{F}\|_{C^4_b} \big(1+\|f\|_{\calC^\alpha} + \|g\|_{\calC^\alpha} \big)^2\,\|f-g\|_{\calC^\alpha}.
$$ \vspace{0.2cm}
    
   \item Fix $p\in(1,\infty)$. For every $\alpha \in \big(\frac{\nu}{p},1-\theta\big)$ and every $f\in W^{\alpha,p}$, we have $\textrm{F}(f)\in W^{\alpha,p}$ and
$$ 
R_F(f) := \textrm{F}(f) - \Pi^{(b)}_{\textrm{F}'(f)}(f) \in W^{2\alpha-\frac{\nu}{p}\, ,\,p}.
$$
More precisely
$$ 
\Big\|\textrm{F}(f) - \Pi^{(b)}_{\textrm{F}'(f)}(f)\Big\|_{W^{2\alpha-\frac{\nu}{p}\,,\,p}} \lesssim \|\textrm{F}\|_{C^3_b} \Big(1+\|f\|_{W^{\alpha,p}}^2\Big).
$$
If $F\in C^4_b$, then the remainder term $R_F(f)$ is Lipschitz with respect to $f$.
\end{enumerate}   }
\end{theorem}

\medskip 

Let us now examine the composition of two paraproducts. Note that for $u\in \calC^\alpha$ and $v\in \calC^\beta$, with $\alpha\in(0,1)$, $\beta\in(0,\alpha]$, we have $uv\in\calC^\beta$.

\medskip

About the composition of paraproducts, Theorem \ref{thm:composition}  in H\"older spaces still holds since it only relies on the Gaussian estimate \eqref{UE}; its Sobolev counterpart also holds.

\begin{theorem}  \label{thm:composition-Cs}    {\sf
\begin{enumerate}
   \item Fix an integer $b\geq2$, $\alpha\in(0,1)$, $\beta\in(0,\alpha]$ and consider $u\in \calC^\alpha$ and $v\in \calC^\beta$. Then for every $f\in \calC^\alpha$, we have
$$
\Pi^{(b)}_{u} \Big( \Pi^{(b)}_{v}(f) \Big) - \Pi^{(b)}_{uv}(f) \in \calC^{\alpha+\beta}
$$
with
$$ 
\left\| \Pi^{(b)}_{u}\Big( \Pi^{(b)}_{v}(f) \Big) - \Pi^{(b)}_{uv}(f) \right\|_{\calC^{\alpha+\beta}} \lesssim \|f\|_{\calC^\alpha} \, \|u\|_{\calC^\alpha} \|v\|_{\calC^\beta}.
$$ \vspace{0.2cm}
   
   \item Fix an integer $b\geq2$ and $p\in(1,\infty)$. For $\alpha\in(0,1)$ and $\beta \in \big(\frac{\nu}{p},\alpha\big]$, consider $u\in W^{\alpha,p}$ and $v\in W^{\beta,p}$. Then for every $f\in W^{\alpha,p}$, we have
$$
\Pi^{(b)}_{u} \Big( \Pi^{(b)}_{v}(f) \Big) - \Pi^{(b)}_{uv}(f) \in W^{\alpha+\beta-2\frac{\nu}{p}\,,\,p}
$$
with
$$ 
\left\| \Pi^{(b)}_{u}\Big( \Pi^{(b)}_{v}(f) \Big) - \Pi^{(b)}_{uv}(f) \right\|_{W^{\alpha+\beta-2\frac{\nu}{p}\,,\,p}} \lesssim \|f\|_{W^{\alpha,p}} \, \|u\|_{W^{\alpha,p}} \|v\|_{W^{\beta,p}}.
$$
\end{enumerate}    }
\end{theorem}

\medskip

\subsubsection{Schauder estimates}

Proposition \ref{prop:schauder1} gives an elementary proof in our setting of a Schauder-type estimate about the regularizing character of the convolution operation with the operators $P^{(b)}_s$. The same properties hold in our minimal setting since hey only rely on Gaussian property \eqref{UE} and the semigroup structure, together with a Sobolev version which we state here without proof as it can be proved along the lines of proof of proposition \ref{prop:schauder1}. (Another approach can be also obtained by interpolating between the trivial case $\epsilon=1$ and the limit case $\epsilon=0$. The latter case $\epsilon=0$, corresponds exactly to the so-called $L^p$ maximal regularity which  has been the topic of a huge literature, see for example \cite{HP} where the Gaussian upper estimates \eqref{UE} are used.)

\medskip

\begin{proposition}  \label{prop:schauder-Cs}    {\sf
\begin{enumerate}
   \item Consider $\beta\in{\mathbb R}$ and $\epsilon\in(0,1)$. For every $T>0$ and $v\in C_T\calC^\beta$, the function $V(t) := \int_0^t P^{(b)}_{t-s} v(s)\, ds$ belongs to $C_T\calC^{\beta+2-2\epsilon}$ with
$$ 
\big\|V(t)\big\|_{\calC^{\beta+2-2\epsilon}} \lesssim T^\epsilon \sup_{s\in[0,t]}\, \big\|v(s)\big\|_{\calC^\beta}
$$
and 
$$
\big\|V\big\|_{C_T^{(\beta+2-2\epsilon)/2}L^\infty} \lesssim T^\epsilon\, \|v\|_{C_T \calC^\beta}.
$$
So 
$$ 
\big\|Jf \big\|_{\LL^\alpha_T} \lesssim T^\epsilon \,\|f\|_{C_T \calC^{\alpha-2+2\epsilon}}.
$$  \vspace{0.2cm}

   \item Consider $\beta\in{\mathbb R}$, $p\in(1,\infty)$ and  $\epsilon\in(0,1)$. For every $T>0$ and $v\in C_T W^{\beta,p}$, the function $V(t) := \int_0^t P^{(b)}_{t-s} v(s)\, ds$ belongs to $C_TW^{\beta+2-2\epsilon,p}$ with
$$ 
\big\|V(t)\big\|_{W^{\beta+2-2\epsilon,p}} \lesssim T^\epsilon \sup_{s\in[0,t]}\, \big\|v(s)\big\|_{W^{\beta,p}}
$$
and 
$$
\big\|V\big\|_{C_T^{\frac{1}{2}\,(\beta+2-2\epsilon)}L^p} \lesssim T^\epsilon\, \|v\|_{C_T W^{\beta,p}}.
$$ 
\end{enumerate}   }
\end{proposition}

\bigskip

\subsection[\hspace{-1cm} {\sf Resolution of PAM in such a 2-dimensional setting}]{Resolution of PAM equation in such a 2-dimensional setting}
\label{SubsectionAppPAM}

Building on the estimates proved in this Appendix, it is elementary to introduce and study paracontrolled distributions in H\"older and Sobolev spaces along the lines of Subsections \ref{SubsectionParaControlledDistributions} and \ref{SubsectionParacontrolledSchauder}, in the present extended setting. Its application to the parabolic Anderson model equation (PAM) is also almost straightforward as we only need to check that the renormalization procedure explained in details in subsection \ref{SubsectionRenormalization} under the (Lip) assumption can be run here as well. This is indeed the case if the exponent $q_0$ in the gradient assumption $(\textrm{G}_{q_0})$ is large enough, as this assumption yields some "$L^{q_1}$-Gaussian" estimates for every $q_1<q_0$.

\ssk

Let us compute, as an example, an integral of type 
$$
I_{s,t} := \int \Gamma_x p_t(x,y) \Gamma_x p_s(x,z) \, d\mu(x),
$$
where $p_t$ is the heat kernel of $e^{-tL}$, and $s<t$. By \eqref{gradlp} with the local Ahlfors regularity, there exists a positive constant $c$ such that we have
\begin{align*}
I_{s,t}& \lesssim  t^{-\frac{\nu}{2q'_1}}s^{-\frac{\nu}{2q'_1}} \left( \int e^{c\frac{d(x,y)^2}{t}} \, e^{c\frac{d(x,z)^2}{s}} \, \mu(dx) \right)^{1-\frac{2}{q_1}} \\
         & \lesssim  t^{-\frac{\nu}{2q_1}}s^{-\frac{\nu}{2q_1}}  \left( \int \G_t(x,y) \, \G_s(x,z)  \, \mu(dx) \right)^{1-\frac{2}{q_1}} \\
         & \lesssim t^{-\frac{\nu}{2q_1}}s^{-\frac{\nu}{2q_1}} \G_{t+s}(y,z)^{\frac{\nu}{2}\big(1-\frac{2}{q_1}\big)}, 
\end{align*}
where we used Lemma \ref{lem:gauss}. So with respect to Subsection \ref{SubsectionRenormalization}, where \eqref{Lipschitz} was assumed and where $I_{s,t}$ would be estimated by $\G_{s+t}(y,z)$, we now have the estimate
$$
I_{s,t} \lesssim  \left(\frac{(t+s)^2}{ts}\right)^\frac{\nu}{2q_1} \G_{t+s}(y,z),
$$
involving an extra factor $\left(\frac{(t+s)^2}{ts}\right)^\frac{\nu}{2q_1}$. Since all the conditions on the exponents were open conditions in Subsection \ref{SubsectionRenormalization}, we may allow a small loss if it is small enough. As a consequence, we deduce that if $q_1$ can be chosen large enough then we may adapt and repeat the renormalization procedure of the white noise in H\"older and Sobolev spaces. The latter condition on $q_1$ is equivalent to taking $q_0$ big enough.

\bigskip

We summarize this result under the following form, which gives an analogue of theorem \ref{thm:pam}. 

\medskip 

\begin{theorem}   \label{thm:pam-bis}   {\sf 
Assume the local Ahlfors regularity of dimension $2$, as well as $(P_2)$, $(\textrm{G}_{q_0})$ and $(D\textrm{G}_{q,\theta})$ for $q_0$ large enough and $\theta$ small enough. Fix $p>2$ a large enough exponent.

Let $\xi$ stand for a time-independent weighted noise in space, and set $\xi^\epsilon := P_\epsilon\xi$, and $X^\epsilon(t) = \int_0^t P_{t-s}\big(\xi^\epsilon\big)\,ds$.  \vspace{0.1cm}
\begin{enumerate}
   \item The pair $\big(\xi^\epsilon , X^\epsilon\big)$ converges in probability in some space (in the H\"older scaling $(\calC^s)_s$ or Sobolev scaling $(W^{s,p})$) to some extended noise $(\zeta,X)$, with $\zeta = \xi$, and $\Pi(X,\zeta)$ well-defined in the above sense.   \vspace{0.1cm}

   \item Furthermore, if $u^\epsilon$ stands for the solution of the renormalized equation 
   \begin{equation} 
   \label{eq:renor-bis}
   \partial_t u^\epsilon + L u^\epsilon = \textrm{F}\big(u^\epsilon\big)\,\xi^\epsilon - c^\epsilon\,\textrm{F}'\big(u^\epsilon\big)\,\textrm{F}(u^\epsilon) ,\qquad u^\epsilon(0)=u_0  
   \end{equation}
where $c^\epsilon(\cdot) := \E\Big[\Pi\big(L^{-1}\xi^\epsilon,\xi^\epsilon\big)(\cdot)\Big]$ is a deterministic real-valued function on $M$, then $u^\epsilon$ converges in probability to the solution $u$ of $(gPAM)$ associated with $(\zeta,X)$, in some space whose definition depends on whether or not $\textrm{F}$ is linear.
\end{enumerate}   }
\end{theorem}



\begin{thebibliography}{AAA}

\bibitem{Aus}
P. Auscher,
\newblock Regularity theorems and heat kernel for elliptic operators, 
\newblock {\it J. London Math. Soc.},  \textbf{54} (1996), no. 2, 284--296.

\bibitem{ACDH} P. Auscher, T. Coulhon, X. T. Duong, and S. Hofmann, 
 Riesz transform on manifolds and heat kernel regularity,
 \textit{Ann. Sci. Ecole Norm. Sup.}, \textbf{37} (2004), no. 4, 911--957.

\bibitem{ADM} 
D. Albrecht, X. T. Duong and A. McIntosh,
\newblock Operator theory and harmonic analysis, in {\em Instructional Workshop on Analysis and Geometry, Part III (Canberra, 1995)}, 
\newblock { Proc. Centre Math. Appl. Austral. Nat. Univ. 34} (1996), 77--136. 

\bibitem{AMT}
P. Auscher, A. McIntosh and P. Tchamitchian,
\newblock Noyau de la chaleur d'op\'rateurs elliptiques complexes,
\newblock {\em Math. Research Lett.} \textbf{1} (1994), 37--45.

\bibitem{BBR} N. Badr, F. Bernicot, and E. Russ, 
Algebra properties for Sobolev spaces-applications to semilinear PDEs on manifolds, 
\newblock{\em J. Anal. Math.}, \textbf{118}, no.2 (2012),  509--544.

\bibitem{BM2Course}
I. Bailleul,
\newblock A flow-based approach to rough differential equations,
\newblock https://perso.univ-rennes1.fr/ismael.bailleul/files/M2Course.pdf

\bibitem{BBF16} I. Bailleul, F. Bernicot and D. Frey,
Higher order paracontrolled calculus, 3d-PAM and multiplicative Burgers equations,
\newblock arXiv:1506.08773, 2015.

\bibitem{BCD}
H. Bahouri, R. Danchin and J.Y. Chemin,
\newblock Fourier analysis and nonlinear partial differential equations,
\newblock {Grundlehren der Mathematischen Wissenschaften 343}, (2011).

\bibitem{BG}
F. Baudoin and N. Garofalo,
\newblock Curvature-dimension inequalities and Ricci lower bounds for sub-Riemannian manifolds with transverse symmetries,
\newblock {\it J. Eur. Math. Soc.} (2015).

\bibitem{B-T1}
F. Bernicot, 
\newblock A T(1)-Theorem in relation to a semigroup of operators and applications to new paraproducts,
\newblock {\it Trans. Amer. Math. Soc.} \textbf{364} (2012), 257--294.

\bibitem{BCF1}
F. Bernicot, T. Coulhon  and D. Frey, Gaussian heat kernel bounds through elliptic Moser iteration,
\newblock {\it J. Math. Pures Appl.} (2016), arXiv:1407.3906.

\bibitem{BCF2}
F. Bernicot, T. Coulhon  and D. Frey, Sobolev algebra through Heat semigroup,
\newblock {\it submitted}, arXiv:1505.01442.

\bibitem{BF2}
F. Bernicot and D. Frey, 
\newblock Riesz transforms through reverse H\"older and Poincar\'e inequalities,
{\it submitted}, arXiv:1503.02508.

\bibitem{BS}
F. Bernicot and Y. Sire, 
\newblock Propagation of low regularity for solutions of nonlinear PDEs on a Riemannian manifold with a sub-Laplacian structure,
\newblock {\it Ann. I. H. Poincar\'e - AN} \textbf{30} (2013), 935--958.

\bibitem{Bony}
J.M. Bony, 
\newblock Calcul symbolique et propagation des singulari\'es pour les \'equations aux d\'eriv\'ees partielles non lin\'eaires, 
\newblock {\it Ann. Sci. Eco. Norm. Sup.} 1\textbf{14} (1981), 209--246.

\bibitem{BCS} 
S. Boutayeb, T. Coulhon and A. Sikora, 
\newblock A new approach to pointwise  heat kernel upper bounds on doubling metric measure spaces, 
\newblock  {\it Adv. in Math.} \textbf{270} (2015), 302--374.

\bibitem{BDY}
H-Q. Bui, X. T. Duong and L. Yan, 
\newblock Calder\'on reproducing formulas and new Besov spaces associated with operators,
\newblock {\it Adv. in Math.} \textbf{229} (2012), 2449--2502.

\bibitem{CC15}
R. Catellier and K. Chouk,
\newblock Paracontrolled distributions and the 3-dimensional stochastic quantization equation,
\newblock To appear in {\em Ann. Probab.}, 2016.

\bibitem{CRT} T. Coulhon, E. Russ, and V. Tardivel-Nachef, Sobolev algebras on Lie groups and Riemannian manifolds,  {\em Amer. J. of Math.}, \textbf{123} (2001), 283--342.

\bibitem{CS} 
T. ~Coulhon and A. ~Sikora,
\newblock Gaussian heat kernel bounds via Phragm\'en-Lindel\"of theorem,  
\newblock {\it Proc. London Math. Soc.} \textbf{96} (2008), 507--544.

\bibitem{CDMY}
M. Cowling, I. Doust, A. McIntosh and A. Yagi,
\newblock Banach space operators with a bounded $H^\infty$ functional calculus, 
\newblock {\it J. Austral. Math. Soc. Ser. A} \textbf{60} (1996), 51--89.

\bibitem{DeG}
E.  De Giorgi, 
Sulla differenziabilita de analiticita delle estremali degli integrali multipli regolari, 
\newblock {\it Mem. Accad. Sci. Torino Cl. Sci. Fis. Mat. Nat.}, \textbf{3}  (1957), no. 3, 25--43.

\bibitem{DR}
X.-T. Duong and D. W. Robinson,
\newblock  Semigroup kernels, Poisson bounds, and holomorphic functional calculus,
\newblock {\it  J. Funct. Anal.} \textbf{142}(1) (1996), 89--128. 

\bibitem{FH}
P. Friz and M. Hairer,
\newblock  A course on rough paths (with an introduction to regularity stuctures),
\newblock Springer, 2014. 


\bibitem{FV10}
P. Friz and N. Victoir, 
\newblock {\it Multidimensional stochastic processes as rough paths; Theory and Applications.} 
\newblock Cambridge University Press, 2010.

\bibitem{FOT} 
M. ~Fukushima, Y. ~Oshima and M. ~Takeda, 
\newblock {\it Dirichlet forms and symmetric Markov processes.} 
\newblock De Gruyter Studies in Mathematics 19, {Walter de Gruyter}, 
\newblock Berlin, 1994.

\bibitem{GS}
I. Gallagher and Y. Sire,
\newblock Besov algebras on Lie groups of polynomial growth,
\newblock {\it Studia Math.} \textbf{212} (2012), no. 2, 119--139. 
 
\bibitem{Gr1} 
A.~Grigor'yan, 
\newblock Gaussian upper bounds for the heat kernel on arbitrary manifolds,
\newblock {\it J. Diff. Geom.} \textbf{45} (1997), 33--52.

\bibitem{Gr2} 
A.~Grigor'yan, 
\newblock Analytic and geometric background of recurrence and non-explosion of the Brownian motion on Riemannian manifolds,
\newblock {\it Bull. Amer. Math. Soc.} \textbf{36} (1999), 135--249.

\bibitem{GL} 
A. Grigor'yan and L. Liu, 
\newblock Heat kernel and Lipschitz-Besov spaces, 
\newblock {\it Forum Math.} (2015).

\bibitem{GControlled}
M. Gubinelli, 
\newblock Controlling rough paths,
\newblock {\it J. Funct. Anal.} \textbf{2016}, no. 1, (2004), 86--140. 

\bibitem{GIP}
M. Gubinelli, P. Imkeller and N. Perkowski,
\newblock Paracontrolled distributions and singular PDEs,
\newblock Preprint 2012, arXiv:1210.2684.

\bibitem{G}
Y. Guivar'ch,
\newblock Croissance polynomiale et p\'eriodes des fonctions harmoniques,
\newblock {\it Bull. Soc. Math. France} \textbf{101} (1973), 333--379.

\bibitem{GSC}
P. Gyrya and L. Saloff-Coste,
\newblock {\it Neumann and Dirichlet heat kernels in inner uniform domains},
\newblock Ast{\'e}risque \textbf{33} (2011), Soc. Math. France.

\bibitem{H}
M. Hairer,
\newblock A theory of regularity structures,
\newblock {\it Invent. Math.}, 198 (2014), no. 2, 269--504.

\bibitem{HCPAM}
M. Hairer,
\newblock Rough Stochastic PDEs,
\newblock {\it Commun. Pure Appl. Math.}, 64 (2011), no. 11, 1547--1585.

\bibitem{HKPZ}
M. Hairer,
\newblock Solving the KPZ equation ,
\newblock {\it Annals of Maths}, 178 (2013), no. 2, 559--664.

\bibitem{HL}
M. Hairer and C. Labb\'e,
\newblock A simple construction of the continuum parabolic Anderson model on ${\mathbb R}^2$,
\newblock {\em Elect. J. Probab.}, 2016.

\bibitem{HLR3}
M. Hairer and C. Labb\'e,
\newblock Multiplicative stochastic heat equations on the whole space,
\newblock Preprint, arXiv:arXiv:1504.07162.

\bibitem{HW}
M. Hairer and H. Weber,
\newblock Rough Burgers-like equations with multiplicative noise,
\newblock {\it Probab. Theory Rel. Fields}, 155 (2013), no. 1, 71--126.

\bibitem{HSc1} 
W. Hebisch and L. Saloff-Coste, 
\newblock Gaussian estimates for Markov chains and random walks on groups, 
\newblock {\it Ann. Probab.} \textbf{21} (1993), no. 2, 673--709.

\bibitem{HS}
W. Hebisch and L. Saloff-Coste,
\newblock On the relation between elliptic and parabolic Harnack inequalities,
\newblock {\it Ann. Inst. Fourier} \textbf{51} (2001), no. 5, 1437--1481.

\bibitem{HP}
M. Hieber and J. Pr\"uss,
\newblock Heat kernels and maximal $L^p$-$L^q$ estimates for parabolic evolution equations,
\newblock {\it Comm. in Partial Differential Equation} \textbf{22} (1997), 1647--1669.

\bibitem{Is} 
S. Ishiwata, 
\newblock A Berry-Esseen type theorem on nilpotent covering graphs, 
\newblock {\it Canad. J. Math.} {\bf 56} (2004), no. 5, 963--982.

\bibitem{Is2} 
S. Ishiwata, 
\newblock Gradient estimate of the heat kernel on modified graphs, 
\newblock {\it Pot. Analysis} {\bf 27} (2007), no. 4, 335--351.

\bibitem{KW}
P. C. Kunstmann and L. Weis,
\newblock{\it  Maximal $L_p$ regularity for parabolic problems, Fourier multiplier theorems and $H^{\infty}$-functional calculus}, 
\newblock {Lect. Notes in Math. 1855}. Springer-Verlag (2004).

\bibitem{Li}
H-Q Li,
\newblock La transformation de Riesz sur les vari\'et\'es coniques,
\newblock {\it J. Funct. Anal.} \textbf{168} (1999), 145--238.

\bibitem{LY}
P. Li and S. T. Yau, 
\newblock On the parabolic kernel of the Schr\"odinger operator,
\newblock {\it Acta Math.} \textbf{156} (1986), 153--201.

\bibitem{LYY}
L. Liu, D. Yang and W. Yuan, 
\newblock Besov-type and Triebel-Lizorkin-type spaces associated with Heat kernels,
\newblock Preprint 2013, arXiv:1309.1366.

\bibitem{Lyons98}
T. Lyons,
\newblock Differential equations driven by rough signals,
\newblock{\em Rev.Mat. Iberoamericana} {\bf 14}(2), 1998.

\bibitem{LyonsStFlour}
M. Caruana, T. L\'evy and T. Lyons,
\newblock 
\newblock{\em Lect. Notes Math.} 1908, 2007.

\bibitem{Meyer} 
Y. Meyer, 
\newblock {\it Wavelets and operators}, 
\newblock Cambridge Studies in advanced mathematics (1992).

\bibitem{NSW}
A. Nagel, E. M. Stein and S. Wainger,
\newblock Balls and metrics defined by vector fields I: Basic properties,
\newblock {\it Acta Math.} \textbf{155} (1985), 103--147.

\bibitem{Qian}    
B. Qian,
\newblock Hamilton type gradient estimate for the sub-elliptic operators,
\newblock {\it Potential Anal.} (2015).

\bibitem{Sa} 
L. ~Saloff-Coste, 
\newblock Analyse sur les groupes de Lie \`a croissance polynomiale, 
\newblock {\it Ark. Mat.} \textbf{28} (1990), 315--331.

\bibitem{S} 
L. ~Saloff-Coste,
\newblock A note on Poincar\'e, Sobolev, and Harnack inequalities, 
\newblock {\it Int. Math. Res. Not.} (1992), no. 2,  27--38.

\bibitem{topics} 
E. M. Stein, 
\newblock {\it Topics in harmonic analysis related to the Littlewood-Paley theory.} 
\newblock Princeton University Press, 1970.

\bibitem{ST1} 
K.T. Sturm, 
\newblock Analysis on local Dirichlet spaces I. Recurrence, conservativeness  and $L^p$-Liouville property, 
\newblock {\it  J. Reine Angew. Math.} \textbf{456} (1994), 173--196.

\bibitem{ST2} 
K.T. Sturm, 
\newblock Analysis on local Dirichlet spaces II. Upper Gaussian estimates for the fundamental solutions of parabolic equations, 
\newblock {\it Osaka J. Math.} \textbf{32} (1995), no. 2, 275--312.

\bibitem{Var0} 
N. Th. Varopoulos,
\newblock Analysis on Lie groups,
\newblock {\it J. Funct. Anal.} \textbf{76} (1988), 346--410.

\bibitem{Var1} 
N. Th. Varopoulos,
\newblock Small time Gaussian estimates of heat diffusion kernels. Part 1. The semigroup technique,
\newblock {\it Bull. Sci. Math.} Ser. 2 \textbf{113} (1989), 253--277.

\bibitem{VSC} 
N. Th. Varopoulos, L. Saloff-Coste and T. Coulhon,
\newblock {\it Analysis and geometry on groups.}
\newblock Cambridge University Press, 1992.

\bibitem{WY}
F.-Y. Wang and L. Yan,
\newblock Gradient estimate on convex domains and applications,
\newblock {\it Proc. Amer. Math. Soc.} \textbf{141} (2013), 1067--1081.

\bibitem{ZZ1} 
R. Zhu and X. Zhu,
\newblock Three dimensional {N}avier-{S}tockes equation driven by space-time white noise,
\newblock arXiv:1406.0047, (2014).

\bibitem{ZZ2} 
R. Zhu and X. Zhu,
\newblock Approximating three-dimensional {N}avier-{S}tokes equations driven by space-time white noise,
\newblock arXiv:1409.4864 , (2014).

\end{thebibliography}
\end{document}